\documentclass[a4paper, 11pt]{article}
\usepackage{etoolbox}
\usepackage{graphicx}
\usepackage[T1]{fontenc}
\usepackage[utf8]{inputenc}
\usepackage{textcomp}
\usepackage[]{amsmath, amssymb}
\usepackage{fancyvrb}
\usepackage{caption}
\usepackage{float}
\usepackage{url}
\usepackage{pdfpages}
\usepackage{listings}
\usepackage{amsmath, amssymb}
\usepackage{amsthm}
\usepackage{grffile}
\usepackage{multirow}
\usepackage{pgfplots}
\usepackage{silence}
\usepackage{todonotes}
\usepackage{subcaption}
\usepackage{gensymb}
\usepackage{epstopdf}

\DefineVerbatimEnvironment{code}{Verbatim}{fontsize=\small, xleftmargin=5mm}
\newcommand{\Z}{{\mathbb Z}}

\newcommand{\abs}[1]{\left|#1\right|}

\newcommand{\dr}{\mathrm{d}}

\newcommand{\dd}[2]{\frac{\dr#1}{\dr#2}}
\newcommand{\pdd}[2]{\frac{\partial #1}{\partial #2}}

\newcommand{\pddn}[3]{\frac{\partial^{#1} #2}{\partial #3^{#1}}}

\newcommand{\bip}[3]{\left( #1,#2 \right)_{#3}}
\newcommand{\ip}[3]{\left( #1, #2 \right)_{#3}}
\newcommand{\norm}[1]{\left\lVert #1 \right\rVert}
\newcommand{\seminorm}[1]{{\left\vert\kern-0.25ex\left\vert\kern-0.25ex\left\vert #1
    \right\vert\kern-0.25ex\right\vert\kern-0.25ex\right\vert}}

\newcommand{\dnorm}[1]{ \norm{#1}_{\Omega,h}^2 }

\newcommand{\myint}[4]{ \int \limits_{#1}^{#2} \! #3 \, \mathrm{d} #4}

\newcommand{\physdom}{\Omega'}
\newcommand{\refdom}{\Omega}
\newcommand{\physboundary}{\Gamma'}
\newcommand{\refboundary}{\Gamma}
\newcommand{\normphys}[1]{\norm{#1}_{\physdom} }
\newcommand{\normphysu}[1]{\norm{#1}_{\physdom_u} }
\newcommand{\normphysv}[1]{\norm{#1}_{\physdom_v} }
\newcommand{\dnormphys}[1]{\norm{#1}_{\physdom, h} }
\newcommand{\dnormphysu}[1]{\norm{#1}_{\physdom_u, h} }
\newcommand{\dnormphysv}[1]{\norm{#1}_{\physdom_v, h} }

\newcommand{\volint}[2]{ \ip{#1}{#2}{\Omega} }
\newcommand{\volintphys}[2]{ \ip{#1}{#2}{\physdom} }
\newcommand{\surfint}[2]{ \bip{#1}{#2}{\Gamma} }
\newcommand{\surfintphys}[2]{ \bip{#1}{#2}{\physboundary} }

\newcommand{\surfintindexphys}[3]{ \bip{#1}{#2}{\physboundary_{#3}} }
\newcommand{\dvolint}[2]{ \ip{#1}{#2}{\Omega,h} }
\newcommand{\dvolintphys}[2]{ \ip{#1}{#2}{\physdom,h} }
\newcommand{\dsurfint}[2]{ \bip{#1}{#2}{\Gamma,h}}
\newcommand{\dsurfintphys}[2]{ \bip{#1}{#2}{\physboundary,h}}

\newcommand{\dsurfintphysindex}[3]{ \bip{#1}{#2}{\physboundary_{#3},h} }
\newcommand{\dsurfintindexsign}[4]{ \bip{#1}{#2}{\Gamma_{#3}^{#4},h} }
\newcommand{\bnorm}[1]{ \norm{#1}_{\Gamma}^2 }

\newcommand{\dbnorm}[1]{ \norm{#1}_{\Gamma,h}^2 }
\newcommand{\dbnormphys}[1]{ \norm{#1}_{\physboundary, h}^2 }

\newcommand{\dbnormphysindex}[2]{ \norm{#1}_{\physboundary_{#2},h}^2 }

\newcommand{\bv}[1]{\mathbf{#1}}

\newcommand{\hphys}{\widetilde{h}}

\newcommand{\db}{\hat{D}}
\newcommand{\ddxi}{\partial_{x_i}}
\newcommand{\ddxii}{\partial_{\xi_i}}
\newcommand{\ddxij}{\partial_{\xi_j}}

\newcommand{\ddi}{\partial_i}
\newcommand{\ddj}{\partial_j}
\newcommand{\ddk}{\partial_k}

\newcommand{\ut}{\dot{u}}

\newcommand{\si}{\sum \limits_{i}}

\newcommand{\sk}{\sum \limits_{k}}

\newcommand{\skjnk}{\sum \limits_{k, j \neq k}}

\newcommand{\bfu}{{\bv{u}}}
\newcommand{\bfv}{{\bv{v}}}
\newcommand{\bfut}{\dot{\bfu}}
\newcommand{\bfvt}{\dot{\bfv}}

\newcommand{\cD}{\mathbb{D}}
\newcommand{\K}{\mathcal{K}}

\theoremstyle{plain}
\newtheorem{theorem}{Theorem}
\newtheorem{lemma}[theorem]{Lemma}

\theoremstyle{definition}

\theoremstyle{remark}

\newlength\figureheight
\newlength\figurewidth

\definecolor{light-gray}{gray}{0.8}

\title{Non-stiff narrow-stencil finite difference approximations of the Laplacian on curvilinear multiblock grids}

\author{Martin Almquist\thanks{Department of Geophysics, Stanford University, Stanford, CA, USA} \and Eric M.\ Dunham$^{*,}$\thanks{Institute for Computational and Mathematical Engineering, Stanford University, Stanford, CA, USA} }

\date{}

\begin{document}
\maketitle

\begin{abstract}
  The Laplacian appears in several partial differential equations used to model wave propagation. Summation-by-parts--simultaneous approximation term (SBP-SAT) finite difference methods are often used for such equations, as they combine computational efficiency with provable stability on curvilinear multiblock grids. However, the existing SBP-SAT discretization of the Laplacian quickly becomes prohibitively stiff as grid skewness increases. The stiffness stems from the SATs that impose inter-block couplings and Dirichlet boundary conditions. We resolve this issue by deriving stable SATs whose stiffness is almost insensitive to grid skewness. The new discretization thus allows for large time steps in explicit time integrators, even on very skewed grids. It also applies to the variable-coefficient generalization of the Laplacian. We demonstrate the efficacy and versatility of the new method by applying it to acoustic wave propagation problems inspired by marine seismic exploration and infrasound monitoring of volcanoes.
\end{abstract}

\section{Introduction}
Partial differential equations (PDEs) that feature the Laplacian include the second order wave equation. For wave-dominated equations, high-order finite difference methods are computationally efficient \cite{KreissOliger72}, assuming that the solution is smooth and the domain not too complex. Finite difference operators with the summation-by-parts (SBP) property \cite{KreissScherer74} lead to energy-stable discretizations when combined with suitable methods for imposing boundary and interface conditions. Weak enforcement via simultaneous approximation terms (SATs) \cite{CarpenterGottlieb94} has proven successful in a wide range of applications \cite{DelReyFernandez2014a,Svard2014}. A key property of the SBP-SAT method is that it lends itself well to curvilinear multiblock grids, which are necessary in most applications.

For second order differential operators in general, narrow-stencil second-derivative SBP operators \cite{MattssonNordstrom04,Mattsson11} provide superior solution accuracy compared to applying a first-derivative operator twice. As a rule of thumb, the global convergence rate is one order higher \cite{MattssonNordstrom04} and the numerical dispersion relation mimics the exact dispersion relation better for marginally resolved modes \cite{Kreiss_wave}. Hence, we only consider narrow-stencil operators in this paper. However, the gain in accuracy comes at the cost of more involved stability analysis.

Key properties of the discrete Laplacian that guarantee energy stability for hyperbolic problems are symmetry and negative semidefiniteness with respect to a discrete inner product (see for example \cite{AlmquistWangWerpers2019}). The current state-of-the-art SBP-SAT discretization of the Laplacian was presented by Virta and Mattsson \cite{VirtaMattsson14}, who were the first to derive a symmetric negative semidefinite discretization for Dirichlet boundary conditions on curvilinear multiblock grids. However, as we demonstrate in Section \ref{sec:spectral_radius}, the spectral radius of the discretization by Virta and Mattsson grows rapidly as grid skewness increases. Hence, it forces explicit time integrators to use very small time steps unless the grid is almost Cartesian. We here improve the method by deriving a new discretization whose spectral radius is much less sensitive to grid skewness. The two methods use the same SBP operators and differ only in the SATs for Dirichlet boundary conditions and inter-block couplings. The SATs for Neumann boundary conditions are identical.

Although we use the acoustic wave equation in two dimensions in the numerical experiments, the new discretization is in no way specific to wave equations or two dimensions. Rather, what we present is a symmetric and negative semidefinite discrete Laplacian that may be used to discretize other PDEs, for example the Schrödinger and heat equations. Furthermore, the new discretization is presented for $d$ spatial dimensions and applies to the generalized Laplace operator $\nabla \cdot b(\vec{x}) \nabla$.

This paper only considers weak enforcement of boundary and interface conditions using SATs, but it is possible to combine SBP operators with strong enforcement by, for example, injection \cite{Duru2014} or ghost points \cite{Sjogreen2012, WangPetersson2018}. In fact, for certain classes of boundary and interface conditions, the ghost point technique yields a significantly smaller spectral radius than SATs \cite{WangPetersson2018}. We stress that the purpose of this paper is not to find the finite difference discretization of the Laplacian with the smallest possible spectral radius for a given set of boundary conditions. Instead, the purpose is to improve the general-purpose SBP-SAT framework by showing how to impose Dirichlet-type boundary conditions and couple grid blocks with less stiffness than with the existing SBP-SAT method.

The rest of this paper is organized as follows. Section \ref{sec:spectral_radius} compares the new method with the current state-of-the-art method by Virta and Mattsson \cite{VirtaMattsson14}. Sections 3 through 8 derive the new method and prove that it is stable. Section \ref{sec:applications} demonstrates the efficacy of the new method by solving two application problems inspired by marine seismic exploration and infrasound monitoring of volcanoes.

\section{Comparison with Virta and Mattsson, 2014.} \label{sec:spectral_radius}
This section presents the main result: the new method, derived in subsequent sections of this paper, is significantly more efficient than the method in Virta and Mattsson (VM) \cite{VirtaMattsson14} when the grid is curved. To compare the two methods, we will consider the second order wave equation with constant coefficients and homogeneous Dirichlet boundary conditions,
\begin{equation} \label{eq:pde_constant_coeff}
\begin{array}{ll}
	\ddot{u} = \nabla^2 u , & \vec{x} \in \physdom,\\
	u = 0, & \vec{x} \in \partial \physdom.
\end{array}
\end{equation}
For the spatial discretization we use diagonal-norm SBP operators \cite{Mattsson11} of interior orders $q=2$, $q=4$ and $q=6$. Near boundaries, the accuracy drops to order $q/2$. As a rule of thumb, one typically observes global convergence rates of order $\min(q, q/2+2)$ for the wave equation \cite{SvardNordstrom06}, but the convergence rate can be both lower and higher in special cases \cite{Wang2016,Wang2018}. By the rule of thumb, we expect convergence rates of orders 2, 4, and 5, respectively, in this paper. The boundary conditions are enforced using either the VM SATs or the new SATs.

Section \ref{subsec:sr} shows that the spectral radius of the discrete Laplacian (including SATs) grows much faster with grid skewness with the VM SATs than with the new SATs. This implies that the VM method requires smaller time steps in explicit time steppers, or more iterations for convergence in implicit time steppers combined with iterative solvers.

Section \ref{subsec:eff} compares the two methods in terms of CPU time required  to satisfy a given error tolerance. For this comparison we use a mildly curved multiblock grid on the unit disk and an explicit time stepper. Although the VM method is slightly more accurate on a given grid, the new method is significantly more efficient because it allows for much larger time steps. In light of Section \ref{subsec:sr}, we anticipate that a more skewed grid would only increase the efficiency gap in favor of the new method.

\subsection{Spectral radius as a function of grid skewness} \label{subsec:sr}

\begin{figure}[h]
    \begin{subfigure}[b]{.5\linewidth}
        \centering
        \includegraphics[width= 1\linewidth]{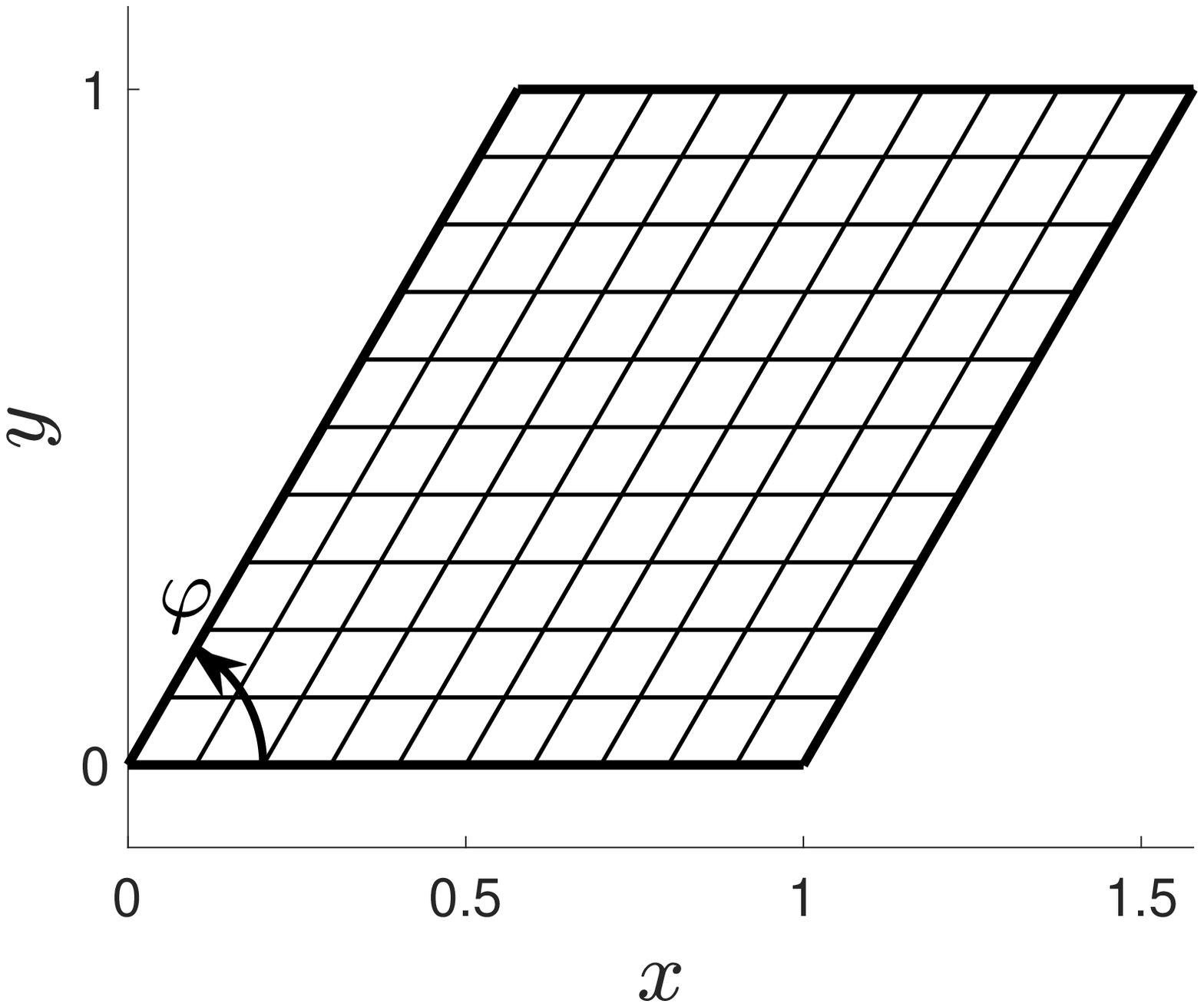}
        \caption{Parallelogram domain}\label{fig:parallelogram}
    \end{subfigure}
    \begin{subfigure}[b]{.5\linewidth}
        \centering
        \includegraphics[width=1\linewidth]{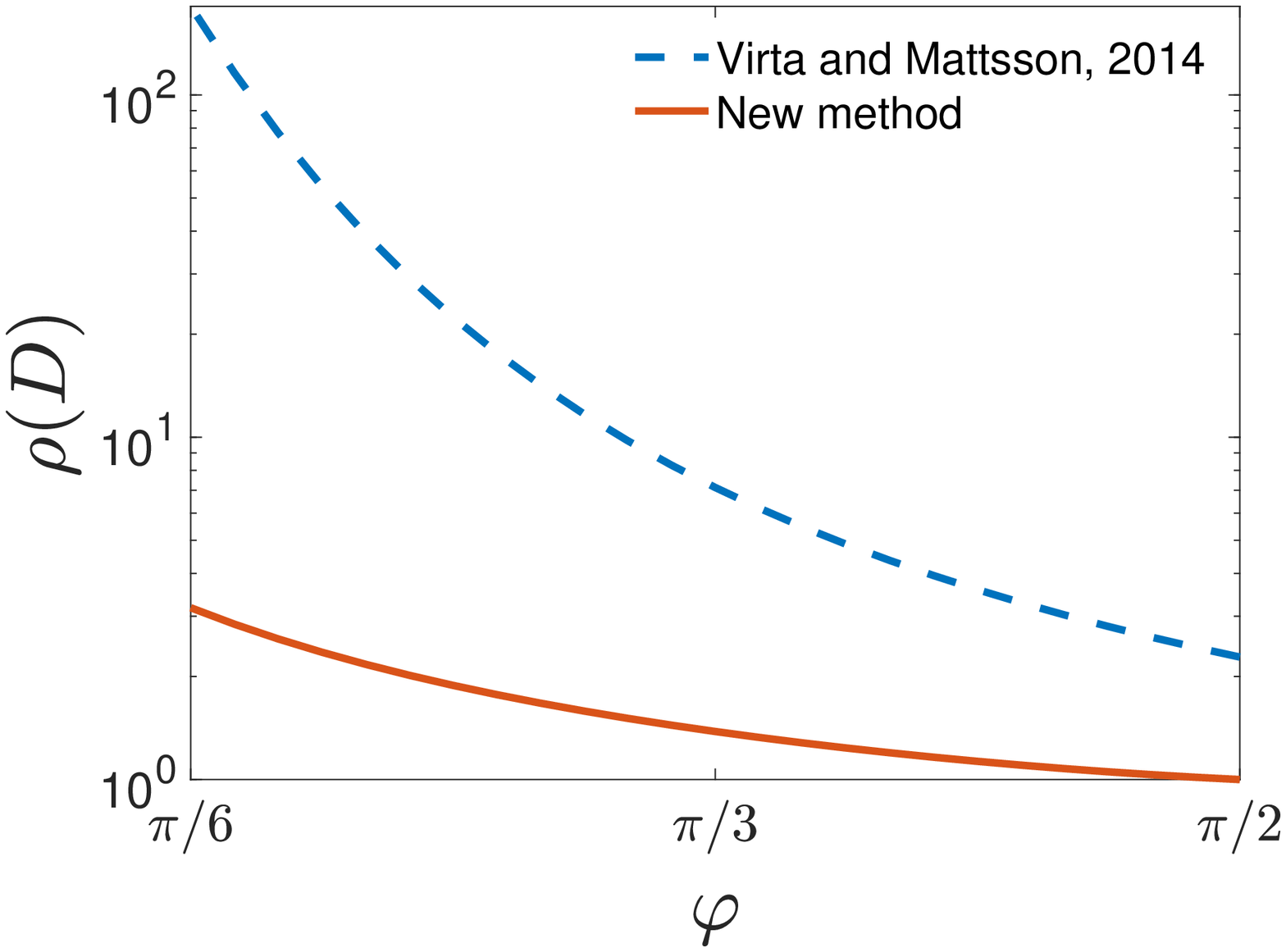}
        \caption{Spectral radius, $q=2$}\label{fig:spectral_radius_2}
    \end{subfigure}
    \caption{(a) The parallelogram-type domain used to compare spectral radii. (b) The spectral radii of the second order discretization matrices, normalized so that the new method has radius 1 for $\varphi = \pi/2$. The new method is much less sensitive to grid skewness.}
    \label{fig:domain_and_spectral_rad_2}
\end{figure}

Let $\physdom$ be the parallelogram domain depicted in Figure \ref{fig:parallelogram}. By varying the angle $\varphi$, we obtain a family of parallelograms whose base and height are of unit length and whose bottom and top boundaries are horizontal. The grid lines are chosen to be equidistant and parallel to the domain boundaries. For $\varphi = \pi/2$, $\physdom$ is square and the grid is Cartesian. As $\varphi$ decreases, grid skewness increases. For the experiments in this subsection, we let the grid consist of $31\times31$ points.

Discretizing \eqref{eq:pde_constant_coeff} in space with the SBP-SAT method yields the semidiscrete equation
\begin{equation}
	\ddot{\bfu} = D \bfu .
\end{equation}
Using the VM SATs and the new SATs derived in later sections of this paper lead to different $D$ matrices, $D_{\text{VM}}$ and $D_{\text{new}}$. Let $\rho(D)$ denote the spectral radius of $D$. Figure \ref{fig:spectral_radius_2} compares $\rho(D_{\text{VM}})$ and $\rho(D_{\text{new}})$ for the second order accurate discretizations, for $\pi/6 \leq \varphi \leq \pi/2$. The spectral radii are normalized by $\rho(D_{\text{new}})$ on the Cartesian grid given by $\varphi = \pi/2$. The new method has a smaller spectral radius for all $\varphi$, but the difference is most apparent for small $\varphi$. As $\varphi$ decreases, the spectral radius of the new method grows much more slowly than that of the VM method.

\begin{figure}[h]
    \begin{subfigure}[b]{.5\linewidth}
        \centering
        \includegraphics[width= 1\linewidth]{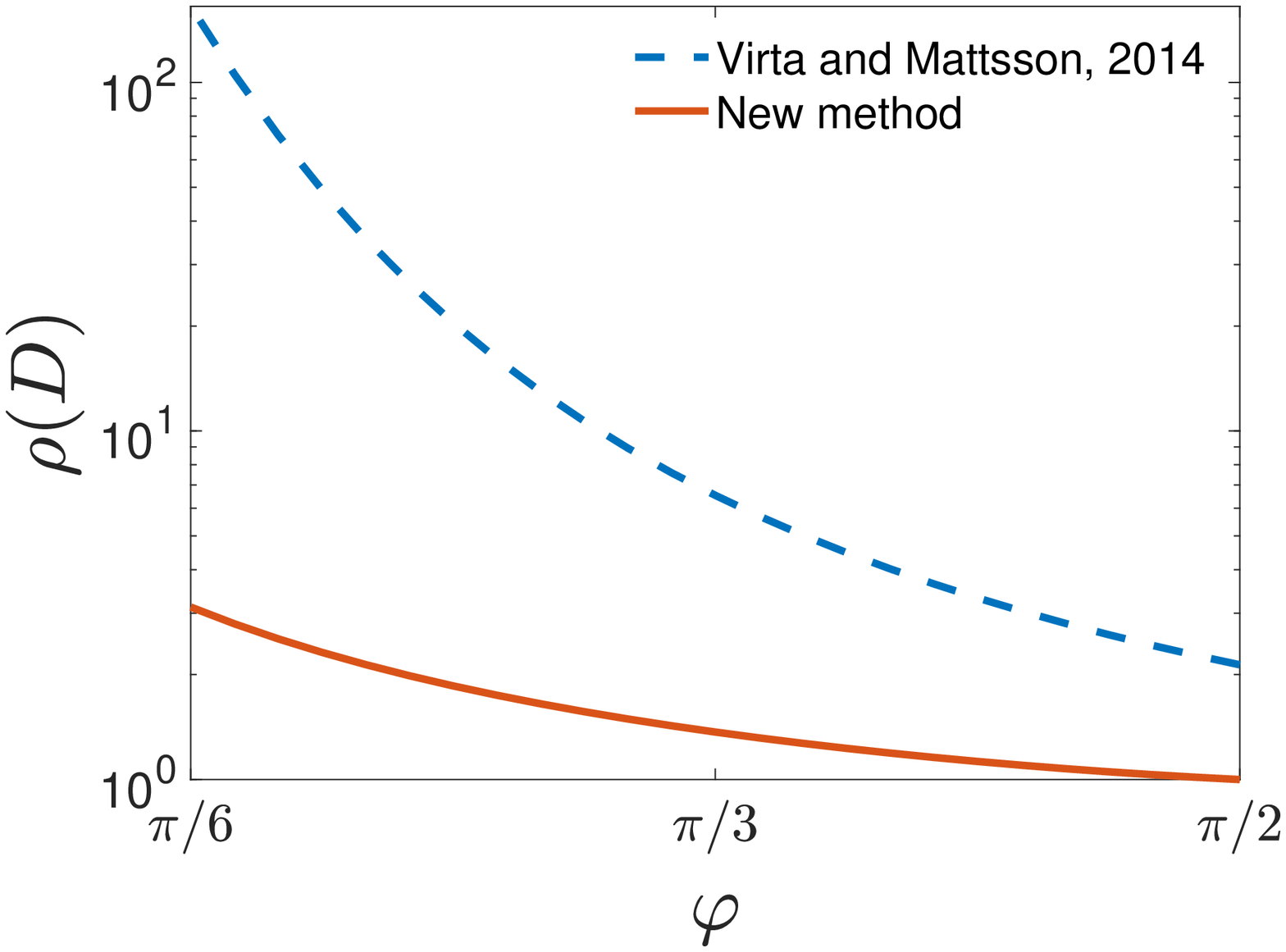}
        \caption{Spectral radius, $q=4$}\label{fig:spectral_radius_4}
    \end{subfigure}
    \begin{subfigure}[b]{.5\linewidth}
        \centering
        \includegraphics[width=1\linewidth]{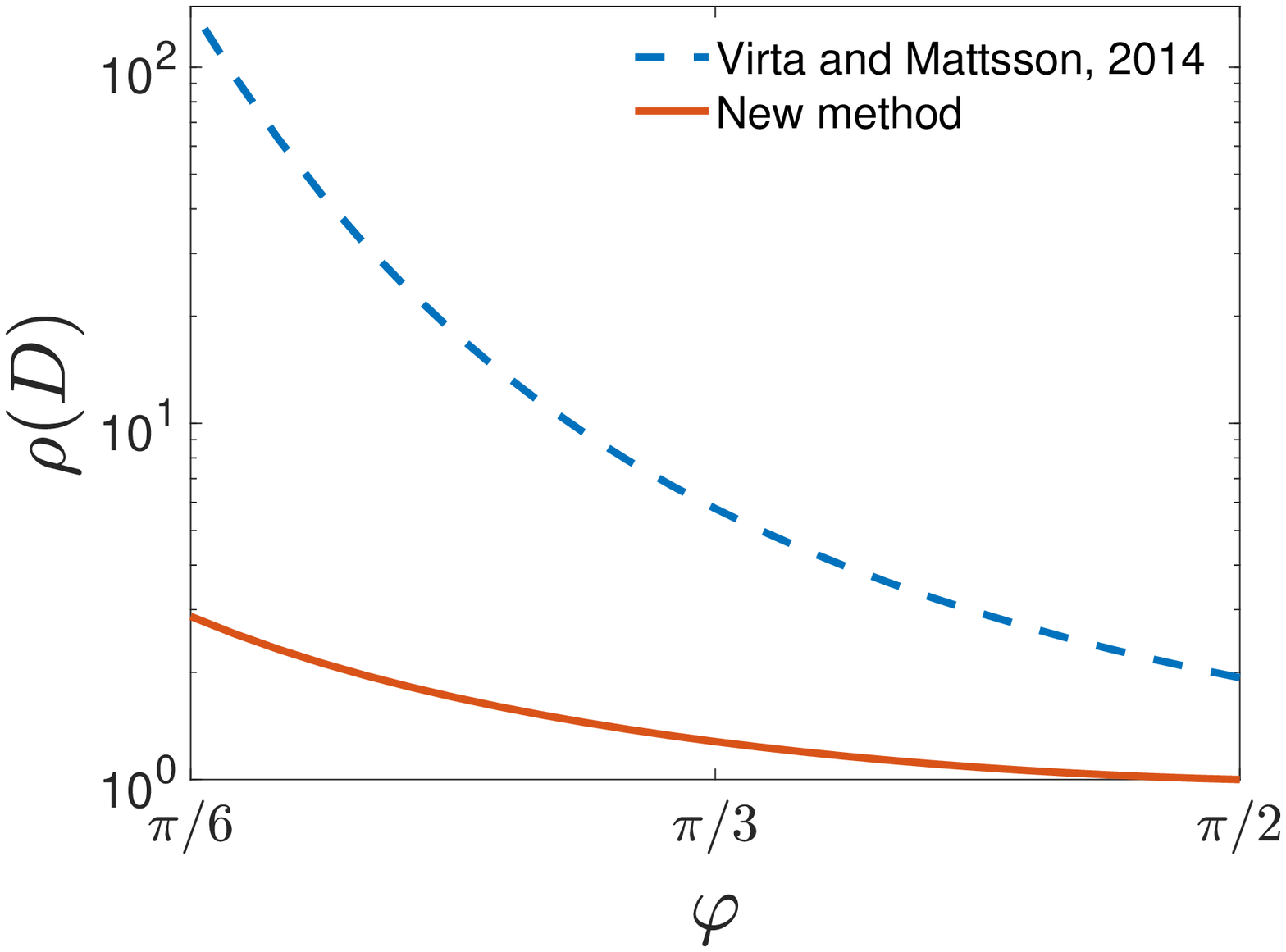}
        \caption{Spectral radius, $q=6$}\label{fig:spectral_radius_6}
    \end{subfigure}
    \caption{The spectral radii of the fourth and sixth order discretization matrices, normalized so that the new method has radius 1 for $\varphi = \pi/2$. The new method is much less sensitive to grid skewness.}
    \label{fig:spectral_rad_46}
\end{figure}

Figure \ref{fig:spectral_rad_46} shows the spectral radii of the fourth and sixth order accurate discretizations, which exhibit similar behavior.

\subsection{Computational efficiency on a multiblock grid} \label{subsec:eff}

\begin{figure}[h]
    \begin{subfigure}[b]{.5\linewidth}
        \centering
        \includegraphics[width=0.97\linewidth,trim={1.8cm 0 2.35cm 0},clip]{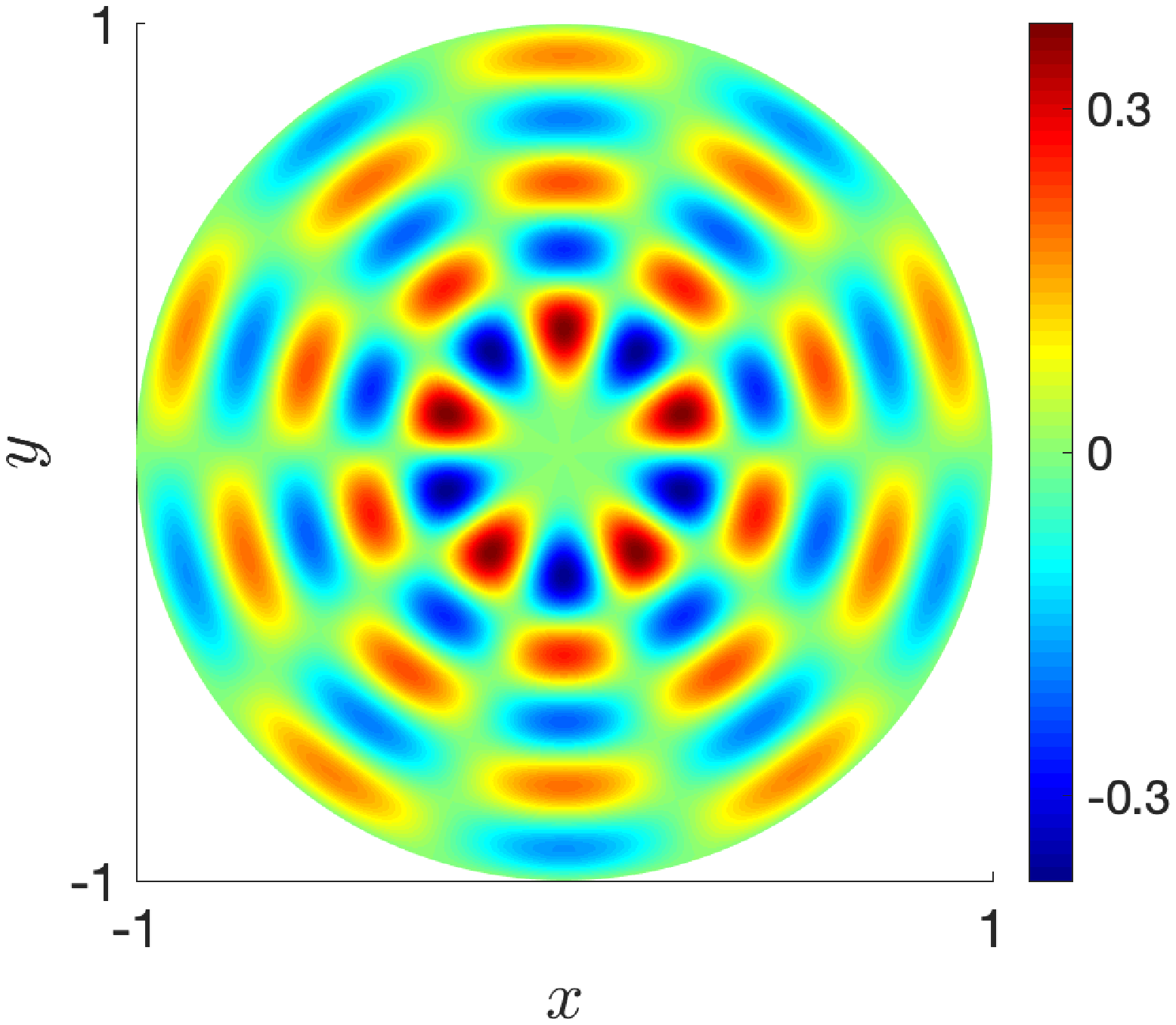}
        \caption{Initial data}\label{fig:bessel}
    \end{subfigure}
    \begin{subfigure}[b]{.5\linewidth}
        \centering
        \includegraphics[width= 0.83\linewidth]{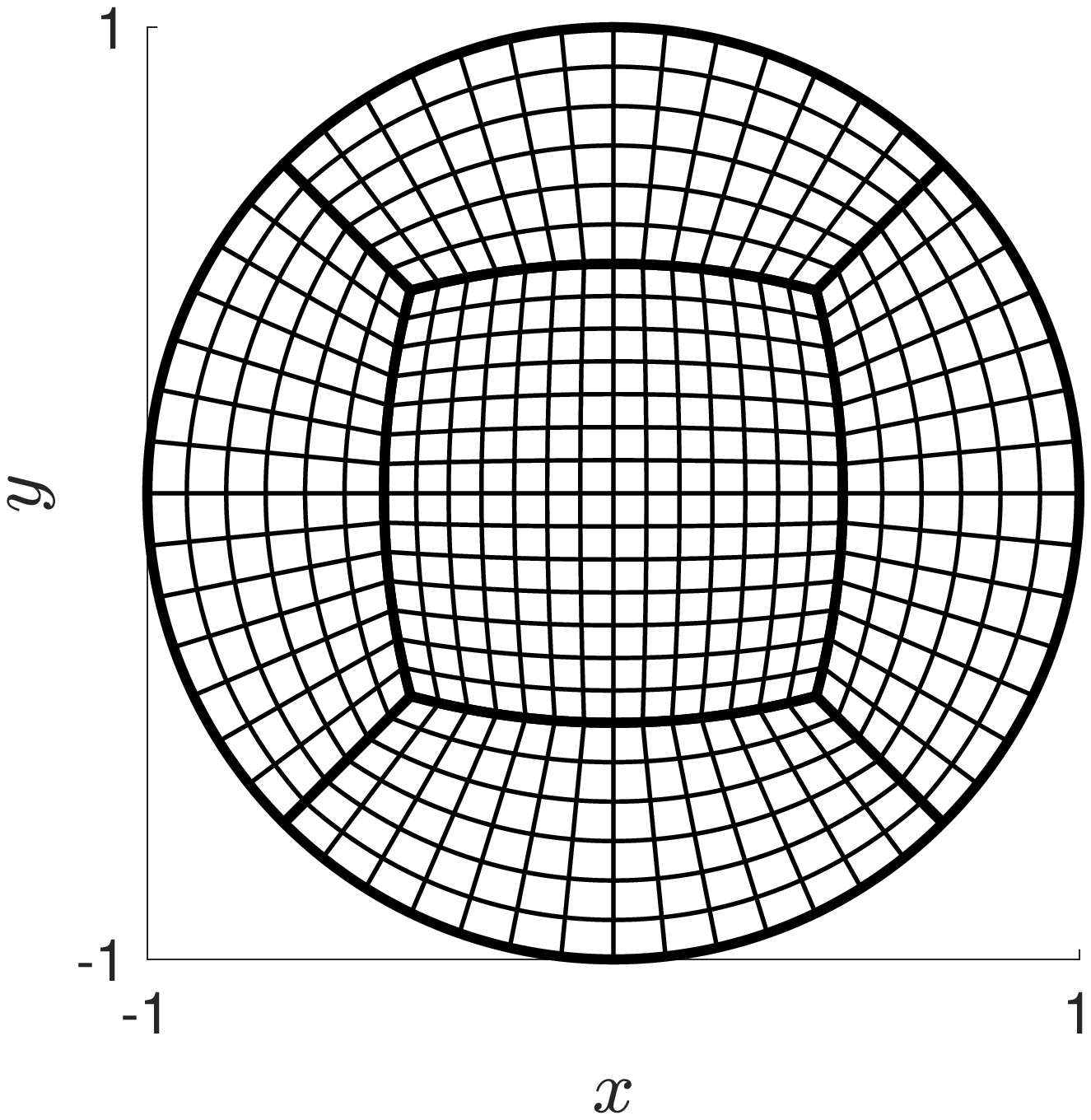}
        \caption{Coarse grid for the unit disk}\label{fig:circle}
    \end{subfigure}
    \caption{(a) Exact solution at time $t=0$. (b) The multiblock grid configuration used in the computations.}
    \label{fig:circle_bessel}
\end{figure}

We here let $\physdom$ be the unit disk, in which case the model problem \eqref{eq:pde_constant_coeff} has exact solutions
\begin{equation}
	u_{\alpha,\beta}(r,\theta,t) =  J_{\alpha}(\beta r) \sin{\alpha \theta} \cos{\beta t},
\end{equation}
where $(r,\theta)$ are the polar coordinates, $\alpha \in \Z$, $J_{\alpha}$ is the order $\alpha$ Bessel function of the first kind, and $\beta$ is a root of $J_{\alpha}$. We set $\alpha = 5$ and $\beta=22.2177998965612$, which is the fifth nonzero root of $J_5$. The resulting solution at time $t=0$ is shown in Figure \ref{fig:bessel}. We use the multiblock grid configuration in Figure \ref{fig:circle}. We rewrite the semidiscrete system of equations as a first order system in time and use the classical fourth order Runge--Kutta method to integrate until time $t=1$, where we measure the relative $\ell^2$ errors. To isolate the accuracy of the spatial discretization, we here use the same time step for the two methods. Let the number of grid points in the center block be $m \times m$ and define the reference grid spacing $h = 1/(m-1)$. We here set $\Delta t = 0.05 h$ for both methods.

\begin{figure}[h]
    \begin{subfigure}[b]{.5\linewidth}
        \centering
        \includegraphics[width= 1\linewidth]{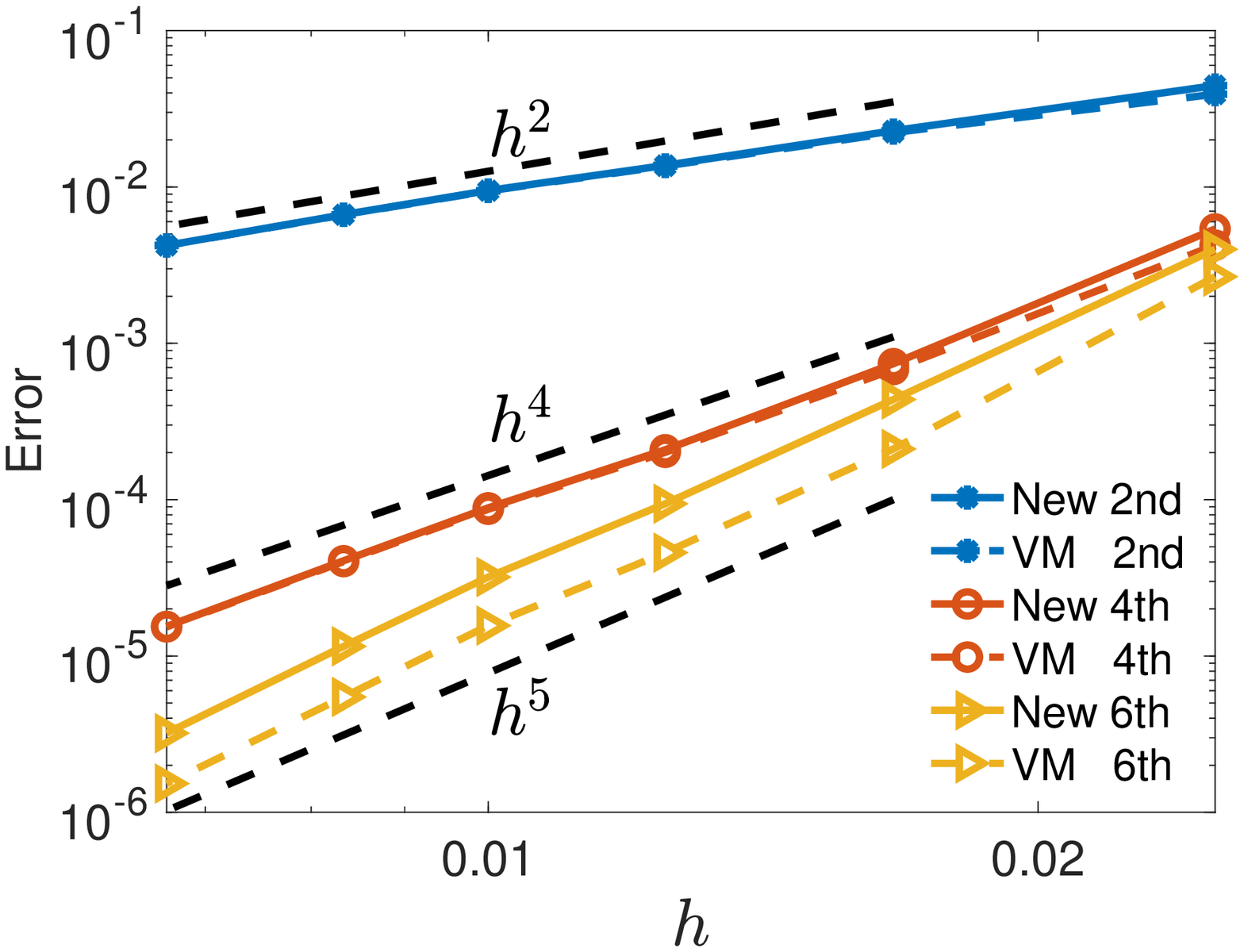}
        \caption{Convergence rates}\label{fig:convergence}
    \end{subfigure}
    \begin{subfigure}[b]{.5\linewidth}
        \centering
        \includegraphics[width=1\linewidth]{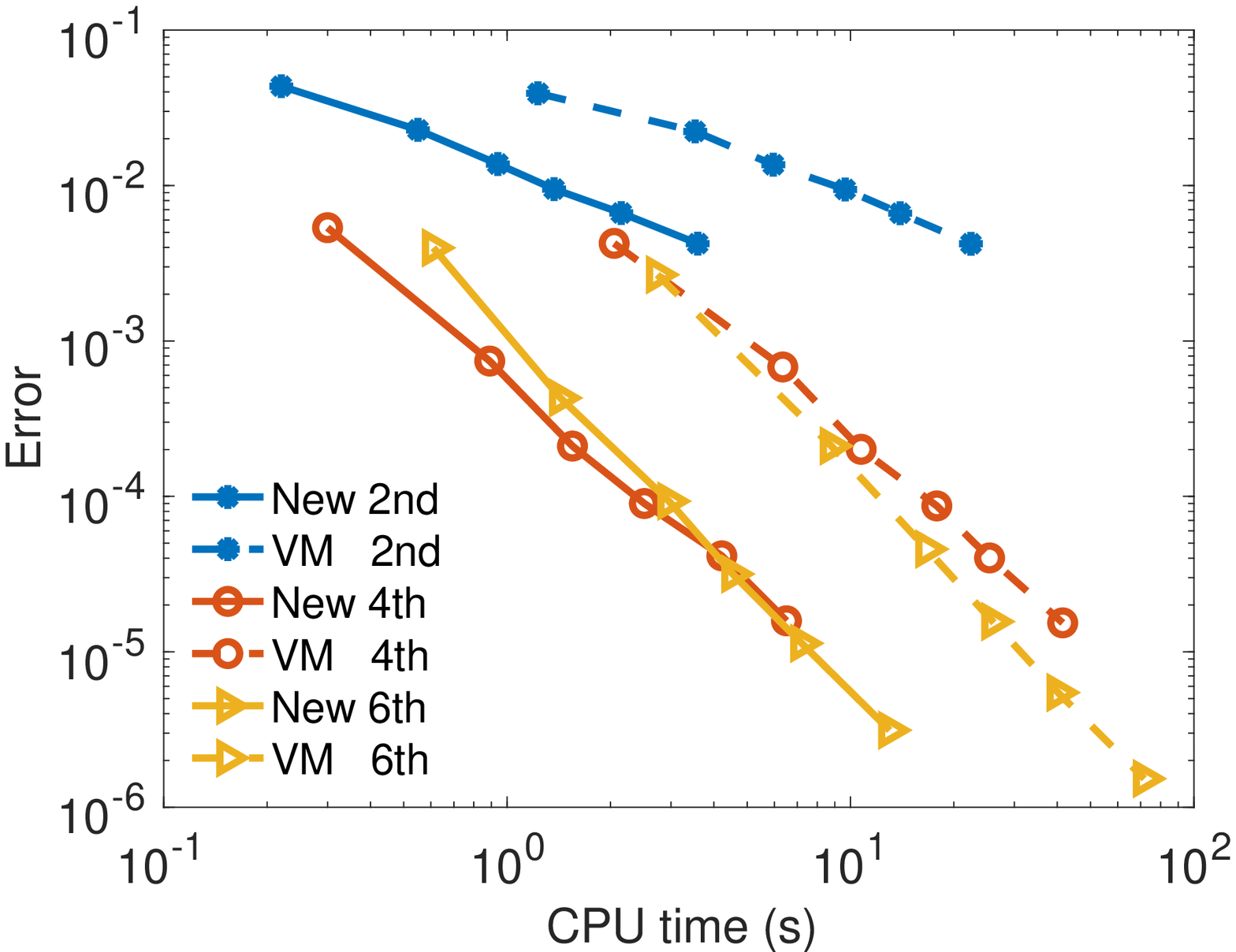}
        \caption{Computational efficiency}\label{fig:efficiency}
    \end{subfigure}
    \caption{$\ell^2$ errors as functions of (a) grid spacing and (b) CPU time.}
    \label{fig:conv_eff}
\end{figure}

Figure \ref{fig:convergence} shows that the two methods converge with the same rates.  The VM method consistently yields the smaller error, which is clearly visible in the sixth order case. For second and fourth order, the difference in error is so small that the curves are nearly indistinguishable.

Of more practical relevance than the error as a function of grid spacing is the CPU time required to satisfy a given error tolerance. To compare the CPU times, we let both methods use the largest stable time step (determined by eigenvalue computations). Since the new method has a smaller spectral radius, it ends up using a significantly larger time step. For an example grid, with the configuration in Figure \ref{fig:circle} and $41 \times 41$ grid points in the center block, the time step ratios are
\begin{equation}
\begin{aligned}
    \text{second order: } &\Delta t_{\text{new}}/\Delta t_{\text{VM}} = 6.69, \\
    \text{fourth order: } &\Delta t_{\text{new}}/\Delta t_{\text{VM}} = 6.46, \\
    \text{sixth order: } &\Delta t_{\text{new}}/\Delta t_{\text{VM}} = 5.84. \\
\end{aligned}
\end{equation}
Figure \ref{fig:efficiency} shows the errors as functions of CPU time. For a given order of accuracy, the larger time steps allow the new method to reach the same error level faster, even though the VM method is slightly more accurate on a given grid.

\section{Notation conventions}
We use the convention to sum over repeated subscript indices so that
\begin{equation}
	u_i v_i = \sum \limits_{i=1}^d u_i v_i,
\end{equation}
where $d$ denotes the number of spatial dimensions. The summation convention applies to inner products too, i.e.,
\begin{equation}
	\ip{u_i}{v_i}{} = \sum \limits_{i=1}^d \ip{u_i}{v_i}{} .
\end{equation}
For variable coefficients $a(\vec{x})$, we use the symbol $a$ also in the discrete case, which then is understood to denote a diagonal matrix with the values of $a(\vec{x})$ on the diagonal. The outward unit normal $\hat{n}$ is regarded as a variable coefficient that takes non-zero values only at boundary points. In the discrete setting, the value of $\hat{n}$ at edge and corner points varies with context. When integrating over a face, $\hat{n}$ is understood to denote the unit normal to that face even at edge and corner points.

\section{Coordinate transformations}
As a model problem that includes the Laplacian, we will consider the second order wave equation,
\begin{equation}
\begin{array}{lll} \label{eq:wave_eq_general}
	a \ddot{u} = \ddxi b \ddxi u, & \vec{x} \in \physdom, & t \in [0, \; T], \\
	Lu = 0, & \vec{x} \in \physboundary, & t \in [0, \; T],
\end{array}
\end{equation}
where $a>0$ and $b>0$ are variable coefficients, $\physboundary = \partial \physdom$, and the linear operator $L$ represents well-posed boundary conditions. Introduce a one-to-one mapping $x_i = x_i(\xi_1,...,\xi_d)$ from $\refdom = [0, \;1]^d$ to $\physdom$. Let $\mathcal{J}$ be the matrix with elements $\mathcal{J}_{ij} = \partial x_j/ \partial \xi_i$, and let $J = \det(\mathcal{J})$. Similarly, let $\K_{ij} = \partial \xi_j / \partial x_i$.
By the chain rule,
\begin{equation} \label{eq:transf_chain}
	\partial_{x_i} = \K_{ij} \partial_{\xi_j}.
\end{equation}
The following metric identities are well known (see \cite{thompson1985numerical}):
\begin{equation} \label{eq:identity}
	J \K_{ij} \partial_{\xi_j} = \partial_{\xi_j} J \K_{ij}.
\end{equation}
Using first \eqref{eq:transf_chain} and then \eqref{eq:identity}, we have
\begin{equation}
	\ddxi b \ddxi  = \K_{ij} \ddxij b \K_{ij} \ddxij = \frac{1}{J} \ddxij J \K_{ij} b \K_{ij} \ddxij .
\end{equation}
For convenience we introduce the notation
\begin{equation}
	\ddi = \ddxii, \quad \alpha_{ijk\ell} = \K_{ij} J \K_{k\ell}.
\end{equation}
It follows that
\begin{equation}
	a \ddot{u} = \ddxi b \ddxi u \Leftrightarrow J a \ddot{u} = \ddj \alpha_{ijik} b \ddk u .
\end{equation}
Henceforth, we will consider the transformed wave equation on the unit cube $\Omega$:
\begin{equation} \label{eq:model_problem_transf}
\begin{array}{lll}
	J a \ddot{u} = \ddj \alpha_{ijik} b \ddk u, & \vec{\xi} \in \refdom, & t \in [0, \; T], \\
	\widetilde{L}u = 0, & \vec{\xi} \in \refboundary, & t \in [0, \; T],
\end{array}
\end{equation}
where $\refboundary = \partial \refdom$ and $\widetilde{L}$ is the transformation of $L$.

\subsection{Integrals}
We will use the $L^2$ inner product and norm:
\begin{equation}
	\volint{u}{v} = \myint{\Omega}{}{uv}{\Omega}, \quad \norm{u}_{\Omega}^2 = \volint{u}{u} .
\end{equation}
Similarly, we write
\begin{equation}
	\surfint{u}{v} = \myint{\Gamma}{}{uv}{\Gamma}, \quad \bnorm{u} = \surfint{u}{u}.
\end{equation}
Note however, that $\surfint{\cdot}{\cdot}$ is not an inner product but a bilinear form, and $\bnorm{\cdot}$ is not a norm but a semi-norm. Further, we define $\gamma$ such that $\mathrm{d} \physboundary = \gamma \, \mathrm{d} \refboundary$ is the surface area element. With this notation, we have
\begin{equation}
	\volintphys{u}{v} = \volint{u}{Jv}
\end{equation}
and
\begin{equation}
 \surfintphys{u}{v} = \surfint{u}{\gamma v}.
\end{equation}

\subsection{Normal derivatives}
Let $\hat{n} = [n_1,\ldots,n_d]$ denote the outward unit normal to $\physdom$, and let $\hat{\nu} = [\nu_1,\ldots,\nu_d] $ denote the unit normal to $\refdom$. Using integration by parts, we have
\begin{equation} \label{eq:ibp1}
	\volintphys{u}{\ddxi v} = \surfintphys{u}{n_i v} - \volintphys{\ddxi u}{v}.
\end{equation}
By instead using the integration-by-parts formula in the reference coordinates, we can also derive
\begin{equation} \label{eq:ibp2}
\begin{aligned}
	\volintphys{u}{\ddxi v} &= \volint{u}{J \ddxi v} = \volint{u}{J \K_{ij} \ddj v} \\
    &= \surfint{u}{J\K_{ij}\nu_j v} - \volint{\ddj J \K_{ij} u}{v} \\
	&= \surfint{u}{J\K_{ij}\nu_j v} - \volint{J \K_{ij} \ddj u}{v} \\
    &= \surfint{u}{\gamma \frac{J}{\gamma}\K_{ij}\nu_j v} - \volint{\ddxi u}{Jv} \\
	&= \surfintphys{u}{\frac{J}{\gamma}\K_{ij}\nu_j v} - \volintphys{\ddxi u}{v}. \\
\end{aligned}
\end{equation}
where we used \eqref{eq:identity}. Since $u$ and $v$ are arbitrary in \eqref{eq:ibp1} and \eqref{eq:ibp2}, it follows that the normal components are related through
\begin{equation}
	n_i = \frac{J}{\gamma}\K_{ij}\nu_j .
\end{equation}
Hence, the normal derivative on $\physdom$ can be expressed as
\begin{equation} \label{eq:normal_derivative_cont}
	\pdd{u}{\hat{n}} = n_i \partial_{x_i} u = \nu_j \frac{J}{\gamma} \K_{ij} \K_{ik} \ddk u = \nu_j \frac{\alpha_{ijik} }{\gamma} \ddk u.
\end{equation}

\section{Summation-by-parts operators}
Consider finite-difference approximations of derivatives using an equidistant grid that includes the boundary points of an interval $[x_{\ell}, x_{r}]$. In this paper, we consider diagonal-norm SBP operators only. That is, the so-called norm matrix $H_x$ has the structure
\begin{equation} \label{eq:structure_H}
    H_{x} = h \, \mbox{diag}(w_1, \ldots, w_{s}, 1, \ldots, 1,  w_{s}, \ldots, w_{1}),
\end{equation}
where $h$ denotes the grid spacing and the $w_i$ are positive dimensionless weights. The first-derivative SBP operators $D_x \approx \partial_x$ have the property
\begin{equation}
	H_x D_x + D_x^T H_x = - e_{\ell} e_{\ell}^T + e_r e_r^T \quad \mbox{(no sums over } x, \ell, r),
\end{equation}
where
\begin{equation} \label{eq:el_er}
e_{\ell} = \begin{bmatrix} 1, 0 , \ldots, 0\end{bmatrix}^T, \quad e_{r} = \begin{bmatrix} 0, \ldots, 0, 1 \end{bmatrix}^T.
\end{equation}
The narrow-stencil second derivative operators $D_{xx}(b) \approx \partial_x b \partial_x$ derived in \cite{Mattsson11} are based on the same norm matrix and have the property
\begin{equation} \label{eq:sbp_property_2nd_der}
\begin{aligned}
	 H_x D_{xx}(b) &= -D_x^T H_x b D_x - R_{xx}(b) - e_{\ell} b_{\ell} e_{\ell}^T \db_x + e_r b_{r} e_r^T \db_x, \\
     &\mbox{(no sums over } x, \ell, r),
\end{aligned}
\end{equation}
where $\db_x$ approximates the first derivative. Note that for the SBP operators considered in this paper, $D_x \neq \db_x$. If $D_x=\db_x$, then the SBP operators $D_x$ and $D_{xx}$ are said to be \emph{fully compatible} \cite{MattssonParisi09}. The difference between the two first-derivative approximations will play an important role in the coming stability analysis and we define
\begin{equation}
    \Delta D_x = D_x - \db_x.
\end{equation}
The SBP operators derived in \cite{Mattsson11} have $\db_x$ that are accurate of order $q/2+1$, i.e., one order higher than the boundary closure of $D_x$. This implies that $\Delta D_x \bfu$ is zero to order $q/2$ for all $\bfu$ that are restrictions of smooth functions to the grid.

The matrix $R_{xx}(b)$ is symmetric positive semidefinite and consists of undivided difference approximations in such a way that $\bfu^T R_{xx}(b) \bfu$ is zero to order $q$ \cite{Mattsson11}.

\subsection{Positivity properties}
It follows immediately from \eqref{eq:structure_H} and \eqref{eq:el_er} that we have
\begin{equation} \label{eq:pos_H_1d}
\begin{aligned}
    H_{x} &= h \, \mbox{diag}(0,w_2, \ldots, w_{s}, 1, \ldots, 1,  w_{s}, \ldots, w_2,0) \\
    &+ h w_{1} e_{\ell} e_{\ell}^T + h w_{1} e_{r} e_{r}^T  \\
    & \geq h w_{1} e_{\ell} e_{\ell}^T + h w_{1} e_{r} e_{r}^T.
\end{aligned}
\end{equation}
We define $\theta_H = w_1$ and write
\begin{equation}
   H_x \geq h \theta_H e_{\ell} e_{\ell}^T + h\theta_H e_{r} e_{r}^T,
\end{equation}
or, equivalently,
\begin{equation} \label{eq:borrowing_H_1d}
   \bfu^T H_x \bfu \geq  h\theta_H (e_{\ell}^T \bfu)^2 + h\theta_H (e_{r}^T \bfu)^2 \quad \mbox{for all } \bfu.
\end{equation}

Let $x_{\ell, m_b}$ and $x_{r, m_b}$  denote the $m_b$th grid points from the left and right boundaries:
\begin{equation}
    x_{\ell, m_b} = x_{\ell}+(m_b-1)h, \quad x_{r, m_b} = x_{r} - (m_b-1)h .
\end{equation}
Let $b_{\ell,min}$ and $b_{r,min}$ denote the minimum of $b$ over the $m_b$ leftmost and rightmost grid points:
\begin{equation}
\begin{aligned}
    b_{\ell,min} &= \min \left( b(x_{\ell}), b(x_{\ell}+h), \ldots, b(x_{\ell,m_b}) \right), \\
    b_{r,min} &= \min \left( b(x_{r,m_b}), b(x_{r,m_b}+h), \ldots, b(x_{r}) \right).
\end{aligned}
\end{equation}
The SBP operators considered in this paper also satisfy the following positivity property, which is essential to the new method.
\begin{lemma} \label{lemma:R}
The second-derivative SBP operators $D_{xx}(b)$ of orders two, four, and six, derived in \cite{Mattsson11}, which can be decomposed as in \eqref{eq:sbp_property_2nd_der}, all satisfy
\begin{equation} \label{eq:borrowing_R_1d}
    \bfu^T R_{xx}(b) \bfu \geq h \theta_R b_{\ell,min} (e_{\ell}^T \Delta D_x \bv{u})^2 + h \theta_R b_{r,min} (e_r^T \Delta D_x \bv{u})^2 \;\mbox{for all } \bfu,
\end{equation}
for some choice of $m_b$ and the dimensionless constant $\theta_R>0$.
\end{lemma}
\begin{proof}
See Appendix.
\end{proof}
A property similar but not identical to \eqref{eq:borrowing_R_1d} was used in \cite{Duru201437}, to prove stability for the elastic wave equation with block interfaces. Because \eqref{eq:borrowing_R_1d} has never before been used in the SBP-SAT literature, we have determined the values of $\theta_R$ numerically. We remark that $\theta_R$ depends on $m_b$. We describe how $m_b$ is chosen and how $\theta_R$ is computed in Appendix. Table \ref{table:borrowing} lists the values of all the dimensionless constants for the SBP operators used in this paper.
\begin{table}[H]
\centering
\begin{tabular}{cccc}
& $\theta_H$ & $\theta_R$ & $m_b$ \\
\hline
second order & $1/2$ & $1.0000$ & $2$  \\
fourth order & $17/48$ & $0.5776$ & $4$  \\
sixth order & $13649/43200$ & $0.3697$ & $7$  \\
\end{tabular}
\caption{Values of $\theta_H$ (exact), $\theta_R$ (approximations, computed in Appendix) and $m_b$. }
\label{table:borrowing}
\end{table}

\subsection{Multi-dimensional first-derivative operators}
Let operators with subscripts ${\xi_i}$ denote one-dimensional operators corresponding to coordinate direction $\xi_i$. The multi-dimensional first derivatives $D_i \approx \ddi $ and $\db_i \approx \ddi$ are constructed using tensor products:
\begin{equation}
	D_i = I_{\xi_1} \otimes I_{\xi_2} \otimes \cdots \otimes I_{\xi_{i-1}} \otimes D_{\xi_i} \otimes \cdots \otimes I_{\xi_{d}} ,
\end{equation}
\begin{equation}
	\db_i = I_{\xi_1} \otimes I_{\xi_2} \otimes \cdots \otimes I_{\xi_{i-1}} \otimes \db_{\xi_i} \otimes \cdots \otimes I_{\xi_{d}} ,
\end{equation}
where the $I_{\xi_i}$ are one-dimensional identity matrices of appropriate sizes. In analogy with the chain rule \eqref{eq:transf_chain}, we define
\begin{equation} \label{eq:chain_discrete}
    D_{x_i} = \K_{ij} D_j .
\end{equation}
We further define the difference between the two first-derivative approximations:
\begin{equation}
	\Delta D_i = D_i - \db_i .
\end{equation}
The multi-dimensional quadrature is
\begin{equation}
	H = H_{\xi_1} \otimes \cdots \otimes H_{\xi_m} .
\end{equation}
Let $\Gamma_{i}^-$ and $\Gamma_i^+$ denote the boundary faces where $\xi_i = 0$ and $\xi_i = 1$, respectively.
For integration over boundary faces, we define
\begin{equation}
	H_{\Gamma_i} = H_{\xi_1} \otimes \cdots \otimes H_{\xi_{i-1}} \otimes H_{\xi_{i+1}} \otimes \cdots \otimes H_{\xi_m}.
\end{equation}
Note that $H_{\Gamma_i}$ can be used to integrate over $\Gamma_i^+$ as well as $\Gamma_i^-$.

For discrete integration over the volume, we write
\begin{equation}
	\dvolint{\bfu}{\bfv} = \bv{u}^T H \bv{v} .
\end{equation}
Let $e_{S}$ denote a restriction operator that picks out only those solution values that reside on the boundary segment $S$. For discrete integration over the face $\Gamma_i^+$, for example, we write
\begin{equation}
	\dsurfintindexsign{\bv{u}}{\bv{v} }{i}{+} = (e_{\Gamma_i^+}^T \bv{u})^T H_{\Gamma_i} (e_{\Gamma_i^+}^T \bv{v}).
\end{equation}
If the surface $S$ is the union of several faces, $S = F_1 \cup ... \cup F_m$, we use the convention that
\begin{equation}
	\bip{\bv{u}}{\bv{v}}{S,h} = \bip{\bv{u}}{\bv{v}}{F_1,h} + \ldots + \bip{\bv{u}}{\bv{v}}{F_m,h} ,
\end{equation}
 i.e., the integration is performed over one face at a time. If the integrand contains the unit normal, its values at edges and corners are defined to be the same as on the remainder of that face.
As in the continuous case, we define the norm and semi-norm
\begin{equation}
	\dnorm{\bfu} = \dvolint{\bfu}{\bfu}, \quad \norm{\bfu}_{S,h}^2 = \bip{\bfu}{\bfu}{S,h}.
\end{equation}

With the notation established in this section, we have the discrete integration-by-parts formula
\begin{equation} \label{eq:sbp_multid_mixed}
	\dvolint{\bfu }{D_i b D_j \bfv} = \dsurfint{\bfu }{\nu_i b D_j \bfv} - \dvolint{D_i \bfu }{b D_j \bfv}.
\end{equation}

\subsection{Multi-dimensional narrow-stencil second-derivative operators}
For any fixed $i$, we construct:
\begin{equation}
	D_{ii}(b) \approx \ddi b \ddi ,
\end{equation}
by using the one-dimensional operator $D_{xx}$ for each grid line. This gives the multi-dimensional SBP property for the second derivative,
\begin{equation} \label{eq:sbp_property_multid_narrow}
\begin{aligned}
	HD_{ii}(b) &= -D_{i}^T Hb D_i - R_{ii}(b) - e_{\Gamma_i^-} b_{\Gamma_i^-} e_{\Gamma_i^-}^T \db_i + e_{\Gamma_i^+} b_{\Gamma_i^+} e_{\Gamma_i^+}^T \db_i, \\
     &\mbox{(no sum over $i$)}.
\end{aligned}
\end{equation}
The structure of $R_{ii}(b)$ in \eqref{eq:sbp_property_multid_narrow} is important and requires a closer look. Let $\widetilde{R}_{ii}(b)$ be the operator that performs the action of $R_{xx}(b)$ for each grid line in direction $i$. Further, define
\begin{equation}
    H_i = I_{\xi_1} \otimes \cdots \otimes I_{\xi_{i-1}} \otimes H_{\xi_i} \otimes I_{\xi_{i+1}} \otimes \cdots \otimes I_{\xi_m} .
\end{equation}
Then,
\begin{equation}
R_{ii}(b) = H^{-1}_i H \widetilde{R}_{ii}(b) \quad \mbox{(no sum over $i$)}.
\end{equation}
The $R_{ii}$ matrices are symmetric positive semidefinite.

The property \eqref{eq:sbp_property_multid_narrow} can be equivalently written as
\begin{equation} \label{eq:sbp_multid_narrow}
\begin{aligned}
	\dvolint{\bfu }{D_{ii}( b ) \bfv} &= \dsurfint{\bfu }{\nu_i b \db_i \bfv} - \dvolint{D_i \bfu }{b D_i \bfv} - \bfu^T R_{ii}(b) \bfv, \\
    &\mbox{(no sum over $i$)},
\end{aligned}
\end{equation}
for all $\bfu$ and $\bfv$.

\subsection{Combining narrow-stencil derivatives and mixed derivatives}
Consider discretizing a term $\ddi b \ddj$. Because we want to use narrow-stencil second derivatives for the non-mixed second derivatives, the following notation will be useful:
\begin{equation}
	\cD_{ij}(b) = \left\{ \begin{array}{cc} D_{ij}(b), & i=j \\ D_i b D_j, & i\neq j \end{array} \right. .
\end{equation}
Combining the two integration-by-parts formulas \eqref{eq:sbp_multid_mixed} and \eqref{eq:sbp_multid_narrow} leads to the integration-by-parts formula
\begin{equation} \label{eq:discr_ibp_multid}
\begin{aligned}
	&\dvolint{\bfu }{\cD_{ij}( b ) \bfv} = \\
	&\left\{
	\def\arraystretch{2.2}
	\begin{array}{lc} \displaystyle\dsurfint{\bfu }{\nu_i b D_j \bfv} - \dvolint{D_i \bfu }{b D_j \bfv}, & i\neq j \\
	\dsurfint{\bfu }{\nu_i b \db_j \bfv}-\dvolint{D_i \bfu }{b D_j \bfv} - \bfu^T R_{ij}(b) \bfv, & i=j \end{array} \right. .
\end{aligned}
\end{equation}

\subsection{Multi-dimensional positivity properties}
 To suppress unnecessary notation, we assume that the grid spacing in the reference domain is $h$ in each dimension (we stress that the analysis does not rely on this assumption). It follows from  \eqref{eq:borrowing_H_1d}
that, for any fixed $j$,
\begin{equation} \label{eq:borrowing_nd_H_local}
    \dnorm{\bfu} \geq h \theta_H \dbnorm{\nu_j \bfu}.
\end{equation}
Using \eqref{eq:borrowing_nd_H_local} we can derive
\begin{equation}
    \dnorm{\bfu} = \frac{1}{d} \sum \limits_{j=1}^d \dnorm{\bfu} \geq \frac{1}{d} \sum \limits_{j=1}^d h\theta_H \dbnorm{\nu_j \bfu} = \frac{h\theta_H}{d} \dbnorm{\bfu},
\end{equation}
which we summarize as
\begin{equation} \label{eq:borrowing_nd_H_global}
    \dnorm{\bfu} \geq \frac{h\theta_H}{d} \dbnorm{\bfu} .
\end{equation}

The multi-dimensional positivity property that follows from \eqref{eq:borrowing_R_1d} is
\begin{equation} \label{eq:borrowing_nd_R}
    \bfu^T R_{ii}(b) \bfu \geq h\theta_R  \dbnorm{\nu_i \sqrt{b_{min}} \Delta D_i \bv{u}} \quad \mbox{(no sum over $i$)},
\end{equation}
where $b_{min}$ denotes the minimum of $b$ over the $m_b$ grid points closest to the boundary $\Gamma$ along each grid line. That is, to evaluate $b_{min}$ at a given grid point on $\Gamma$, one needs to evaluate $b$ at the $m_b$ first grid points on the same grid line orthogonal to $\Gamma$.

\section{Continuous energy balance}
We return to the transformed equation \eqref{eq:model_problem_transf},
\begin{equation} \label{eq:model_problem_transf_repeated}
	J a \ddot{u} = \ddj \alpha_{ijik} b \ddk u.
\end{equation}
The energy method, i.e.\ multiplying \eqref{eq:model_problem_transf_repeated} by $\dot{u}$ and integrating over $\Omega$, yields
\begin{equation} \label{eq:cont_energy_derivation}
\begin{aligned}
	\volint{\dot{u}}{Ja \ddot{u}} &=  \volint{\ut}{\ddj \alpha_{ijik} b \ddk u}
	= \surfint{\ut}{\nu_j \alpha_{ijik} b \ddk u} - \volint{\ddj \ut}{\alpha_{ijik} b \ddk u} \\
	&= \surfint{\ut}{\nu_j \alpha_{ijik} b \ddk u} - \volint{\K_{ij} \ddj \ut}{J \K_{ik} b \ddk u} \\
	& = \surfintphys{\ut}{b \nu_j \frac{\alpha_{ijik} }{\gamma} \ddk u} - \volintphys{\K_{ij} \ddj \ut}{b\K_{ik} \ddk u} \\
	&= \surfintphys{\ut}{b \pdd{u}{\hat{n}} } - \volintphys{\K_{ij} \ddj \ut}{b\K_{ik} \ddk u} .
\end{aligned}
\end{equation}
The volume integrals satisfy
\begin{equation}
	\volint{\dot{u}}{Ja \ddot{u}} = \frac{1}{2} \dd{}{t} \normphys{\sqrt{a}\dot{u}}^2,
\end{equation}
and
\begin{equation}
\begin{aligned}
	\volintphys{\K_{ij} \ddj \ut}{\K_{ik} b \ddk u} &= \frac{1}{2} \dd{}{t}  \si \normphys{\sqrt{b}\K_{ij}\ddj u}^2  \\
    &= \frac{1}{2} \dd{}{t}  \si \normphys{\sqrt{b}\ddxi u}^2.
\end{aligned}
\end{equation}
We define the energy $\mathcal{E}$ as
\begin{equation} \label{eq:cont_energy_def}
	\mathcal{E} = \frac{1}{2} \normphys{\sqrt{a}\dot{u}}^2 + \frac{1}{2} \si \normphys{\sqrt{b}\ddxi u}^2,
\end{equation}
which allows us to write \eqref{eq:cont_energy_derivation} as
\begin{equation} \label{eq:energy_balance_cont}
	\dd{\mathcal{E}}{t} = \surfintphys{\ut}{b \pdd{u}{\hat{n}}}.
\end{equation}

\section{Semidiscrete energy balance}
We discretize the transformed problem \eqref{eq:model_problem_transf} in space as
\begin{equation} \label{eq:semidiscrete_general}
	J a \ddot{\bv{u}} = \cD_{jk}(\alpha_{ijik} b) \bv{u}
	+ SAT ,
\end{equation}
where $SAT$ imposes the boundary conditions and will be specified later. The discrete energy method amounts to multiplying \eqref{eq:semidiscrete_general} by $\dot{\bfu}^T H$ from the left, which leads to
\begin{equation}
    \dvolint{\dot{\bfu}}{J a \ddot{\bv{u}}} = \dvolint{\dot{\bfu}}{\cD_{jk}(\alpha_{ijik} b) \bv{u}}
    + \dot{\bfu}^T H (SAT).
\end{equation}
We have
\begin{equation}
	\dvolint{\bfut }{Ja \ddot{\bv{u}}} =
\frac{1}{2}  \dd{}{t} \si \norm{\sqrt{Ja} \bfut }^2_{\Omega,h} = \frac{1}{2} \dd{}{t} \dnormphys{\sqrt{a} \bfut}^2 .
\end{equation}
The integration-by-parts formula \eqref{eq:discr_ibp_multid} yields
\begin{equation}
\begin{aligned}
	&\dvolint{ \bfut }{\cD_{jk}(\alpha_{ijik} b) \bfu}= \\
    &= \sk \dsurfint{\bfut }{\nu_k \alpha_{ikik}  b \db_{k} \bfu} + \skjnk \dsurfint{\bfut }{\nu_j \alpha_{ijik}  b D_{k} \bfu} \\
	&- \dvolint{D_j  \bfut }{\alpha_{ijik}  b D_k \bfu} - \sk \bfut^T R_{kk}(\alpha_{ikik}  b) \bfu \\
    &= \dsurfint{\bfut }{\nu_j \alpha_{ijik}  b D_{k} \bfu} - \sk \dsurfint{\bfut }{\nu_k \alpha_{ikik}  b \Delta D_{k} \bfu} \\
    &- \dvolint{\K_{ij} D_j \bfut }{ J b \K_{ik}  D_k \bfu} - \sk \bfut^T R_{kk}(\alpha_{ikik}  b) \bfu \\
	&= \dsurfint{\bfut }{\nu_j \alpha_{ijik}  b D_{k} \bfu} - \sk \dsurfint{\bfut }{\nu_k \alpha_{ikik}  b \Delta D_{k} \bfu} \\
	&- \frac{1}{2} \dd{}{t} \left(\si \dnormphys{\sqrt{b} D_{x_i} \bfu}^2 + \sk \bfu^T R_{kk}(\alpha_{ikik}  b) \bfu \right).
\end{aligned}
\end{equation}
We define the discrete normal derivative $D_{\hat{n}}$ as
\begin{equation} \label{eq:normal_derivative_discr}
	D_{\hat{n}} = \nu_j \frac{\alpha_{ijik}}{\gamma} D_{k} - \sk \nu_k \frac{\alpha_{ikik}}{\gamma} \Delta D_{k},
\end{equation}
and the discrete energy $E$ as
\begin{equation}
	E = \frac{1}{2} \dnormphys{\sqrt{a} \bfut}^2 + \frac{1}{2} \si \dnormphys{\sqrt{b} D_{x_i} \bfu}^2 + \frac{1}{2} \sk \bfu^T R_{kk}(\alpha_{ikik}  b) \bfu,
\end{equation}
so that the discrete energy rate reads
\begin{equation} \label{eq:energy_rate_discr}
    \dd{E}{t} = \dsurfintphys{\bfut }{b D_{\hat{n}} \bfu} + \bfut^T H (SAT).
\end{equation}
By comparing \eqref{eq:normal_derivative_discr} with \eqref{eq:normal_derivative_cont} and recalling that $\Delta D_k \bfu$ is zero to order $q/2$, we conclude that $D_{\hat{n}}$ approximates $\partial / \partial {\hat{n}}$. Further, since the $R_{kk}$ are positive semidefinite and $\bfu^T R_{kk} \bfu$ is zero to order $q$, the discrete energy $E$ is nonnegative and approximates the continuous energy $\mathcal{E}$ defined in \eqref{eq:cont_energy_def}.

\subsection{Enforcing Neumann boundary conditions}
Homogeneous Neumann boundary conditions, i.e.,
\begin{equation} \label{eq:bc_neumann_cont}
	\pdd{u}{\hat{n}} = 0, \quad \vec{x} \in \physboundary,
\end{equation}
are energy-conserving for the continuous equation, c.f.\ \eqref{eq:energy_balance_cont}. By inspecting \eqref{eq:energy_rate_discr}, we see that the discrete energy $E$ will be conserved, and hence the scheme will be energy-stable, if we can find a $SAT$ that is consistent with \eqref{eq:bc_neumann_cont} and satisfies
\begin{equation} \label{eq:energy_preserving_sat}
	\bfut^T H (SAT) = - \dsurfintphys{\bfut }{b D_{\hat{n}} \bfu} .
\end{equation}
Let $F$ denote the set of all faces, i.e.\ $F= \{\Gamma_1^-,\ldots,\Gamma_d^-,\Gamma_1^+, \ldots, \Gamma_d^+\}$. We achieve \eqref{eq:energy_preserving_sat} by setting
\begin{equation}
	SAT = -H^{-1} \sum \limits_{f\in F} e_f \gamma H_f b (e_f^T D_{\hat{n}} \bfu - 0) .
\end{equation}

\subsection{Enforcing Dirichlet boundary conditions}
Homogeneous Dirichlet conditions, i.e.,
\begin{equation} \label{eq:bc_dirichlet_cont}
	u = 0, \quad \vec{x} \in \physboundary,
\end{equation}
are energy-conserving for the continuous problem. As opposed to the case with Neumann conditions, we cannot choose consistent SATs that make the boundary terms in the energy rate vanish. Instead, we will show that a slightly different energy $E_D \approx E$ is conserved.

General SATs of the form
\begin{equation}
	SAT = H^{-1} \sum \limits_{f \in F} B^T e_{f} \gamma H_{f} (e_{f}^T \bfu - 0)
\end{equation}
lead to
\begin{equation}
	\bfut^T H (SAT) = \dsurfintphys{B \bfut }{ \bfu} .
\end{equation}
We set $B =  C - A$, where $A$ is symmetric. Inserting this ansatz in the discrete energy rate \eqref{eq:energy_rate_discr} yields
\begin{equation}
	\dd{E}{t} = \dsurfintphys{\bfut }{b D_{\hat{n}} \bfu} + \dsurfintphys{C \bfut }{ \bfu} - \dsurfintphys{A \bfut }{ \bfu} .
\end{equation}
For symmetry we will need
\begin{equation} \label{eq:penalty_condition}
\begin{aligned}
	\dsurfintphys{C \bfut }{ \bfu} &= \dsurfintphys{b D_{\hat{n}} \bfut }{ \bfu},
\end{aligned}
\end{equation}
which implies $C = b D_{\hat{n}}$. We then obtain
\begin{equation}
\begin{aligned}
	\dd{E}{t} &= \dsurfintphys{\bfut }{b D_{\hat{n}} \bfu} + \dsurfintphys{b D_{\hat{n}} \bfut }{ \bfu} - \dsurfintphys{A \bfut }{ \bfu} \\
	&=  \dd{}{t} \left( \dsurfintphys{b D_{\hat{n}} \bfu }{\bfu} - \frac{1}{2} \dsurfintphys{A \bfu }{\bfu} \right).
\end{aligned}
\end{equation}
We now obtain the energy balance
\begin{equation}
	\dd{E_{D}}{t} = 0,
\end{equation}
where
\begin{equation} \label{eq:e_d}
	E_{D} = E - \dsurfintphys{b D_{\hat{n}} \bfu }{\bfu} + \frac{1}{2} \dsurfintphys{A \bfu }{\bfu}.
\end{equation}
To prove stability, it remains to prove that we can choose $A$ so that $E_D$ is a nonnegative quantity. Loosely speaking, we will need $A$ to be sufficiently positive. Note that $E_D$ (just like $E$) is a high-order approximation of the continuous energy $\mathcal{E}$ because the boundary integrals are zero to the order of accuracy due to the boundary condition.

To prove that $E_D$ is nonnegative, we need to use the positivity of $E$. Since the indefinite term in \eqref{eq:e_d} is a surface integral, we shall bound $E$ from below by relevant surface integrals. The terms in $E$ that stem from the Laplacian are
\begin{equation}
 	\frac{1}{2} \underbrace{ \si \dnormphys{\sqrt{b} D_{x_i} \bfu}^2 }_{\text{Term } 1}  \quad \mbox{and} \quad \frac{1}{2} \underbrace{ \sk \bfu^T R_{kk}(\alpha_{ikik}  b) \bfu }_{\text{Term } 2} .
 \end{equation}
Applying the positivity property \eqref{eq:borrowing_nd_H_global} to Term 1 yields
\begin{equation} \label{eq:borrowing_H_multiD}
\begin{aligned}
	\si \dnormphys{\sqrt{b} D_{x_i} \bfu}^2 & = \si \dnorm{\sqrt{Jb}D_{x_i} \bfu} \\
	&\geq \frac{h\theta_H}{d} \si \dbnorm{\sqrt{Jb}D_{x_i} \bfu} \\
	&= \frac{h\theta_H}{d} \si \dbnormphys{\sqrt{\frac{Jb}{\gamma} }D_{x_i} \bfu}.
\end{aligned}
\end{equation}
Now consider Term 2. To simplify the notation, we define
\begin{equation}
    \eta_k = \alpha_{ikik} b \quad \mbox{(no sum over $k$)}.
\end{equation}
Applying the positivity property \eqref{eq:borrowing_nd_R} yields
\begin{equation} \label{eq:borrowing_R_multiD}
\begin{aligned}
	\sk \bfu^T R_{kk}(\eta_k) \bfu & \geq h\theta_R \sk  \dbnorm{\nu_k \left( \sqrt{\eta_k} \right)_{min} \Delta D_k \bfu} \\
	&= h\theta_R \sk  \dbnormphys{\nu_k \frac{\left( \sqrt{\eta_k} \right)_{min}}{\sqrt{\gamma}} \Delta D_k \bfu} .
\end{aligned}
\end{equation}
Using \eqref{eq:borrowing_H_multiD} and \eqref{eq:borrowing_R_multiD}, we obtain
\begin{equation}
\begin{aligned}
2E_D &\geq \frac{h\theta_H}{d} \si \dbnormphys{\sqrt{\frac{Jb}{\gamma} } D_{x_i} \bfu} +h\theta_R \sk  \dbnormphys{\nu_k \frac{\left( \sqrt{\eta_k} \right)_{min}}{\sqrt{\gamma}} \Delta D_k \bfu} \\
&- 2\dsurfintphys{b D_{\hat{n}} \bfu }{\bfu} + \dsurfintphys{A \bfu }{\bfu} .
\end{aligned}
\end{equation}
Expanding the indefinite term $\dsurfintphys{b D_{\hat{n}} \bfu }{\bfu}$ yields
\begin{equation} \label{eq:indefinite_term_expanded}
\begin{aligned}
	&\dsurfintphys{b D_{\hat{n}} \bfu }{\bfu} = \\
    &\dsurfintphys{b \left(\nu_j \frac{\alpha_{ijik}}{\gamma} D_{k} - \sk \nu_k \frac{\alpha_{ikik}}{\gamma} \Delta D_{k} \right) \bfu}{\bfu} =\\
	&\dsurfintphys{b \nu_j \frac{\K_{ij} J \K_{ik}}{\gamma} D_{k} \bfu}{\bfu} - \dsurfintphys{\sk \nu_k \frac{\eta_k}{\gamma} \Delta D_{k} \bfu}{\bfu} =\\
	&\dsurfintphys{ \frac{bJ}{\gamma}  D_{x_i} \bfu}{\nu_j \K_{ij} \bfu} - \sk \dsurfintphys{\nu_k \frac{\eta_k}{\gamma} \Delta D_{k} \bfu}{ \abs{\nu_k} \bfu}.
\end{aligned}
\end{equation}
We obtain
\begin{equation}
    2E_D \geq \dsurfintphys{A \bfu }{\bfu} + \dsurfintphys{y_i}{\widetilde{Y} y_i} + \sk \dsurfintphys{z_{k}}{\widetilde{Z}_k z_{k}},
\end{equation}
where
\begin{equation}
    y_i = \begin{bmatrix} \nu_k \K_{ik} \bfu \\
      D_{x_i} \bfu \\
     \end{bmatrix},
     \quad
    \widetilde{Y} = \begin{bmatrix}
         0 & -1 \\
         -1 & \frac{h\theta_H}{d}  \\
    \end{bmatrix} \otimes \frac{bJ}{\gamma} ,
\end{equation}
\begin{equation}
    z_{k} = \begin{bmatrix} \abs{\nu_k} \bfu \\
      \nu_k \Delta D_k \bfu
     \end{bmatrix}
    \; \mbox{(no sum over $k$)} ,
\end{equation}
and
\begin{equation}
    \widetilde{Z}_k = \begin{bmatrix}
        \medskip
        \displaystyle
         0 &\displaystyle -\frac{\eta_k}{\gamma}  \\
         \displaystyle
         -\frac{\eta_k}{\gamma} &\displaystyle h\theta_R \frac{\left( \eta_k \right)_{min}}{\gamma}
    \end{bmatrix}.
\end{equation}
We need to choose $A$ so that $\dsurfintphys{A \bfu }{\bfu}$ compensates for the zero entries in $\widetilde{Y}$ and $\widetilde{Z}_k$. We make the ansatz
\begin{equation}
    A =  \frac{1}{\gamma} \sk \abs{\nu_k} \left(\tau_H \eta_k  + \tau_R \frac{\eta_k^2}{ \left( \eta_k \right)_{min} } \right),
\end{equation}
where $\tau_{H,R}$ are penalty parameters that we have yet to choose. By noting that
\begin{equation} \label{eq:tau_H_to_sum}
\begin{aligned}
    \sk \abs{\nu_k}\eta_k &= b\sk \abs{\nu_k} \alpha_{ikik} = bJ\sk \abs{\nu_k} \K_{ik} \K_{ik} = bJ \sk \nu_k \nu_k \K_{ik} \K_{ik} \\
    &= bJ \nu_j \K_{ij} \nu_k \K_{ik},
\end{aligned}
\end{equation}
we can write
\begin{equation}
    A =   b \frac{J}{\gamma} \tau_H \nu_j \K_{ij}  \nu_k \K_{ik} + \frac{1}{\gamma} \sk \abs{\nu_k} \tau_R \frac{\eta_k^2}{ \left( \eta_k \right)_{min} } .
\end{equation}
We now obtain
\begin{equation}
	2E_D \geq  \dsurfintphys{y_i}{Y y_i} + \sk \dsurfintphys{z_{k}}{Z_k z_{k}},
\end{equation}
where
\begin{equation}
	Y = \begin{bmatrix}
		 \tau_H & -1 \\
		 -1 & \frac{h\theta_H}{d}  \\
	\end{bmatrix} \otimes \frac{bJ}{\gamma} ,
    \quad
    Z_k = \begin{bmatrix}
        \medskip
        \displaystyle
         \tau_R \frac{\eta_k^2}{ \gamma \left( \eta_k \right)_{min} } &\displaystyle -\frac{\eta_k}{\gamma}  \\
         \displaystyle
         -\frac{\eta_k}{\gamma} &\displaystyle h\theta_R \frac{\left( \eta_k \right)_{min}}{\gamma}
    \end{bmatrix}.
\end{equation}
The matrices $Y$ and $Z_k$ are positive semidefinite if
\begin{equation}
	\tau_H \geq \frac{d}{h\theta_H} \quad \mbox{and} \quad \tau_R \geq \frac{1}{h\theta_R} ,
\end{equation}
respectively.
We have now proven the following theorem.
\begin{theorem} \label{theorem:dirichlet}
The scheme
\begin{equation}
	J a \ddot{\bv{u}} = \cD_{jk}(\alpha_{ijik} b) \bv{u} + SAT,
\end{equation}
with
\begin{equation} \label{eq:displ_sat}
\begin{aligned}
    SAT &= H^{-1} \sum \limits_{f \in F} (bD_{\hat{n}})^T e_{f} \gamma H_{f} (e_{f}^T \bv{u} - 0) \\
    &- H^{-1} \sum \limits_{f \in F} \sk \frac{1}{\gamma} \abs{\nu_k} \left( \tau_H \eta_k + \tau_R \frac{\eta_k^2}{ \left( \eta_k \right)_{min} } \right) e_{f} \gamma H_{f} (e_{f}^T \bv{u} - 0),
\end{aligned}
\end{equation}
is stable if
\begin{equation}
    \tau_H \geq \frac{d}{h \theta_H}, \quad  \tau_R \geq \frac{1}{h \theta_R} .
\end{equation}
\end{theorem}
In all simulations in this paper, we set $\tau_H = d/(h \theta_H)$ and $\tau_R = 1/(h\theta_R)$, i.e., right on the limit of provable stability. The drawback to using larger values of $\tau_{H,R}$ is that it increases the spectral radius.

We remark that the fundamental difference compared to the approach of Virta and Mattsson \cite{VirtaMattsson14} is that they did not utilize Lemma \ref{lemma:R}, and hence not the positivity property \eqref{eq:borrowing_R_multiD}, which is a consequence of Lemma \ref{lemma:R}. Hence, they were forced to use a larger $SAT$ to prove that $E_D$ is nonnegative, which in turn caused the increase in spectral radius.

If the SBP operators $D_{x}$ and $D_{xx}$ are fully compatible, then $\Delta D_x=0$ and hence $\Delta D_i=0$, which implies that the scheme in Theorem \ref{theorem:dirichlet} is stable with $\tau_R = 0$.

\section{Interfaces}
Let $\physboundary_I$ denote the interface between two domains $\physdom_u$ and $\physdom_v$. We use superscripts $u$ and $v$ to denote variable coefficients that correspond to the two different sides of the interface. We consider the problem
\begin{equation} \label{eq:interface_phys}
\begin{array}{rlr}
    a^u \ddot{u} - \ddxi b^u \ddxi u = 0, & \vec{x} \in \physdom_u, & t \in [0, \; T], \\
    a^v \ddot{v} - \ddxi b^v \ddxi v = 0, & \vec{x} \in \physdom_v, & t \in [0, \; T], \\
    u - v = 0, & \vec{x} \in \physboundary_I, & t \in [0, \; T], \\
    b^u \pdd{u}{\hat{n}^u} + b^v \pdd{v}{\hat{n}^v} = 0, & \vec{x} \in \physboundary_I, & t \in [0, \; T], \\
\end{array}
\end{equation}
augmented with suitable boundary conditions. The system \eqref{eq:interface_phys} transforms into
\begin{equation} \label{eq:interface_ref}
\begin{array}{rlr}
	J^u a^u \ddot{u} - \ddj \alpha^u_{ijik} b^u \ddk u = 0, & \vec{\xi} \in \refdom_u, & t \in [0, \; T], \\
	J^v a^v \ddot{v} - \ddj \alpha^v_{ijik} b^v \ddk v = 0, & \vec{\xi} \in \refdom_v, & t \in [0, \; T], \\
	u - v = 0, & \vec{\xi} \in \refboundary_I, & t \in [0, \; T], \\
	b^u \pdd{u}{\hat{n}^u} + b^v \pdd{v}{\hat{n}^v} = 0, & \vec{\xi} \in \refboundary_I, & t \in [0, \; T]. \\
\end{array}
\end{equation}
Define the energies
\begin{equation}
    \mathcal{E}_u = \frac{1}{2} \normphysu{\sqrt{a^u}\dot{u}}^2 + \frac{1}{2} \si \normphysu{\sqrt{b^u} \ddxi u}^2,
\end{equation}
\begin{equation}
    \mathcal{E}_v = \frac{1}{2} \normphysv{\sqrt{a^v}\dot{v}}^2 + \frac{1}{2} \si \normphysv{\sqrt{b^v} \ddxi v}^2.
\end{equation}
Assuming energy-conserving boundary conditions, the energy method yields
\begin{equation}
    \dd{}{t} \left( \mathcal{E}_u + \mathcal{E}_v \right) = \surfintindexphys{\ut}{b^u \pdd{u}{\hat{n}^u}}{I} + \surfintindexphys{\dot{v}}{b^v \pdd{v}{\hat{n}^v}}{I} = 0.
\end{equation}

\begin{theorem} \label{thm:interface_stab}
The scheme
\begin{equation}
\begin{aligned}
	J^u a^u \ddot{\bv{u}} - \cD_{jk}(\alpha_{ijik}^u b^u) \bv{u} = &- H^{-1} A e_{\Gamma_I} \gamma H_{\Gamma_I} (e_{\Gamma_I}^T \bv{u} - e_{\Gamma_I}^T \bv{v} ) \\
	&+ \frac{1}{2} H^{-1} (b^u D^{u}_{\hat{n}})^T e_{\Gamma_I} \gamma H_{\Gamma_I} (e_{\Gamma_I}^T \bv{u} - e_{\Gamma_I}^T \bv{v} ) \\
	&- \frac{1}{2} H^{-1} e_{\Gamma_I} \gamma H_{\Gamma_I} (e_{\Gamma_I}^T b^u D^{u}_{\hat{n}} \bfu + e_{\Gamma_I}^T b^v D_{\hat{n}}^v \bfv ), \\
\end{aligned}
\end{equation}
\begin{equation}
\begin{aligned}
	J^v a^v \ddot{\bv{v}} - \cD_{jk}(\alpha_{ijik}^v b^v) \bv{v} = &- H^{-1} A e_{\Gamma_I} \gamma H_{\Gamma_I} (e_{\Gamma_I}^T \bv{v} - e_{\Gamma_I}^T \bv{u} ) \\
	&+ \frac{1}{2} H^{-1} (b^v D^{v}_{\hat{n}})^T e_{\Gamma_I} \gamma H_{\Gamma_I} (e_{\Gamma_I}^T \bv{v} - e_{\Gamma_I}^T \bv{u} ) \\
	&- \frac{1}{2} H^{-1} e_{\Gamma_I} \gamma H_{\Gamma_I} ( e_{\Gamma_I}^T b^v D_{\hat{n}}^v \bfv + e_{\Gamma_I}^T b^u D^{u}_{\hat{n}} \bfu), \\
\end{aligned}
\end{equation}
where
\begin{equation}
\begin{aligned}
    A &= \tau_H \sk \left( \frac{1}{\gamma^u} \abs{\nu^u_k} \eta^u_k + \frac{1}{\gamma^v} \abs{\nu^v_k} \eta^v_k\right) \\
       &+ \tau_R \sk \left( \abs{\nu^u_k} \frac{1}{\gamma^u}\frac{(\eta_k^u)^2}{ \left( \eta^u_k \right)_{min} } + \abs{\nu^v_k} \frac{1}{\gamma^v}\frac{(\eta_k^v)^2}{ \left( \eta^v_k \right)_{min} } \right),
\end{aligned}
\end{equation}
is a stable discretization of \eqref{eq:interface_ref} if
\begin{equation}
    \tau_H \geq \frac{d}{4h\theta_H}, \quad \tau_R \geq \frac{1}{4h\theta_R}.
\end{equation}
\end{theorem}

\begin{proof}
Define the discrete energies
\begin{equation}
    E_u = \frac{1}{2} \dnormphysu{\sqrt{a^u} \bfut}^2 + \frac{1}{2} \si \dnormphysu{\sqrt{b^u} D_{x_i}^u \bfu}^2 + \frac{1}{2} \sk \bfu^T R_{kk}(\alpha^u_{ikik}  b^u) \bfu
\end{equation}
and
\begin{equation}
    E_v = \frac{1}{2} \dnormphysv{\sqrt{a^v} \bfvt}^2 + \frac{1}{2} \si \dnormphysv{\sqrt{b^v} D_{x_i}^v \bfv}^2 + \frac{1}{2} \sk \bfv^T R_{kk}(\alpha^v_{ikik}  b^v) \bfv,
\end{equation}
where
\begin{equation}
    D_{x_i}^u = \K^u_{ij} D_j, \quad D_{x_i}^v = \K^v_{ij} D_j .
\end{equation}
The discrete energy method yields
\begin{equation} \label{eq:semidiscr_energy_rate_interface}
\begin{aligned}
	\dd{}{t} \left(E_u + E_v \right) &= \dsurfintphysindex{\bfut}{b^u D^{u}_{\hat{n}} \bfu}{I} + \dsurfintphysindex{\bfvt}{b^v D^{v}_{\hat{n}} \bfv}{I} \\
	&+ SAT_u + SAT_v,
\end{aligned}
\end{equation}
where
\begin{equation}
\begin{aligned}
SAT_u = &-\dsurfintphysindex{A \bfut}{\bfu - \bfv}{I} + \frac{1}{2} \dsurfintphysindex{b^u D_{\hat{n}}^u \bfut}{\bfu - \bv{v}}{I} \\
& - \frac{1}{2}  \dsurfintphysindex{\bfut}{b^u D_{\hat{n}}^u \bfu + b^v D_{\hat{n}}^v \bfv}{I} \\
= &-  \dsurfintphysindex{A \bfut}{\bfu - \bfv}{I} + \frac{1}{2} \dsurfintphysindex{b^u D_{\hat{n}}^u \bfut }{\bfu - \bv{v}}{I} \\
&- \frac{1}{2}  \dsurfintphysindex{\bfut}{b^u D_{\hat{n}}^u \bfu + b^v D_{\hat{n}}^v \bfv}{I}
\end{aligned}
\end{equation}
and
\begin{equation}
\begin{aligned}
SAT_v = &-  \dsurfintphysindex{A \bfvt}{\bv{v} - \bfu}{I} + \frac{1}{2}  \dsurfintphysindex{b^v D_{\hat{n}}^v \bfvt}{\bv{v} - \bfu}{I} \\
&- \frac{1}{2}  \dsurfintphysindex{\bfvt}{b^v D_{\hat{n}}^v \bfv + b^u D_{\hat{n}}^u \bfu}{I}	\,.
\end{aligned}
\end{equation}
Moving all terms in \eqref{eq:semidiscr_energy_rate_interface} to the left-hand side yields
\begin{equation}
	\dd{}{t} E_I = 0,
\end{equation}
where
\begin{equation} \label{eq:energy_interface}
\begin{aligned}
	E_I &= E_u + E_v + \frac{1}{2} \Bigg[ \dbnormphysindex{\sqrt{A}(\bfu - \bfv)}{I} - \dsurfintphysindex{b^u D_{\hat{n}}^u \bfu}{\bv{u}}{I} \\
	&- \dsurfintphysindex{b^v D_{\hat{n}}^v \bfv}{\bv{v}}{I}
	+ \dsurfintphysindex{b^u D_{\hat{n}}^u \bfu}{\bv{v}}{I} + \dsurfintphysindex{b^v D_{\hat{n}}^v \bfv}{\bv{u}}{I}  \Bigg] \,.
\end{aligned}
\end{equation}
Note that the discrete energy $E_I$ approximates the continuous energy $\mathcal{E}_u$ + $\mathcal{E}_v$ because the interface terms in \eqref{eq:energy_interface} would be zero if the interface conditions were fulfilled exactly.
It remains to show that $E_I \geq 0$. We start by expanding the indefinite terms. By \eqref{eq:indefinite_term_expanded} we have
\begin{equation}
    \dsurfintphysindex{b D_{\hat{n}} \bfu }{\bfu}{I} = \dsurfintphysindex{\beta D_{x_i} \bfu}{\nu_j \K_{ij} \bfu}{I} - \sk \dsurfintphysindex{\frac{1}{\gamma} \eta_k \Delta D_{k} \bfu}{ \nu_k \bfu}{I} ,
\end{equation}
where we have defined $\beta = b\frac{J}{\gamma}$.
Similarly, we have
\begin{equation*}
    \dsurfintphysindex{b^u D^u_{\hat{n}} \bfu }{\bfu}{I} = \dsurfintphysindex{ \beta^u  D_{x_i}^u \bfu}{\nu^u_j \K_{ij}^u \bfu}{I} - \sk \dsurfintphysindex{\frac{1}{\gamma^u} \eta^u_k \Delta D_{k} \bfu}{ \nu^u_k \bfu}{I} ,
\end{equation*}
\begin{equation*}
    \dsurfintphysindex{b^v D^v_{\hat{n}} \bfv }{\bfv}{I} = \dsurfintphysindex{\beta^v D_{x_i}^v \bfv}{\nu^v_j \K_{ij}^v \bfv}{I} - \sk \dsurfintphysindex{\frac{1}{\gamma^v} \eta^v_k \Delta D_{k} \bfv}{ \nu^v_k \bfv}{I} ,
\end{equation*}
\begin{equation*}
    \dsurfintphysindex{b^u D^u_{\hat{n}} \bfu }{\bfv}{I} = \dsurfintphysindex{\beta^u D_{x_i}^u \bfu}{\nu^u_j \K_{ij}^u \bfv}{I} - \sk \dsurfintphysindex{\frac{1}{\gamma^u} \eta^u_k \Delta D_{k} \bfu}{ \nu^u_k \bfv}{I} ,
\end{equation*}
and
\begin{equation*}
    \dsurfintphysindex{b^v D^v_{\hat{n}} \bfv }{\bfu}{I} = \dsurfintphysindex{\beta^v D_{x_i}^v \bfv}{\nu^v_j \K_{ij}^v \bfu}{I} - \sk \dsurfintphysindex{\frac{1}{\gamma^v} \eta^v_k \Delta D_{k} \bfv}{ \nu^v_k \bfu}{I} .
\end{equation*}
The positivity properties \eqref{eq:borrowing_H_multiD} and \eqref{eq:borrowing_R_multiD} yield
\begin{equation}
\begin{aligned}
2E_u &\geq \frac{h\theta_H}{d} \si \dbnormphysindex{\sqrt{\beta^u } D_{x_i}^u \bfu}{I} \\
&+h\theta_R \sk  \dbnormphysindex{\nu^u_k \frac{\left( \sqrt{\eta_k^u} \right)_{min}}{\sqrt{\gamma^u}} \Delta D_k \bfu}{I} ,
\end{aligned}
\end{equation}
and
\begin{equation}
\begin{aligned}
2E_v &\geq \frac{h\theta_H}{d} \si \dbnormphysindex{\sqrt{\beta^v } D_{x_i}^v \bfv}{I} \\
&+ h\theta_R \sk  \dbnormphysindex{\nu^v_k \frac{\left( \sqrt{\eta_k^v} \right)_{min}}{\sqrt{\gamma^v}} \Delta D_k \bfv}{I} .
\end{aligned}
\end{equation}
Using \eqref{eq:tau_H_to_sum}, we can show that, for example,
\begin{equation}
    \sk \frac{1}{\gamma^u} \abs{\nu^u_k} \eta^u_k = \frac{b_u J_u}{\gamma^u} \nu^u_j \K_{ij}^u \nu^u_k \K_{ik}^u = \beta^u \nu^u_j \K_{ij}^u \nu^u_k \K_{ik}^u,
\end{equation}
and hence
\begin{equation}
\begin{aligned}
A &= \tau_H \left( \beta^u \nu^u_j \K^u_{ij} \nu^u_k \K^u_{ik}  + \beta^v \nu^v_j \K^v_{ij} \nu^v_k \K^v_{ik} \right) \\
       &+ \tau_R \sk \left( \frac{\abs{\nu^u_k}}{\gamma^u}\frac{(\eta_k^u)^2}{ \left( \eta^u_k \right)_{min} } + \frac{\abs{\nu^v_k}}{\gamma^v}\frac{(\eta_k^v)^2}{ \left( \eta^v_k \right)_{min} } \right).
   \end{aligned}
\end{equation}
Define
\begin{equation}
    y_i = \begin{bmatrix}
     \sqrt{\beta^u} \nu^u_j \K^u_{ij} \bfu \\
     \sqrt{\beta^v} \nu^v_j \K^v_{ij} \bfu \\
     \sqrt{\beta^v} \nu^v_j \K^v_{ij} \bfv \\
     \sqrt{\beta^u} \nu^u_j \K^u_{ij} \bfv \\
     \sqrt{\beta^u} D_{x_i}^u \bfu \\
     \sqrt{\beta^v} D_{x_i}^v \bfv \\
    \end{bmatrix},
    \quad
    z_{k} = \begin{bmatrix}
    \abs{\nu^u_k} \bfu \\
    \abs{\nu^v_k} \bfu \\
    \abs{\nu^v_k} \bfv \\
    \abs{\nu^u_k} \bfv \\
    \nu^u_k \Delta D_k \bfu \\
    \nu^v_k \Delta D_k \bfv \\
    \end{bmatrix} \; \mbox{(no sum over $k$)}.
\end{equation}
After some algebra we obtain
\begin{equation}
	2E_I \geq \dsurfintphysindex{y_i}{Y y_i}{I} +  \sk \dsurfintphysindex{z_{k}}{Z_k z_{k}}{I},
\end{equation}
where
\begin{equation}
    Y = \begin{bmatrix}
\smallskip
\tau_H & 0 & 0 & -\tau_H &  -\frac{1}{2} & 0 \\
\smallskip
0 & \tau_H & -\tau_H & 0 & 0 &  \frac{1}{2} \\
\smallskip
0 & -\tau_H & \tau_H & 0 & 0 & -\frac{1}{2} \\
\smallskip
-\tau_H & 0 & 0 & \tau_H &  \frac{1}{2} & 0\\
\smallskip
 -\frac{1}{2} & 0 & 0 &  \frac{1}{2} &  \frac{h\theta_H}{d} & 0 \\
0 &  \frac{1}{2} &  -\frac{1}{2} & 0 & 0 &  \frac{h\theta_H}{d} \\
\end{bmatrix} \otimes I,
\end{equation}
and
\begin{equation*}
Z_k = \begin{bmatrix}
\smallskip
 \frac{\tau_R}{\gamma^u} \frac{\left(\eta_k^u\right)^2}{\left(\eta_k^u\right)_{min}} & 0 & 0 &  -\frac{\tau_R}{\gamma^u} \frac{\left(\eta_k^u\right)^2}{\left(\eta_k^u\right)_{min}} &  \frac{\eta_k^u}{2\gamma^u}  & 0 \\
\smallskip
0 &  \frac{\tau_R}{\gamma^v} \frac{\left(\eta_k^v\right)^2}{\left(\eta_k^v\right)_{min}} &  -\frac{\tau_R}{\gamma^v} \frac{\left(\eta_k^v\right)^2}{\left(\eta_k^v\right)_{min}} & 0 & 0 & - \frac{\eta_k^v}{2\gamma^v} \\
\smallskip
0 & -\frac{\tau_R}{\gamma^v} \frac{\left(\eta_k^v\right)^2}{\left(\eta_k^v\right)_{min}} &  \frac{\tau_R}{\gamma^v} \frac{\left(\eta_k^v\right)^2}{\left(\eta_k^v\right)_{min}} & 0 & 0 &  \frac{\eta_k^v}{2\gamma^v} \\
\smallskip
- \frac{\tau_R}{\gamma^u} \frac{\left(\eta_k^u\right)^2}{\left(\eta_k^u\right)_{min}} & 0 & 0 &  \frac{\tau_R}{\gamma^u} \frac{\left(\eta_k^u\right)^2}{\left(\eta_k^u\right)_{min}} & - \frac{\eta_k^u}{2\gamma^u} & 0\\
\smallskip
 \frac{\eta_k^u}{2\gamma^u} & 0 & 0 & - \frac{\eta_k^u}{2\gamma^u} & \frac{h\theta_R \left( \eta_k^u \right)_{min}}{\gamma^u} & 0 \\
\smallskip
0 & - \frac{\eta_k^v}{2\gamma^v} &  \frac{\eta_k^v}{2\gamma^v} & 0 & 0 & \frac{h\theta_R \left( \eta_k^v \right)_{min}}{\gamma^v} \\
\end{bmatrix} .
\end{equation*}
The matrices $Y$ and $Z_k$ are positive semidefinite if $\tau_H \geq \frac{d}{4 \theta_H}$ and $\tau_R \geq \frac{1}{4 \theta_R}$, respectively.
\end{proof}

In all experiments with interfaces in this paper, we set $\tau_H = d/(4 \theta_H)$ and $\tau_R = 1/(4 \theta_R)$, i.e., right on the limit of provable stability.

\section{Applications} \label{sec:applications}

In this section we will use the new method to solve two application problems. To choose the time-step without having to compute eigenvalues of the spatial discretization we use the procedure described below.

Consider first a single-block grid. Let $c = \sqrt{b/a}$ denote the wave speed. Let $h_i$ denote the grid spacing in the $\xi_i$-direction in the reference domain. We compute approximate physical grid spacings $\hphys_j$ as
\begin{equation}
    \hphys_j = \sqrt{\si \left( \mathcal{J}_{ij} h_j \right)^2 } \quad \mbox{no sum over }j,
\end{equation}
and define
\begin{equation}
    \hphys = \min \limits_j \hphys_j.
\end{equation}
Note that both $\hphys$ and $c$ are non-constant on the grid, in general. Let $\mathcal{G}$ denote the set of all grid points. We compute the time step $\Delta t$ as
\begin{equation}
    \left. \Delta t = CFL \times \min \limits_{\vec{y} \in \mathcal{G} } \frac{\hphys}{c} \right|_{\vec{x}=\vec{y}} ,
\end{equation}
where $CFL$ is a dimensionless constant that depends on the spatial order of accuracy and the choice of time-integrator. For a multiblock grid, it only remains to loop over all grid blocks and choose the smallest $\Delta t$ suggested by any block.

For the example applications in this section we use the classical fourth order Runge--Kutta method with the $CFL$ values listed in Table \ref{table:cfl}. The numbers have been tuned to yield time-steps that are slightly smaller than the largest stable time-step, for the grid configurations that we have tried. While we cannot guarantee that these $CFL$ values yield stability for every possible grid, they should at least provide a reasonable starting point.
\begin{table}[H]
\centering
\begin{tabular}{c|ccc}
& second order & fourth order & sixth order \\
\hline
$CFL$ & $0.5$ & $0.35$ & $0.27$  \\
\end{tabular}
\caption{CFL numbers used in the classical fourth order Runge--Kutta time-integrator in the application problems, for different orders of spatial accuracy. }
\label{table:cfl}
\end{table}

\renewcommand{\ddi}{\partial_{x_i}}
\subsection{Marine seismic exploration with ocean-bottom nodes}
\begin{figure}[h]
        \centering
        \includegraphics[width= 0.9 \linewidth]{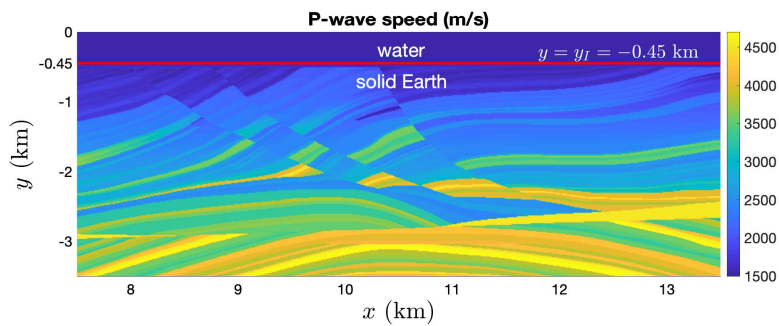}
        \caption{P-wave speed in the Marmousi2 velocity model. The horizontal red line marks the fluid-solid interface at the seafloor.}
    \label{fig:vp}
\end{figure}

The Marmousi2 velocity model \cite{Martin2006} is a fully elastic velocity model. It is an extension of the acoustic Marmousi model \cite{Versteeg1994,Versteeg1991}, which has become a standard benchmark problem for seismic imaging and full waveform inversion. Figure \ref{fig:vp} shows the P-wave speed $c_p$ in Marmousi2. The red line at $y=y_I=-0.45$ km marks the seafloor. Above the red line there is water, and below the red line there is solid Earth.
In this paper, we shall approximate the response of both ocean and solid as acoustic. The acoustic equations can be written on first order form as
\begin{equation} \label{eq:acoustics_first_order}
    \left\{ \begin{array}{rr}
    \displaystyle
    \medskip
    \rho \pdd{v_i}{t} + \ddi p = F_i \\
    \displaystyle
    K^{-1} \pdd{p}{t} + \ddi v_i = G
    \end{array} \right. ,
\end{equation}
where $p$ is the pressure perturbation, $v_i$ is the particle velocity vector, $\rho$ is density, $K$ is bulk modulus, and $F_i$ and $G$ are source terms. The Marmousi2 dataset provides density $\rho$ and P-wave speed $c_p$.
We set the bulk modulus to be $K = \rho c_p^2$, so that the wave speed in \eqref{eq:acoustics_first_order} matches the P-wave speed in the elastic model. Both $\rho$ and $K$ are spatially variable within the solid. The system \eqref{eq:acoustics_first_order} may be written as a second-order equation for $p$ only:
\begin{equation} \label{eq:p_scalar}
    K^{-1} \pddn{2}{p}{t} - \ddi \rho^{-1} \ddi p = \pdd{G}{t} - \ddi \rho^{-1} F_i .
\end{equation}
In seismic exploration, one attempts to solve optimization problems of the form \cite{Gauthier1986,Virieux2009}
\begin{equation}
    \min \limits_{\rho,K} \, M,
\end{equation}
where
\begin{equation}
\begin{aligned}
    M &= \frac{1}{2} \si \myint{0}{T}{w_{v,i}(t)\abs{v_i(\vec{x}_v,t) - v_{i,data}(t)}^2}{t}\\
    &+ \frac{1}{2} \myint{0}{T}{w_p(t)\abs{p(\vec{x}_p,t) - p_{data}(t)}^2}{t},
\end{aligned}
\end{equation}
and $v_i$ and $p$ satisfy \eqref{eq:acoustics_first_order} and appropriate boundary, interface, and initial conditions. The data $v_{i,data}$ and $p_{data}$ are measured by geophones and hydrophones, respectively, at the locations $\vec{x}_v$ and $\vec{x}_p$. The weight functions $w_{v,i}$ and $w_p$ may be used to window the signal in time. In practice one often regularizes the optimization problem by minimizing $M+S$ rather than $M$, where $S$ is a regularization term. In this paper, we set $S=0$ for simplicity.

We introduce the residuals
\begin{equation}
\begin{aligned}
    \widetilde{v}_i(t) &= w_{v,i}(t) \left(v_{i}(\vec{x}_v,t) - v_{i,data}(t)\right), \\
    \widetilde{p}(t) &= w_p(t) \left(p(\vec{x}_p,t) - p_{data}(t) \right) .
\end{aligned}
\end{equation}
As part of computing the gradient of $M$ with respect to the material parameters, one may solve the adjoint equations \cite{Gauthier1986,Virieux2009} for $v_i^{\dagger}$ and $p^{\dagger}$,
\begin{equation}\label{eq:acoustics_first_order_adj}
    \left\{ \begin{array}{rlr}
    \displaystyle
    \medskip
    \rho \pdd{v^{\dagger}_i}{\tau} + \ddi p^{\dagger} & =  &-\widetilde{v}_i(\tau) \delta(\vec{x} - \vec{x}_v) \\
    \displaystyle
    K^{-1} \pdd{p^{\dagger}}{\tau} + \ddi v_i^{\dagger} & = & \widetilde{p}(\tau) \delta(\vec{x} - \vec{x}_p)
    \end{array} \right. ,
\end{equation}
where $\tau = T - t$ is reverse time. Notice that the residuals $\widetilde{v}_i$ and $\widetilde{p}$ appear as singular source terms at the receiver locations. The second-order equation for $p^{\dagger}$ only is
\begin{equation} \label{eq:p_adj_scalar}
\begin{aligned}
    \frac{1}{K} \pddn{2}{{p}^{\dagger}}{\tau}  - \ddi \rho^{-1} \ddi p^{\dagger}
    &= \pdd{\widetilde{p}}{\tau}(\tau) \delta(\vec{x} - \vec{x}_p) + \widetilde{v}_i(\tau) \ddi \rho^{-1} \delta(\vec{x} - \vec{x}_v).
\end{aligned}
\end{equation}
Notice that \eqref{eq:p_adj_scalar} is of the form \eqref{eq:wave_eq_general}, with added source terms. We use the technique developed in this paper to discretize the operator $\ddi \rho^{-1} \ddi$. Pressure and the vertical velocity component are continuous across the fluid-solid interface. The horizontal velocity component, however, may be discontinuous. To avoid differentiating across the discontinuity, we use two grid blocks: one for the water and one for the solid Earth. Expressed in terms of $p^{\dagger}$, the interface conditions are
\begin{equation} \label{eq:interface_conditions}
    \begin{array}{rll}
    \displaystyle
    \medskip
    p^{\dagger}_{w} - p^{\dagger}_{s} & = 0, & y = y_I, \\
    \displaystyle
    \rho^{-1}_{w} \pdd{p^{\dagger}_{w}}{\hat{n}_{w}} + \rho^{-1}_{s} \pdd{p^{\dagger}_{s}}{\hat{n}_{s}} &= 0, & y = y_I ,
    \end{array}
\end{equation}
where subscripts $w$ and $s$ denote water and solid sides, respectively.
The interface conditions are of the same form as those considered in \eqref{eq:interface_phys} and we may use the interface discretization presented in Theorem \ref{thm:interface_stab}.

References \cite{Petersson2016,Sjogreen2014} describe how to discretize $\delta$ and $\ddi \delta$. When using $q$th order finite difference stencils in the interior, the discretizations of $\delta$ and $\ddi \delta$ satisfy $q$ and $q+1$ moment conditions, respectively. The distribution $\ddi \delta$ satisfies
\begin{equation}
    \ip{u_i}{\ddi \delta(\vec{x} - \vec{x}_v)}{\physdom} = -\ddi u_i |_{\vec{x} = \vec{x}_v},
\end{equation}
for all smooth fields $ u_i$. Let $ \bfu_i$ denote the restriction of $ u_i$ to the grid.
A discrete approximation $\mathbf{d}_{x_i,\vec{x}_v}' \approx \ddi \delta(\vec{x} - \vec{x}_v) $ is required to satisfy the moment conditions \cite{Sjogreen2014}
\begin{equation}
    \dvolintphys{\bfu_i}{\mathbf{d}_{x_i,\vec{x}_v}'} = -\ddi u_i |_{\vec{x} = \vec{x}_v},
\end{equation}
for all $u_i$ that are polynomials of degree less than or equal to $q$. To solve \eqref{eq:p_adj_scalar}, we actually need to consider moment conditions for $\ddi \rho^{-1} \delta(\vec{x} - \vec{x}_v)$. Integration by parts yields
\begin{equation}
    \ip{u_i}{\ddi \rho^{-1} \delta(\vec{x} - \vec{x}_v)}{\physdom} = -\rho^{-1} \ddi u_i |_{\vec{x} = \vec{x}_v}.
\end{equation}
The moment conditions on $\mathbf{d}_{x_i,\vec{x}_v, \rho^{-1}}' \approx \ddi \rho^{-1} \delta(\vec{x} - \vec{x}_v)$ are thus
\begin{equation}
    \dvolintphys{\bfu_i}{\mathbf{d}_{x_i,\vec{x}_v, \rho^{-1}}'} = -\rho^{-1} \ddi u_i |_{\vec{x} = \vec{x}_v},
\end{equation}
for all $u_i$ that are polynomials of degree less than or equal to $q$. Given $\mathbf{d}_{x_i,\vec{x}_v}'$ that satisfies the appropriate moment conditions, we note that we may set
\begin{equation}
    \mathbf{d}_{x_i,\vec{x}_v, \rho^{-1}}' = \rho^{-1}(\vec{x}_v) \mathbf{d}_{x_i,\vec{x}_v}',
\end{equation}
and satisfy the required moment conditions on $\mathbf{d}_{x_i,\vec{x}_v, \rho^{-1}}'$. Because we solve the equations on second order form, no smoothness conditions are required.

The far-field boundaries are treated with super-grid-scale absorbing layers \cite{Appelo09_2,petersson_sjogreen_2014}. At the water surface, we impose the Dirichlet condition $p=0$.
We use the classical fourth order Runge--Kutta method for time integration.

We attempt to set up simulations that are representative of the adjoint equations in marine seismic exploration. Let $W$ denote the Ricker wavelet \cite{Ricker1943, Ricker1944} with peak frequency $f$ centered at time $t_0$, i.e.,
\begin{equation}
    W(t) = (1-2 \pi^2 f^2 (t-t_0)^2 ) e^{-\pi^2 f^2 (t-t_0)^2}.
\end{equation}
We choose $f = 5$ Hz and $t_0 = 0.2$ s, and set
\begin{equation}
\widetilde{v}_i(\tau) = A_v W(\tau), \, i = 1,2 \quad \mbox{and} \quad \dd{\widetilde{p}}{\tau}(\tau) = A_p W(\tau).
\end{equation}
 To illustrate the difference between the two different adjoint sources, we run two separate simulations with only one active source in each. We place the sources at the same location, but on opposite sides of the water-solid interface. That is, we set $\vec{x}_v = \begin{bmatrix} 10.5 \mbox{ km}, &y_I- \end{bmatrix}$ and $\vec{x}_p = \begin{bmatrix} 10.5 \mbox{ km}, & y_I+ \end{bmatrix}$.

\begin{figure}[H]
\centering
    \begin{subfigure}[b]{.45\linewidth}
        \centering
        \includegraphics[width=1\linewidth,trim={0.8cm 1cm 1.5cm 1.5cm},clip]{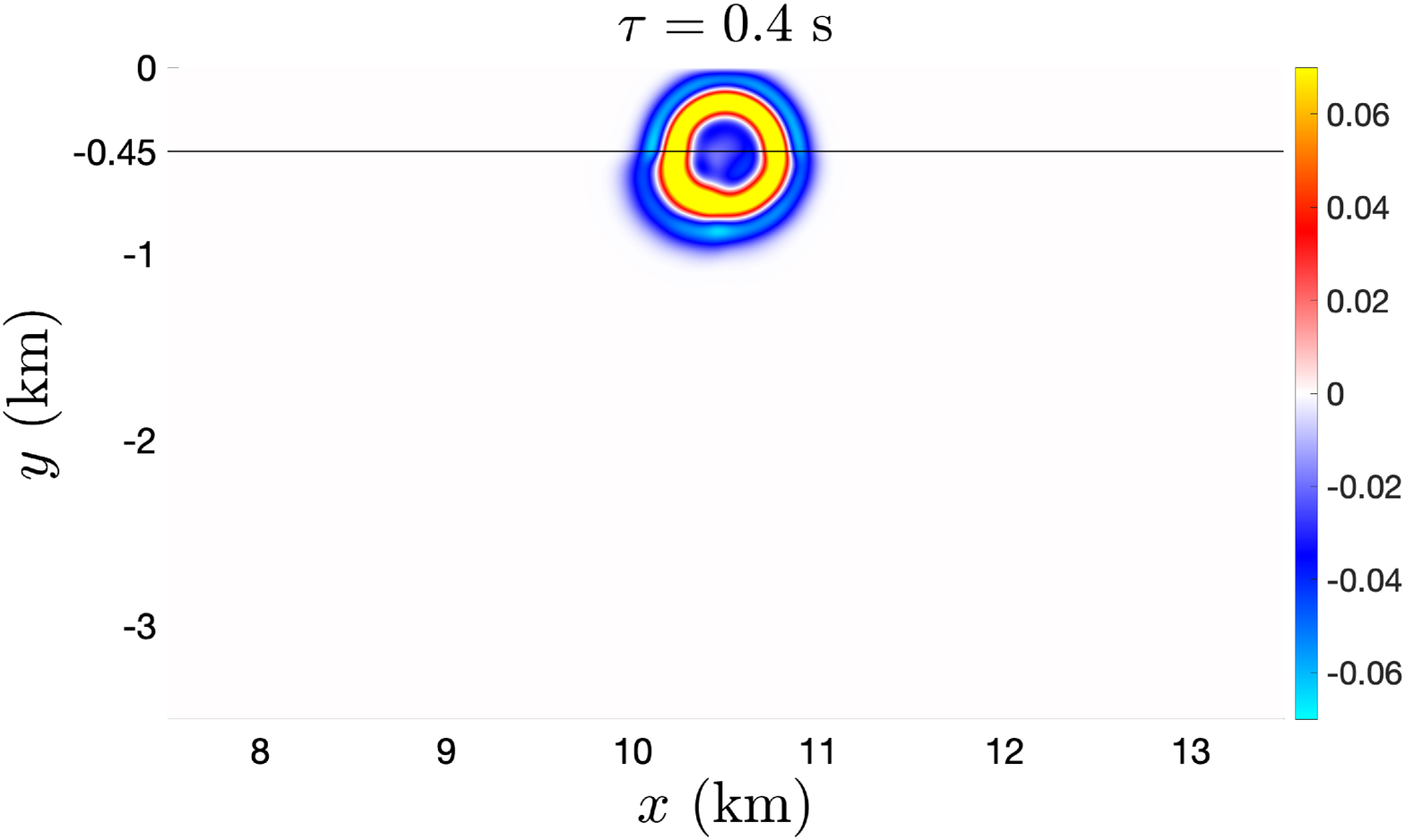}
    \end{subfigure}
     \begin{subfigure}[b]{.45\linewidth}
        \centering
        \includegraphics[width=1\linewidth,trim={0.8cm 1cm 1.5cm 1.5cm},clip]{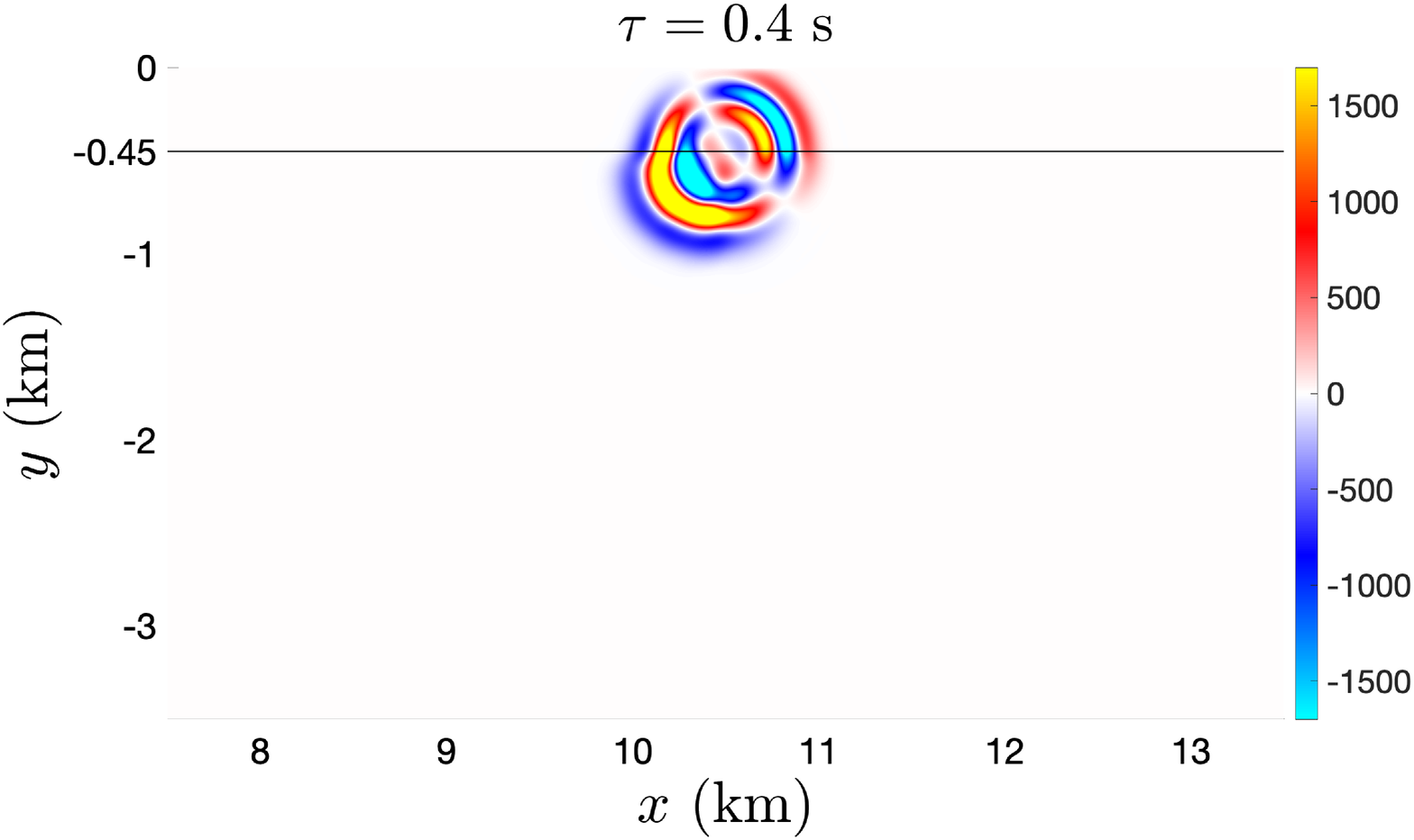}
    \end{subfigure}
    \\
    \begin{subfigure}[b]{.45\linewidth}
        \centering
        \includegraphics[width=1\linewidth,trim={0.8cm 1cm 1.5cm 1.5cm},clip]{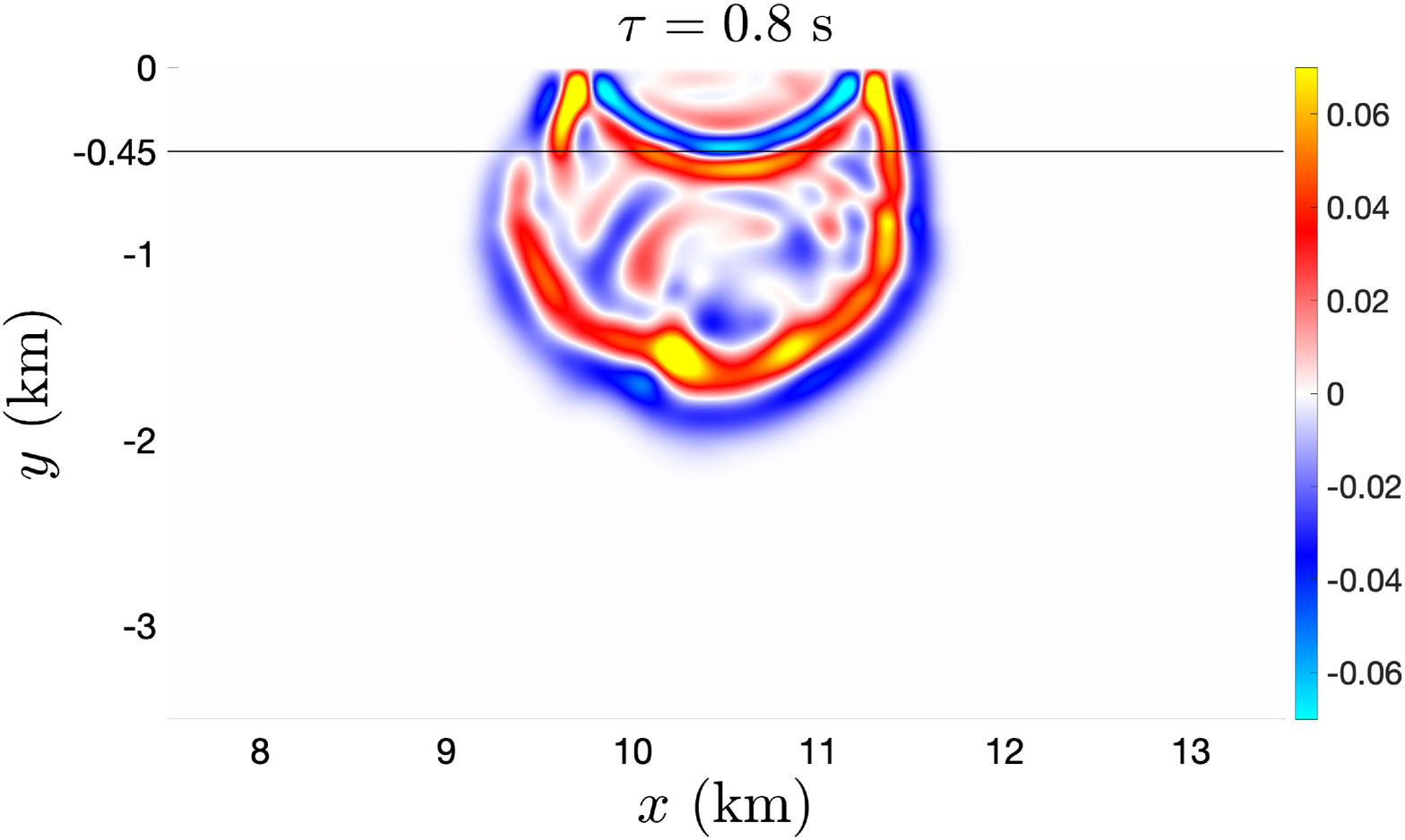}
    \end{subfigure}
    \begin{subfigure}[b]{.45\linewidth}
        \centering
        \includegraphics[width=1\linewidth,trim={0.8cm 1cm 1.5cm 1.5cm},clip]{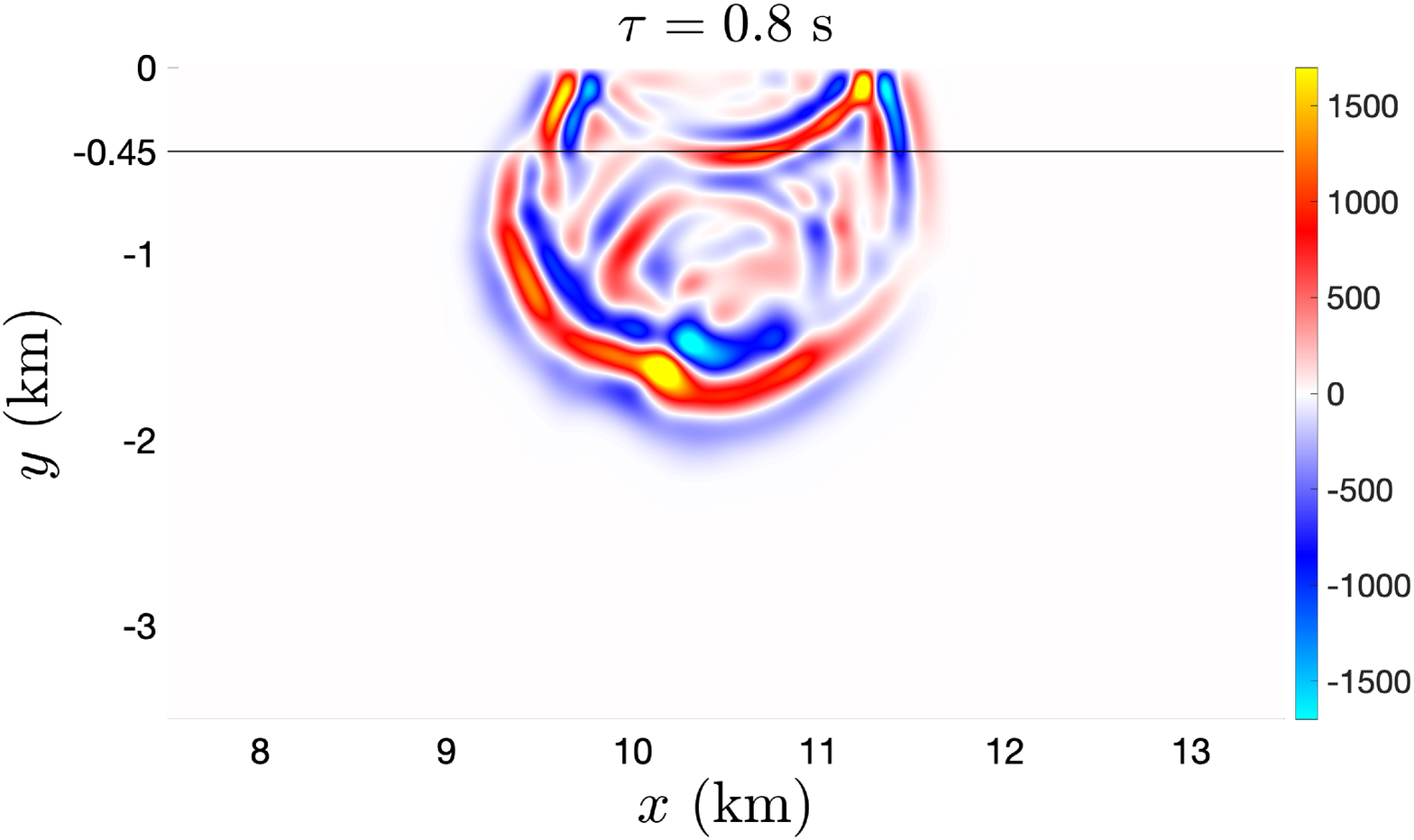}
    \end{subfigure} \\
    \begin{subfigure}[b]{.45\linewidth}
        \centering
        \includegraphics[width=1\linewidth,trim={0.8cm 1cm 1.5cm 1.5cm},clip]{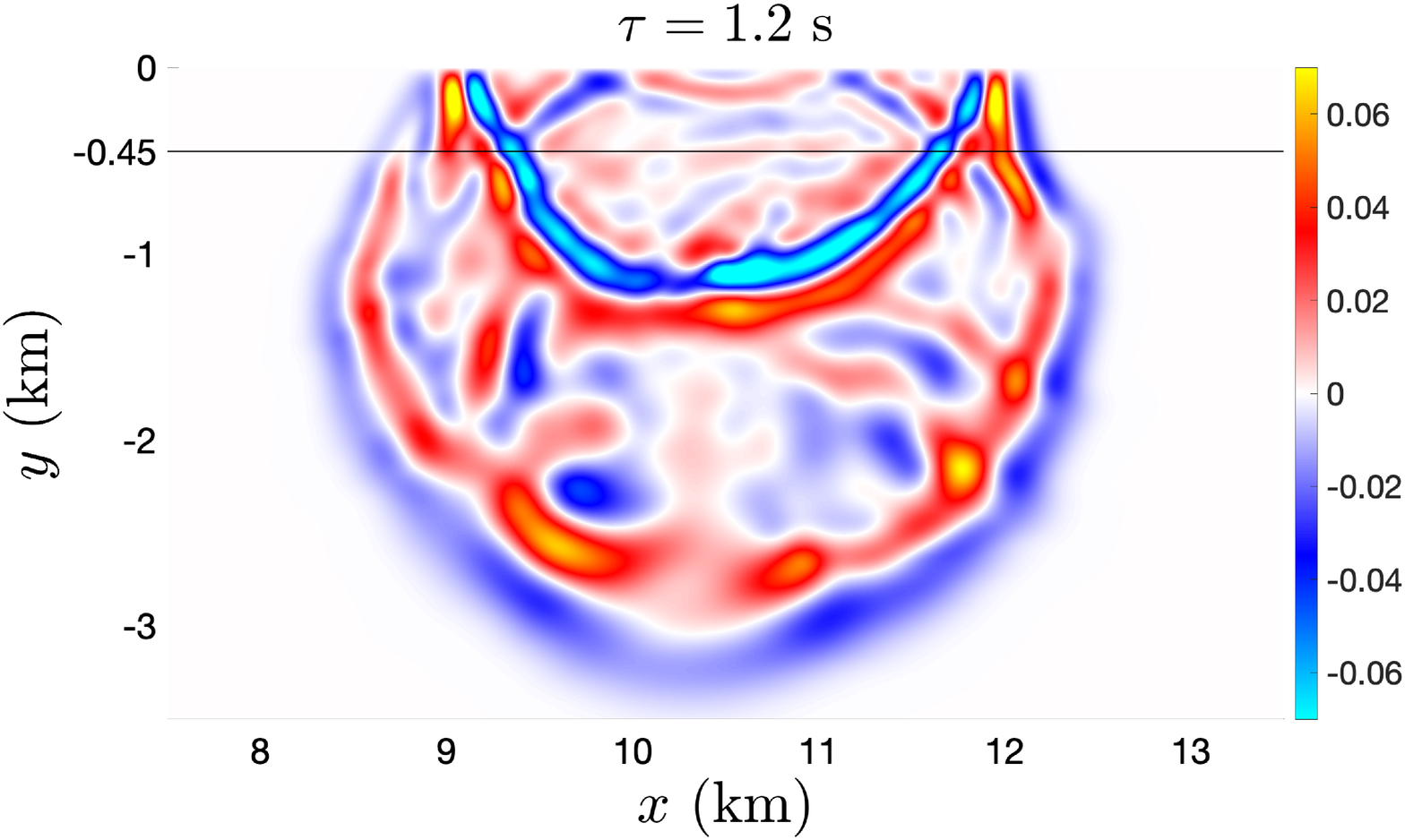}
    \end{subfigure}
    \begin{subfigure}[b]{.45\linewidth}
        \centering
        \includegraphics[width=1\linewidth,trim={0.8cm 1cm 1.5cm 1.5cm},clip]{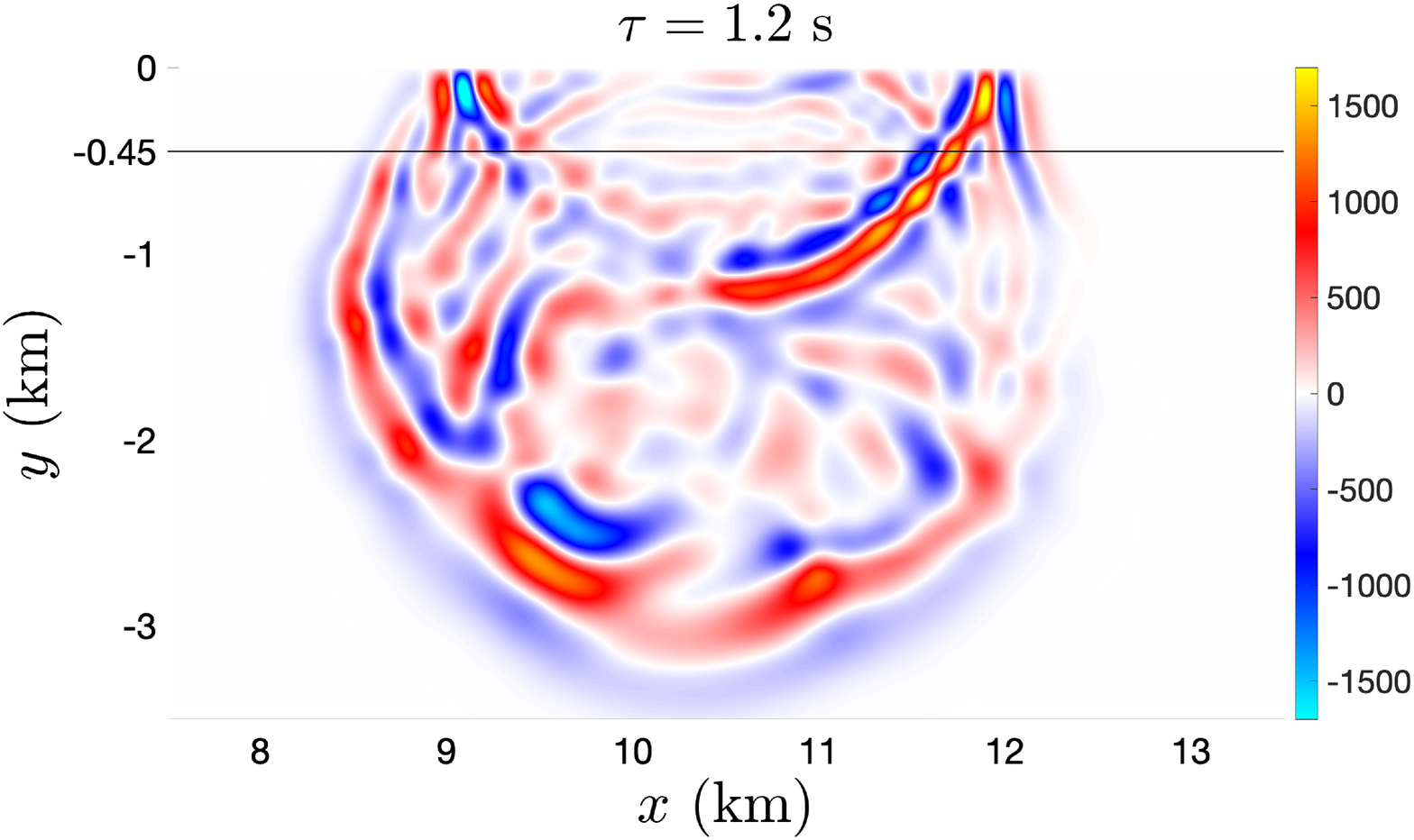}
    \end{subfigure} \\
    \begin{subfigure}[b]{.45\linewidth}
        \centering
        \includegraphics[width=1\linewidth,trim={1cm 0cm 1.5cm 0.5cm},clip]{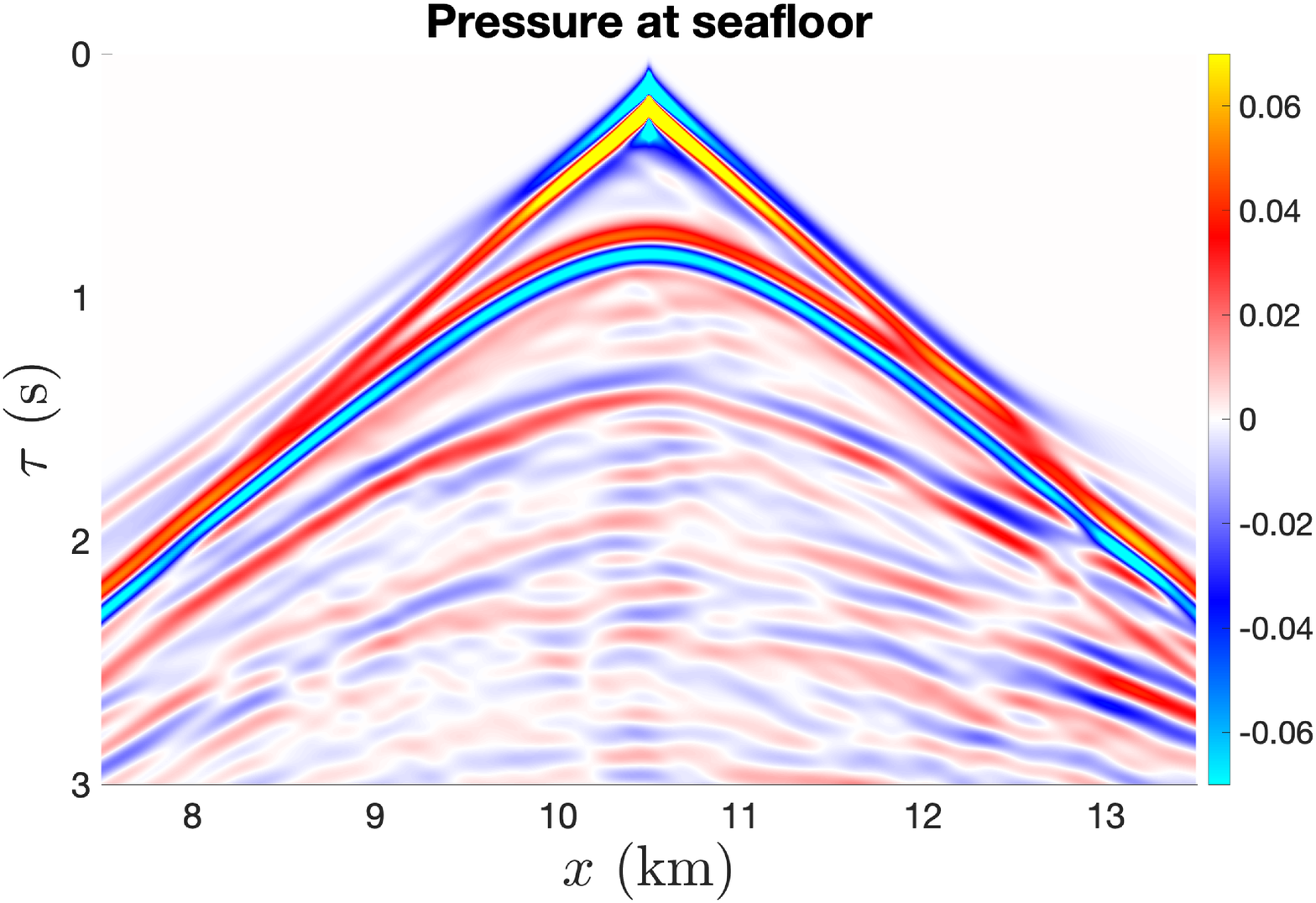}
    \end{subfigure}
    \begin{subfigure}[b]{.45\linewidth}
        \centering
        \includegraphics[width=1\linewidth,trim={1cm 0cm 1.5cm 0.5cm},clip]{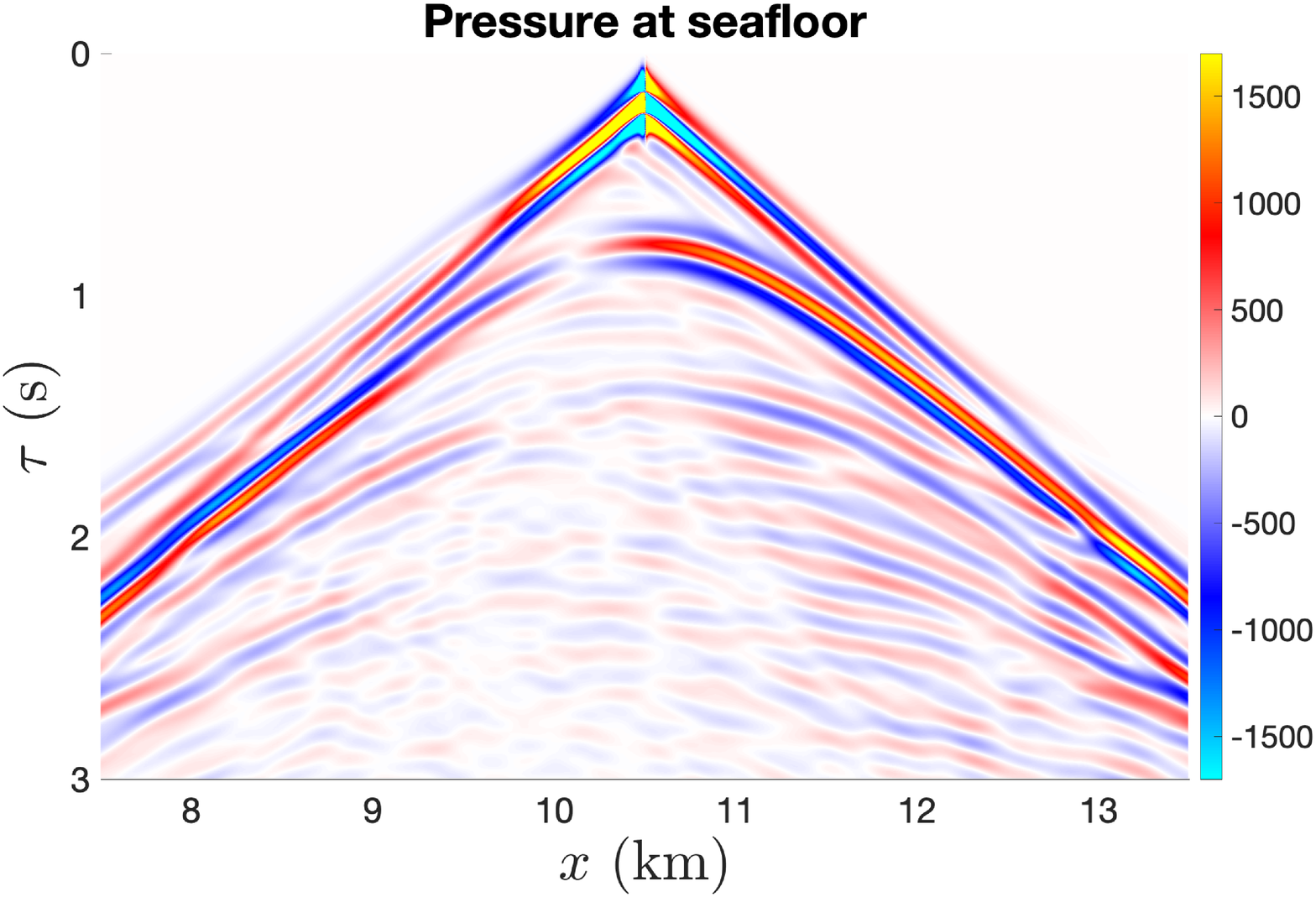}
    \end{subfigure}
    \caption{\emph{Left column:} Hydrophone data source. \emph{Right column:} Geophone data source. The top three rows show snapshots of pressure. The bottom row shows space-time plots of pressure in the water at the seafloor, i.e., at $y=y_I+$.}
    \label{fig:snapshotsMarmousiPressure}
\end{figure}

\begin{figure}[H]
\centering
    \begin{subfigure}[b]{.45\linewidth}
        \centering
        \includegraphics[width=1\linewidth,trim={0.8cm 1cm 1.5cm 1.5cm},clip]{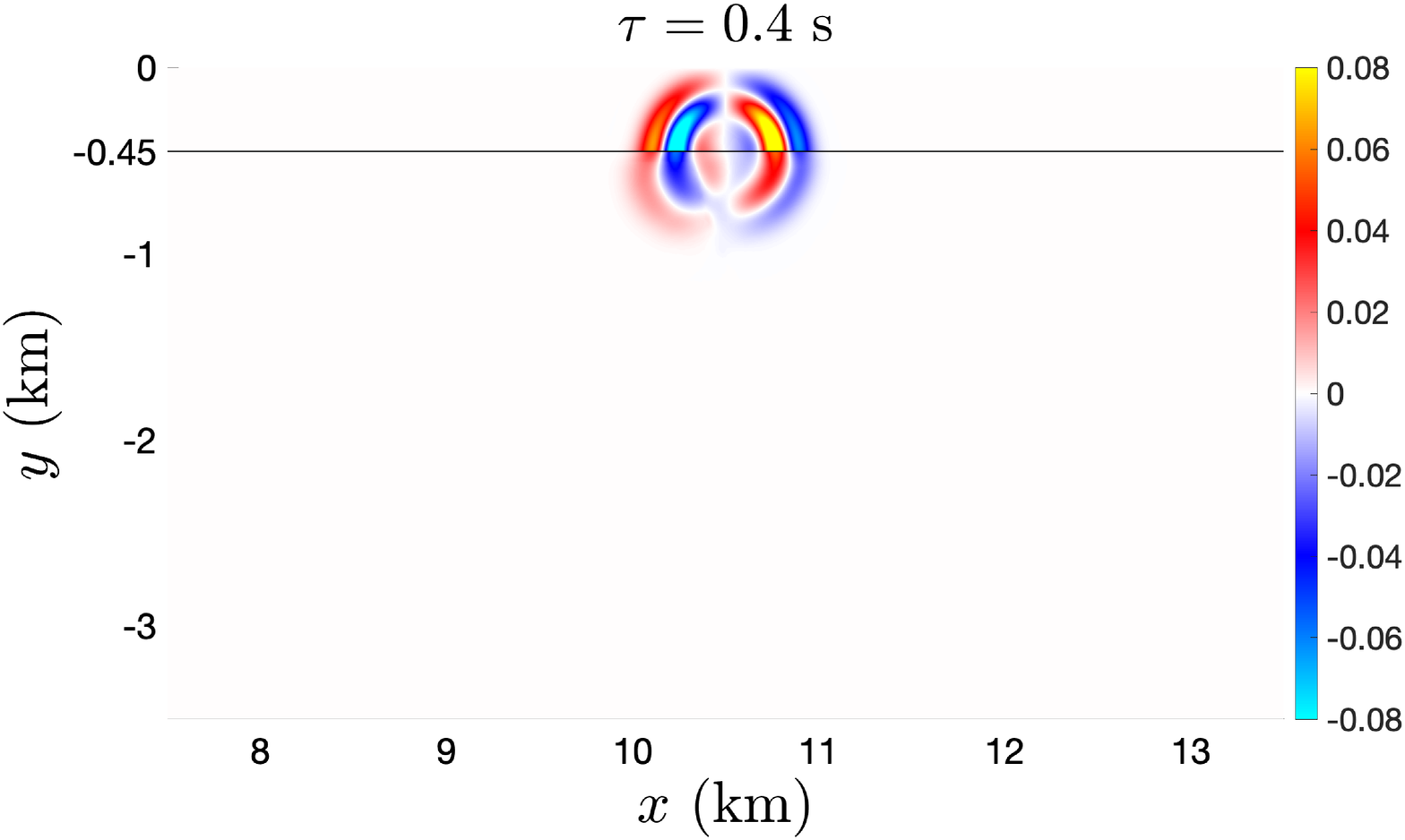}
    \end{subfigure}
     \begin{subfigure}[b]{.45\linewidth}
        \centering
        \includegraphics[width=1\linewidth,trim={0.8cm 1cm 1.5cm 1.5cm},clip]{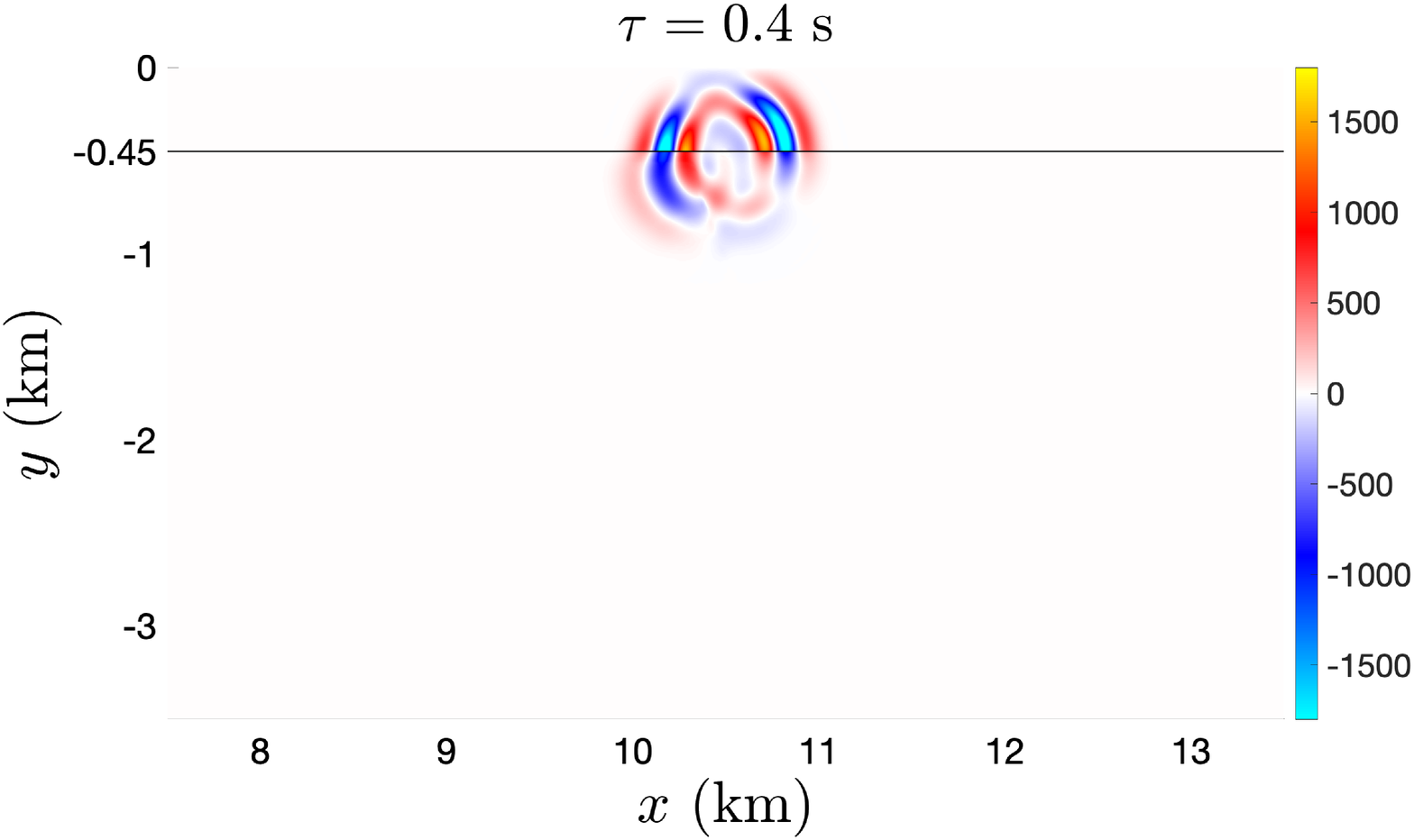}
    \end{subfigure}
    \\
    \begin{subfigure}[b]{.45\linewidth}
        \centering
        \includegraphics[width=1\linewidth,trim={0.8cm 1cm 1.5cm 1.5cm},clip]{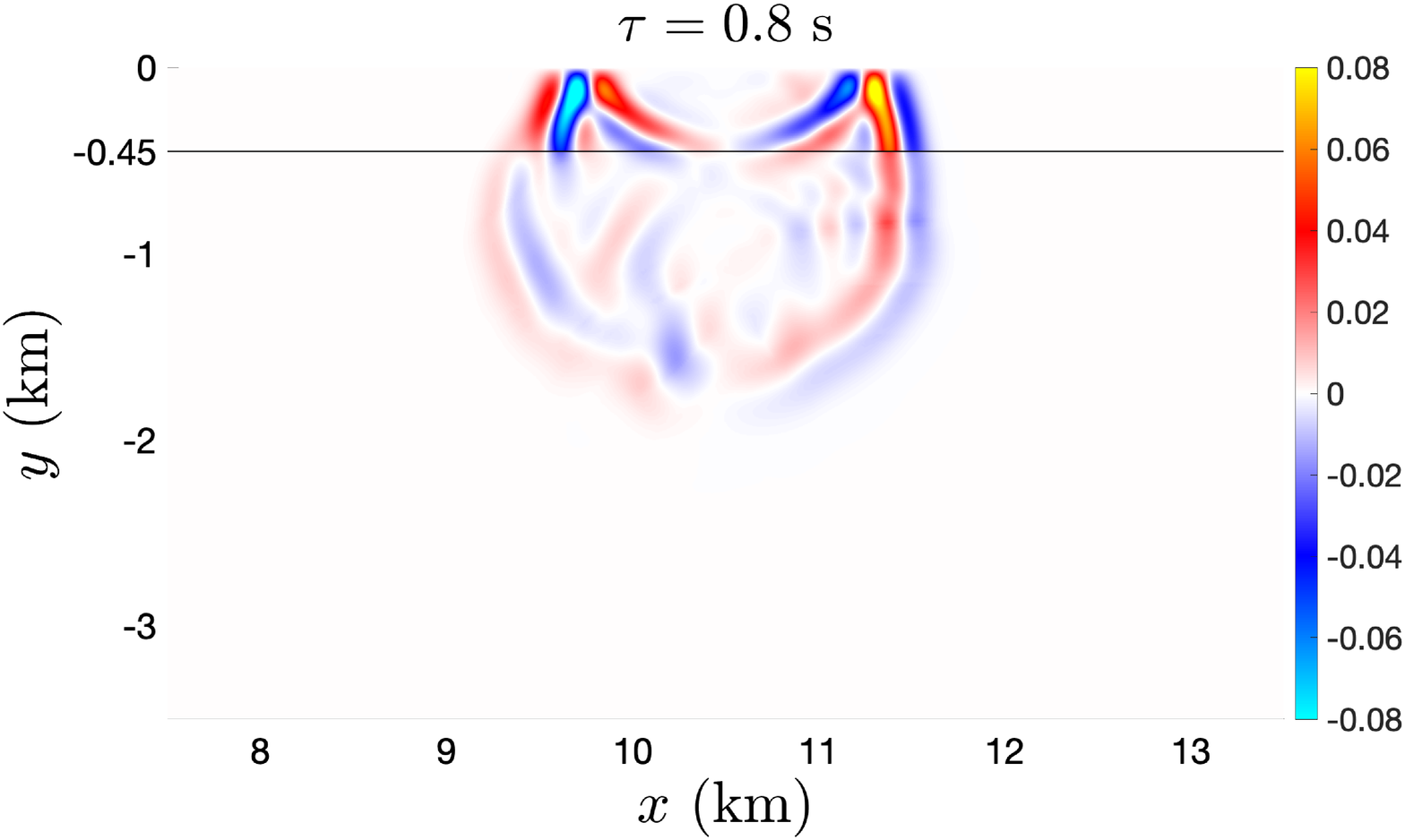}
    \end{subfigure}
    \begin{subfigure}[b]{.45\linewidth}
        \centering
        \includegraphics[width=1\linewidth,trim={0.8cm 1cm 1.5cm 1.5cm},clip]{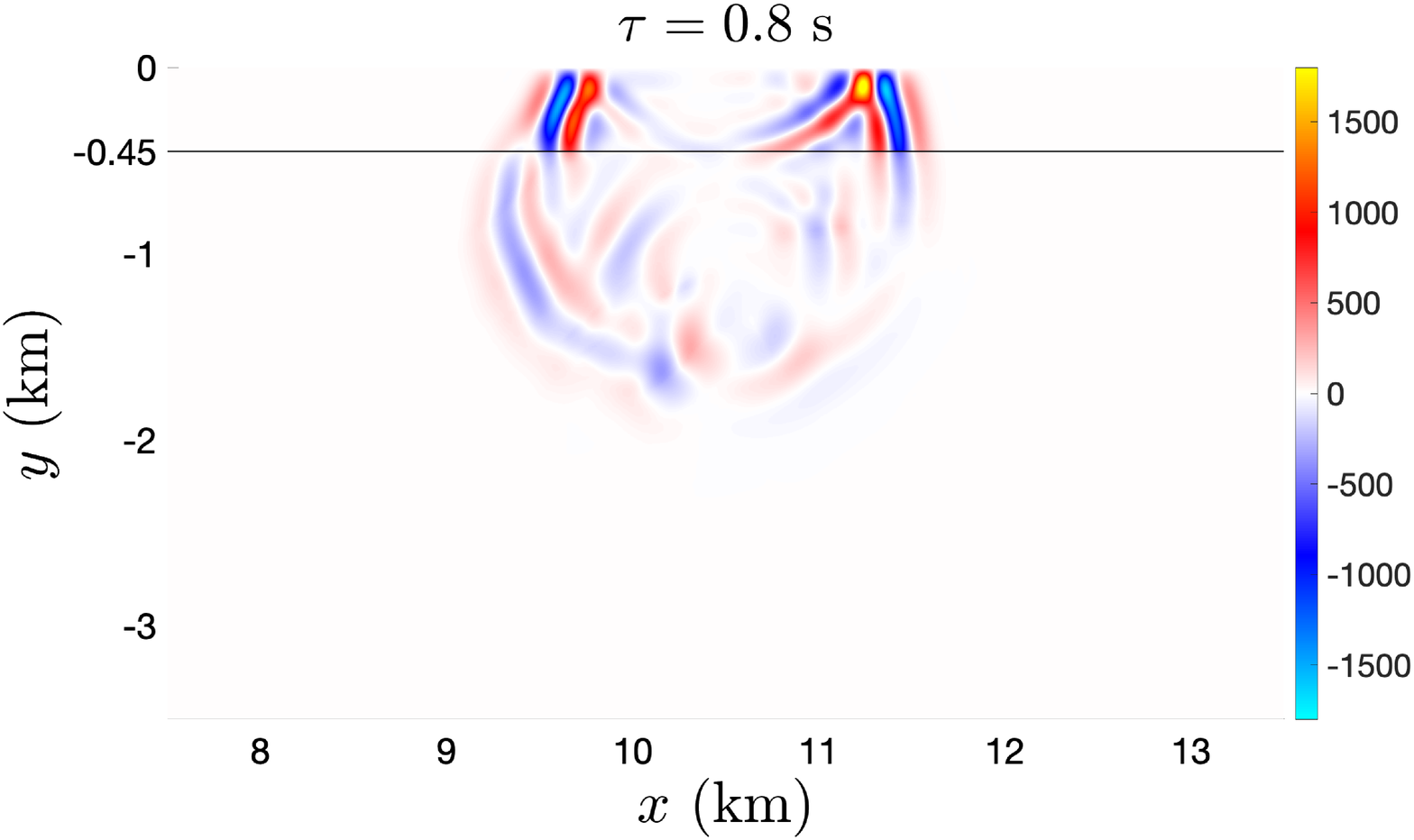}
    \end{subfigure} \\
    \begin{subfigure}[b]{.45\linewidth}
        \centering
        \includegraphics[width=1\linewidth,trim={0.8cm 1cm 1.5cm 1.5cm},clip]{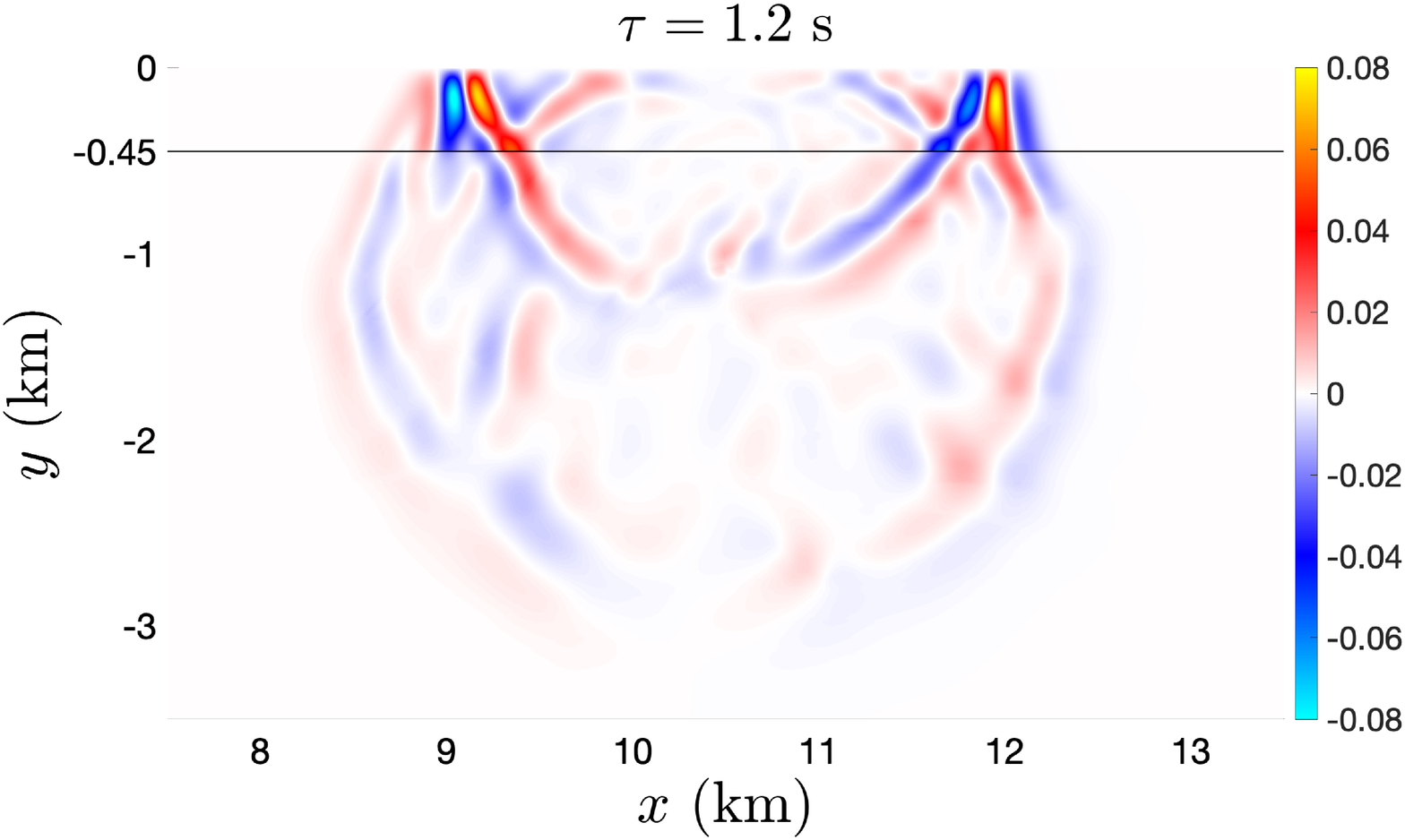}
    \end{subfigure}
    \begin{subfigure}[b]{.45\linewidth}
        \centering
        \includegraphics[width=1\linewidth,trim={0.8cm 1cm 1.5cm 1.5cm},clip]{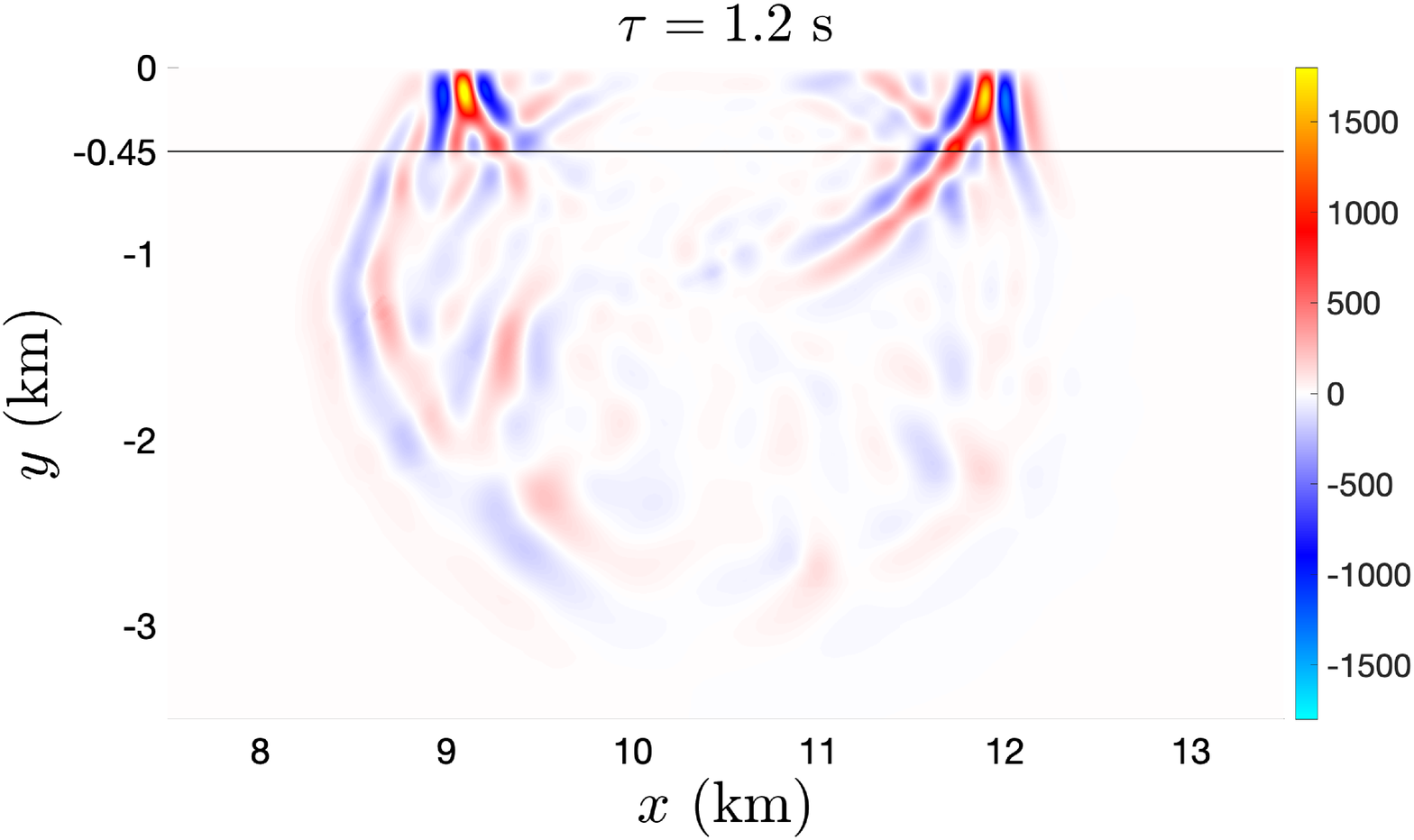}
    \end{subfigure} \\
    \begin{subfigure}[b]{.45\linewidth}
        \centering
        \includegraphics[width=1\linewidth,trim={1cm 0cm 1.5cm 0.5cm},clip]{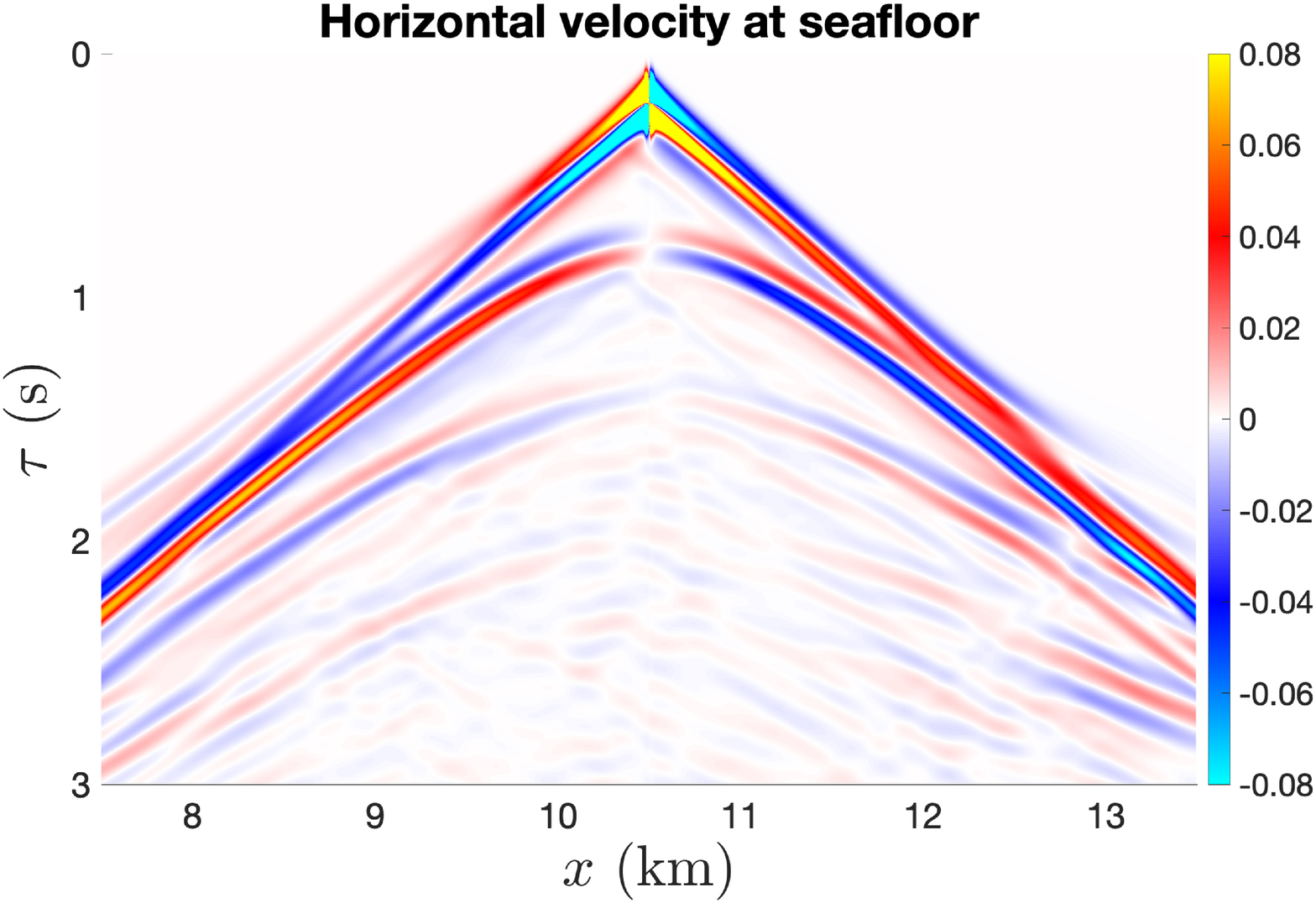}
    \end{subfigure}
    \begin{subfigure}[b]{.45\linewidth}
        \centering
        \includegraphics[width=1\linewidth,trim={1cm 0cm 1.5cm 0.5cm},clip]{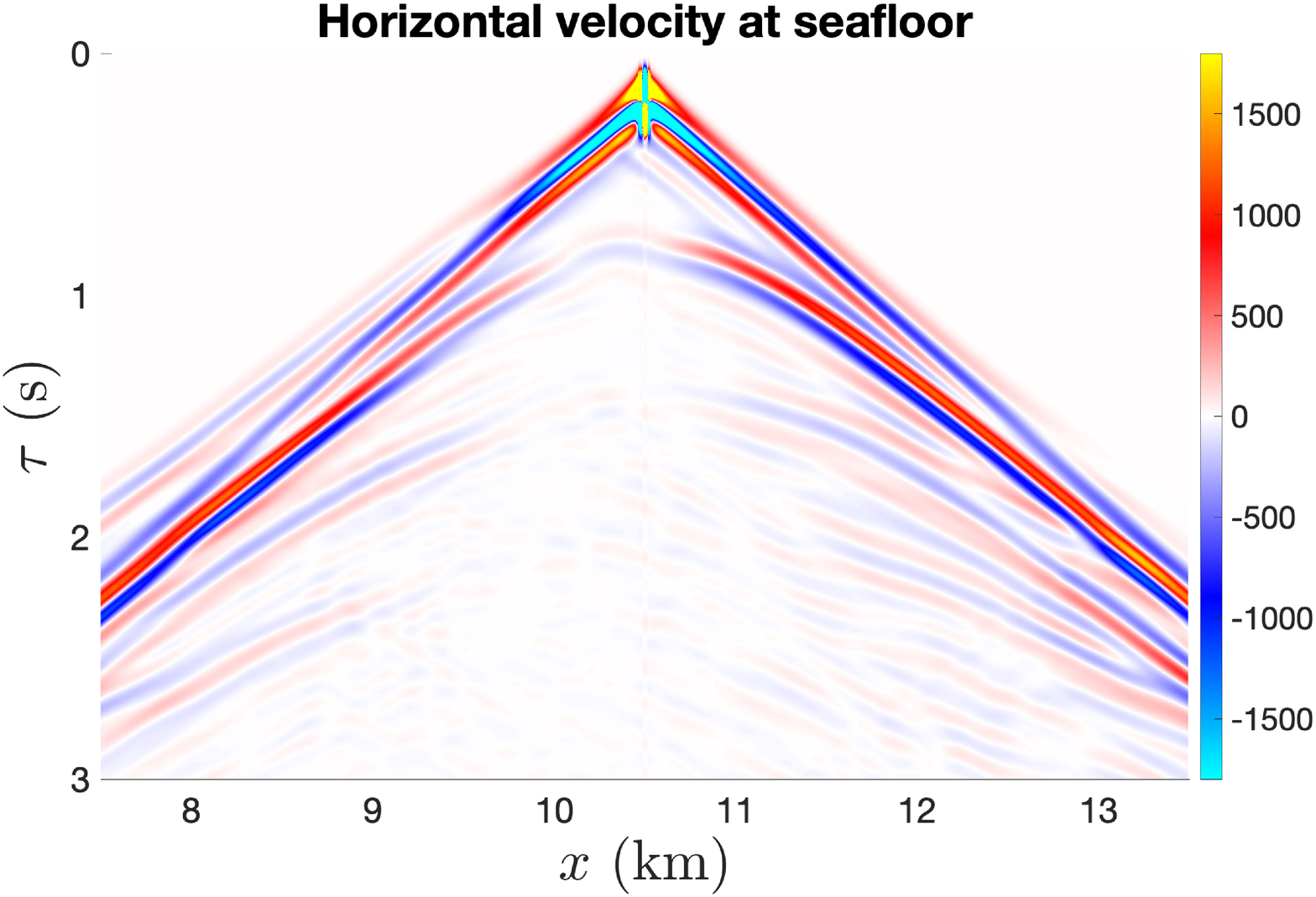}
    \end{subfigure}
    \caption{\emph{Left column:} Hydrophone data source. \emph{Right column:} Geophone data source. The top three rows show snapshots of horizontal velocity. The bottom row shows space-time plots of horizontal velocity in the water at the seafloor, i.e., at $y=y_I+$.}
    \label{fig:snapshotsMarmousiVx}
\end{figure}

 Let $c_0 = 1500$ m/s and $\rho_0=1010$ kg/m$^3$ denote reference values for sound speed and density of water. We define non-dimensional pressures and velocities as
 \begin{equation}
     p^{\dagger}_{nd} = \frac{p^{\dagger}}{\rho_0 A_p}, \quad v^{\dagger}_{i,nd}  = \frac{v^{\dagger}_i c_0}{A_p},
 \end{equation}
 in the case of only hydrophone data,
 and
 \begin{equation}
     p^{\dagger}_{nd} = \frac{c_0 p^{\dagger}}{f A_v}, \quad v^{\dagger}_{i,nd} = \frac{v^{\dagger}_i \rho_0 c_0^2}{f A_v},
 \end{equation}
 in the case of only geophone data. Figure \ref{fig:snapshotsMarmousiPressure} shows snapshots of $p^{\dagger}_{nd}$, computed with the fourth order ($q=4$) spatial discretization. The left column corresponds to having only hydrophone data ($\widetilde{v}_i = 0, \, i =1,2$), and the right column corresponds to having only geophone data ($\widetilde{p} = 0$). The bottom row shows space-time plots of pressure at the seafloor, on the ocean side ($y= y_I+$). Figure \ref{fig:snapshotsMarmousiVx} shows the horizontal component of $v^{\dagger}_{i,nd}$ in the same way.

Figures \ref{fig:snapshotsMarmousiPressure} and \ref{fig:snapshotsMarmousiVx} show that, despite placing the singular sources right on the fluid-solid interface, there are no visible artifacts from the source discretizations. Further, Figure \ref{fig:snapshotsMarmousiPressure} shows that the simulated pressure is continuous across the interface, as prescribed by \eqref{eq:interface_conditions}. Because the interface conditions are imposed weakly, the discrete pressure values on the two sides of the interface are actually only equal to the order of accuracy, but this small difference is not visible to the eye. Figure \ref{fig:snapshotsMarmousiVx}, on the other hand, shows that a discontinuity in the horizontal velocity component develops at the interface. Despite the discontinuity, there are no signs of numerical artifacts arising from the interface treatment.

\subsection{Resonant infrasound tones of Volcán Villarrica}

\begin{figure}[h]
    \centering
    \begin{subfigure}[b]{.95\linewidth}
        \centering
        \includegraphics[width= 1\linewidth]{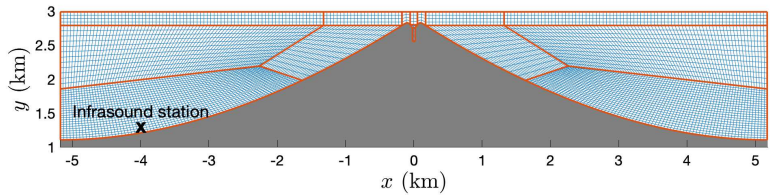}
        \caption{Coarse grid, full domain}\label{fig:grid_full}
    \end{subfigure} \\
    \begin{subfigure}[b]{.6\linewidth}
        \centering
        \includegraphics[width=1\linewidth]{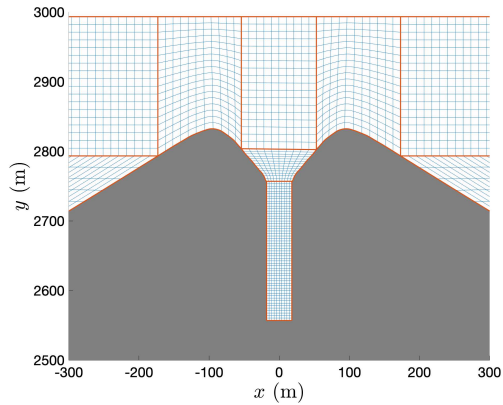}
        \caption{Finer grid, zoomed-in view of the crater region}\label{fig:grid_zoom}
    \end{subfigure}
    \caption{(a) A coarse example grid, with only $m=3$ grid points in the horizontal direction inside the conduit. (b) A zoom on the crater region, showing the $m=11$ grid. The computations used $m=26$.}
    \label{fig:grid}
\end{figure}

In \cite{Johnson2018} it was demonstrated how infrasound can be used to track the surface of the lava lake in the crater of Volcán Villarrica (Chile). The position of the lava lake may in turn be used to forecast volcanic eruptions. The nearest infrasound station is located about 4 km northwest of the crater.
To obtain a relevant crater topography we use a two-dimensional average of the three-dimensional crater geometry in \cite{Johnson2018}, see $x \in [-173 \mbox{ m}, \; 173 \mbox{ m}]$ in Figure \ref{fig:grid_zoom}. The conduit, centered at $x=0$ m, is approximated to have perfectly vertical sides. Figure \ref{fig:grid_full} shows the full computational domain. Outside of the crater region, we let the topography be described by a degree five polynomial that approximately matches the elevation drop of the Villarrica mountain in the northwest direction.
We assume constant density $\rho = 1 \, \mbox{kg/m}^3$  in the atmosphere and introuduce the velocity potential $\phi$ such that
\begin{equation}
    v_i = \ddi \phi, \quad p = -\rho \dot{\phi} .
\end{equation}
We further approximate the sound speed $c = \sqrt{K/\rho}$ as constant and set $c = 331$ m/s.
The acoustic equations \eqref{eq:acoustics_first_order} with $F_i=0$, $G=0$ are then satisfied if
\begin{equation} \label{eq:vel_pot}
     \ddot{\phi} - c^2 \ddi \ddi \phi = 0.
 \end{equation}
The Laplacian in \eqref{eq:vel_pot} is discretized using the scheme developed in this paper, and we choose fourth order accuracy ($q=4$) in the numerical experiments. We use the classical fourth order Runge--Kutta method for time-integration. The far-field boundaries are again treated with super-grid-scale absorbing layers \cite{Appelo09_2,petersson_sjogreen_2014}. The conduit walls and the mountain sides are modeled as rigid walls, resulting in the boundary condition $\vec{v} \cdot \hat{n} = 0$. At the conduit bottom we simulate movement of the lava lake by imposing the boundary condition
\begin{equation}
    \vec{v} \cdot \hat{n} = -v_0 e^{-\frac{(t-t_0)^2}{2\sigma^2}},
\end{equation}
with $t_0=1$ s and $\sigma=0.13$ s. The same boundary data were used in \cite{Johnson2018}, with the motivation that they include frequencies that are characteristic of lava lake events.

\begin{figure}[h]
    \begin{subfigure}[b]{.5\linewidth}
        \centering
        \includegraphics[width= 1\linewidth]{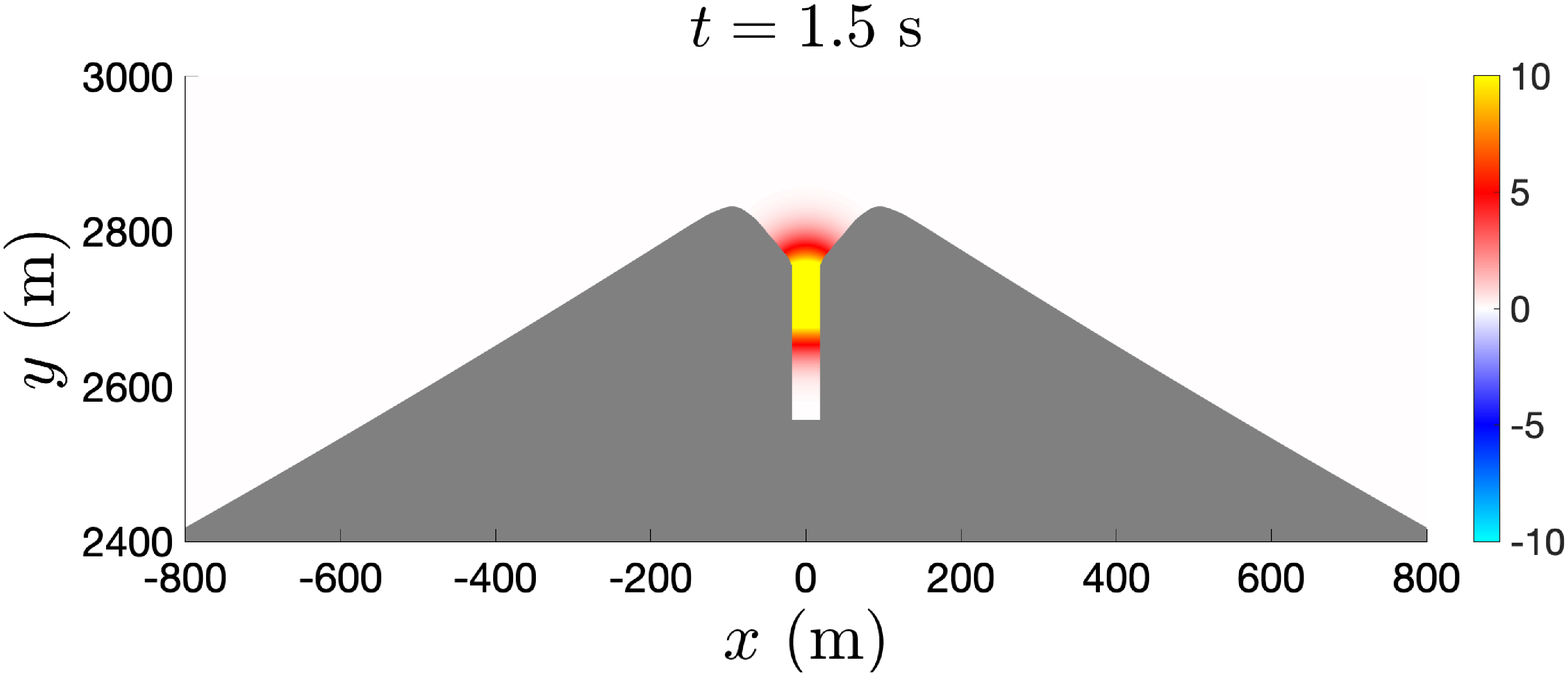}
    \end{subfigure}
    \begin{subfigure}[b]{.5\linewidth}
        \centering
        \includegraphics[width=1\linewidth]{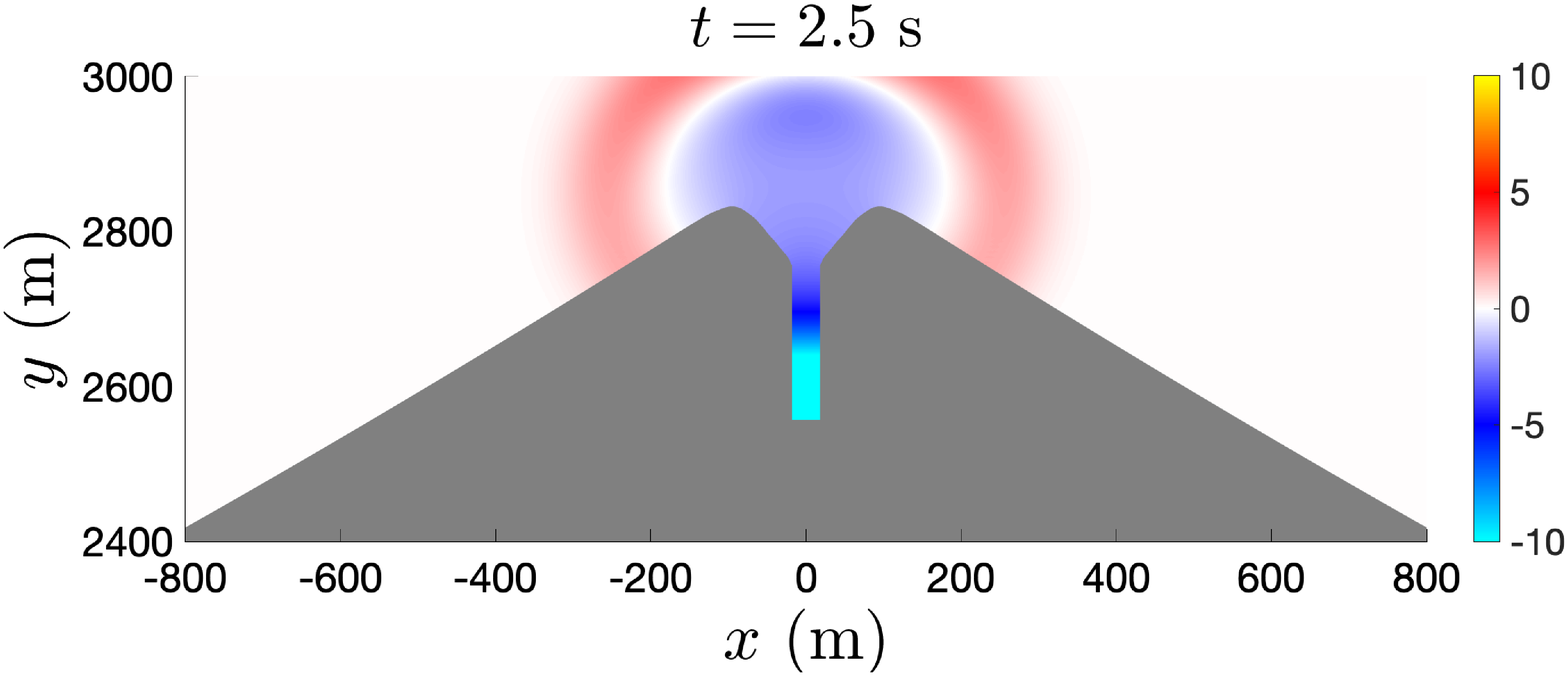}
    \end{subfigure} \\
    \begin{subfigure}[b]{.5\linewidth}
        \centering
        \includegraphics[width= 1\linewidth]{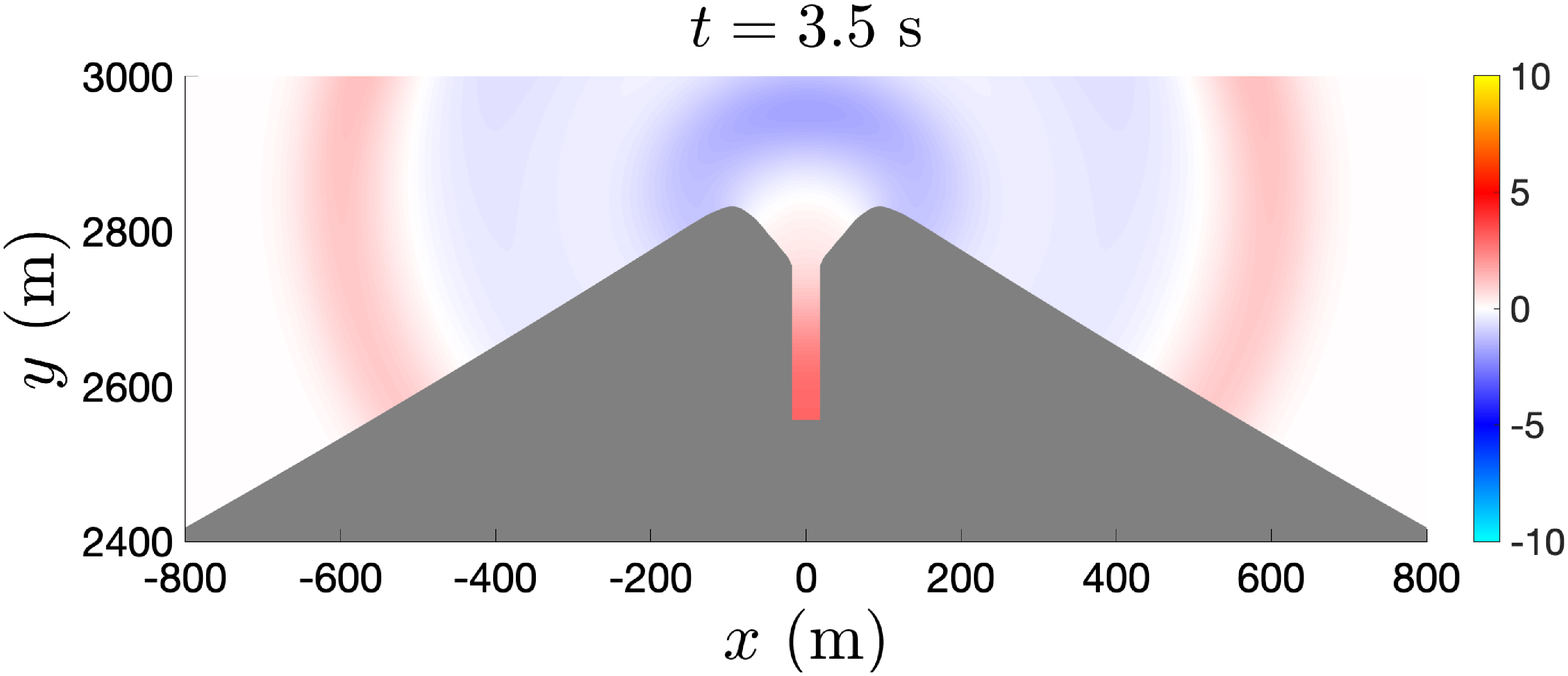}
    \end{subfigure}
    \begin{subfigure}[b]{.5\linewidth}
        \centering
        \includegraphics[width=1\linewidth]{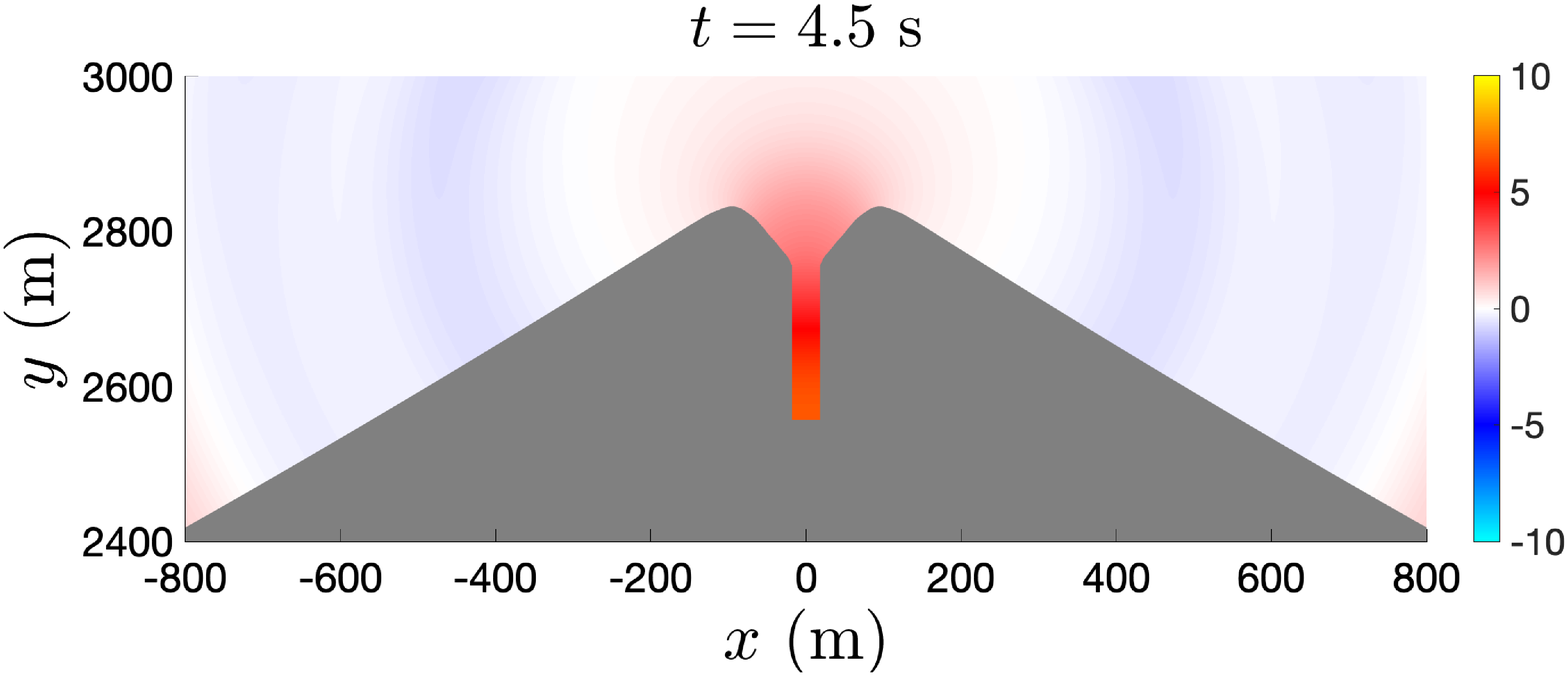}
    \end{subfigure}
    \caption{Snapshots of non-dimensional pressure $p_{nd}$}
    \label{fig:snapshotsVillarrica}
\end{figure}

\begin{figure}[h]
        \centering
        \includegraphics[width= 0.8 \linewidth]{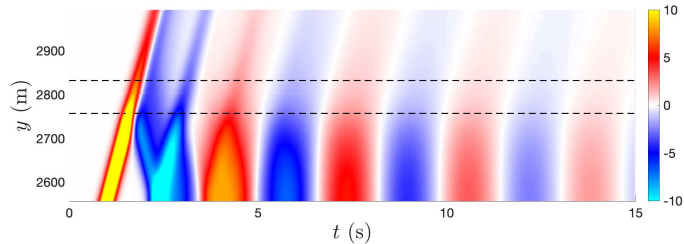}
        \caption{Space-time plot of pressure at horizontal position $x=0$ m, i.e., within and immediately above the conduit. The dashed black lines mark the top of the conduit (lower line) and the crater rim (upper line).}
    \label{fig:spacetime}
\end{figure}

\begin{figure}[h]
    \begin{subfigure}[b]{.5\linewidth}
        \centering
        \includegraphics[width= 1\linewidth]{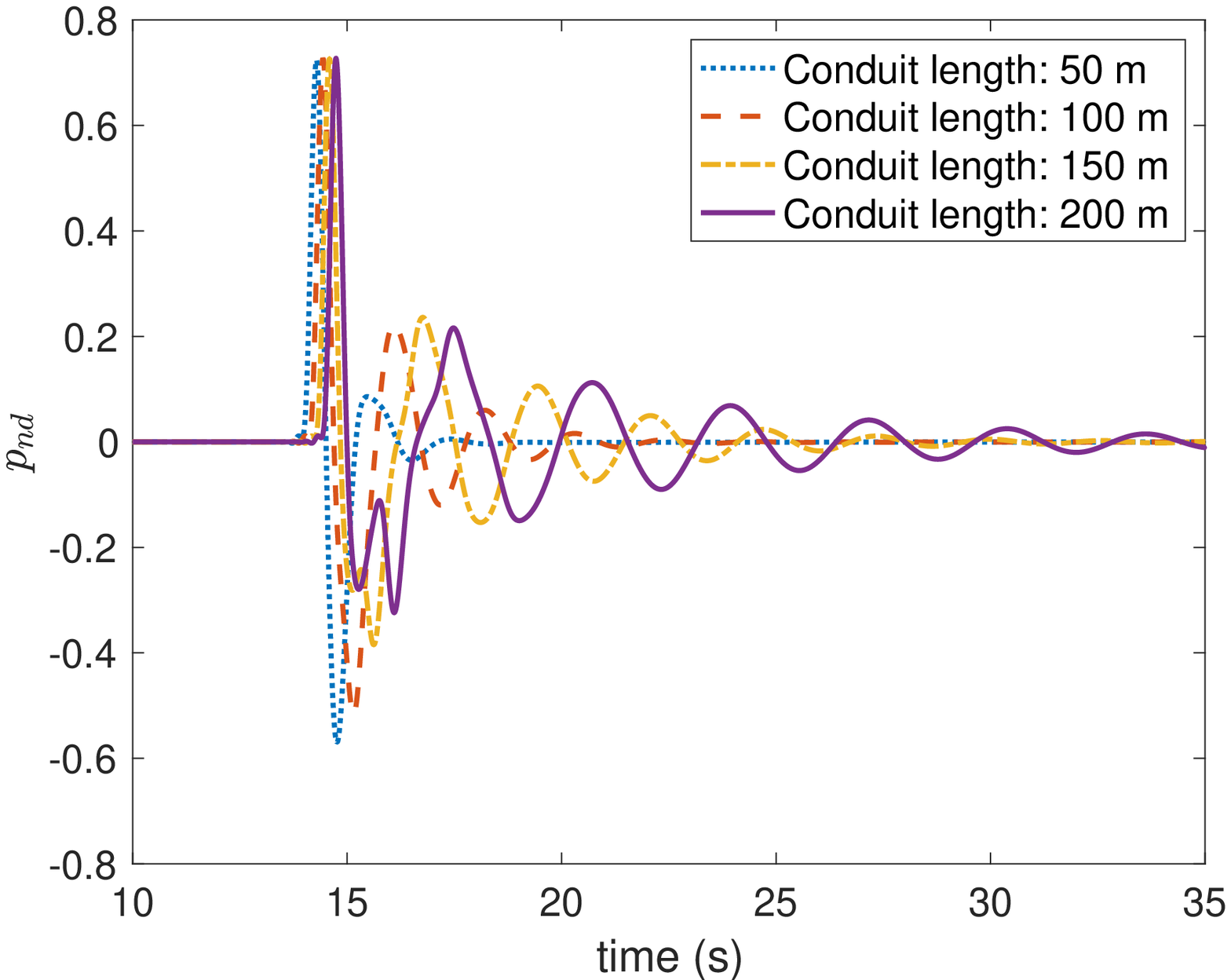}
        \caption{Pressure time series}\label{fig:timeSeries}
    \end{subfigure}
    \begin{subfigure}[b]{.5\linewidth}
        \centering
        \includegraphics[width=1\linewidth]{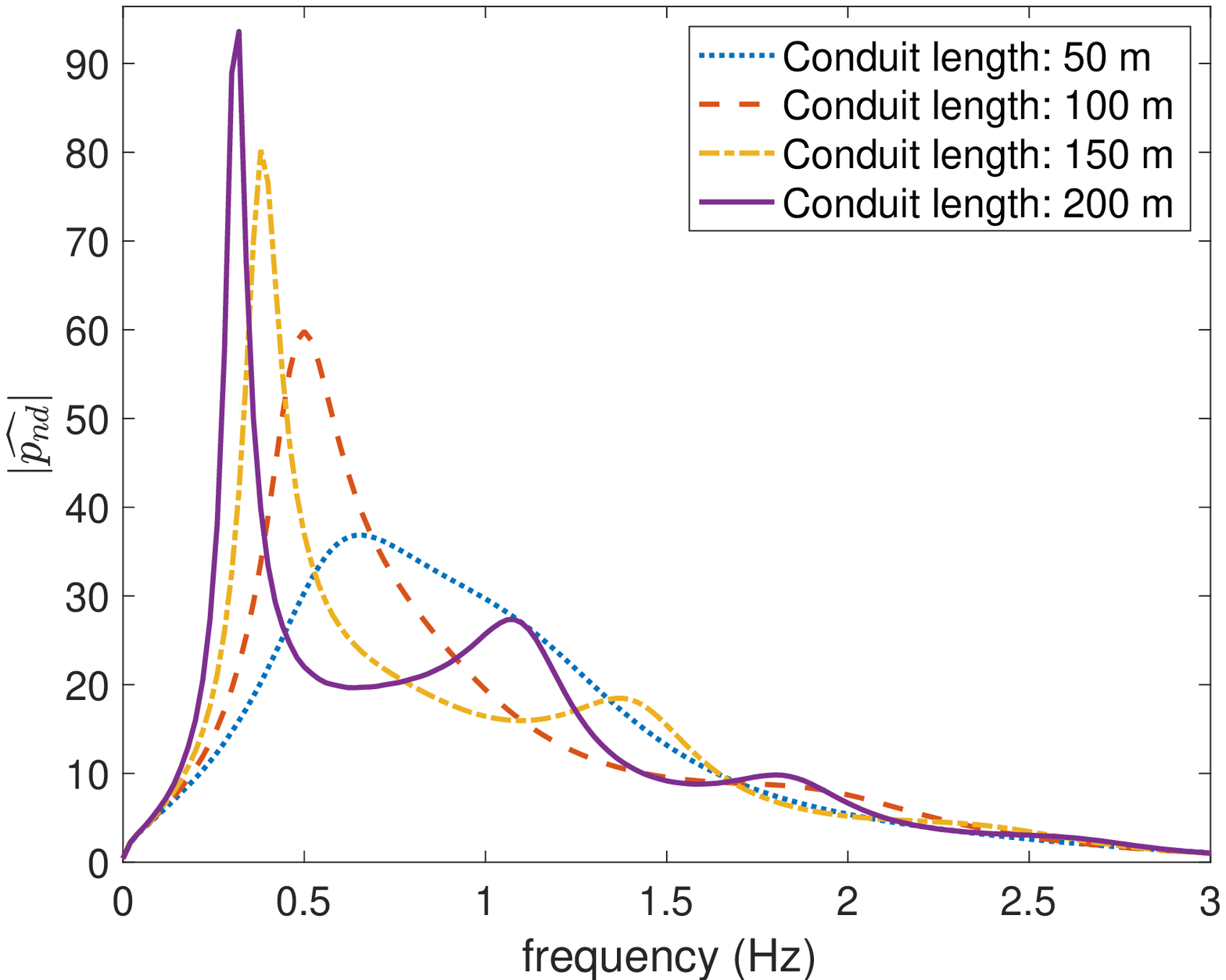}
        \caption{Pressure frequency spectrum}\label{fig:fourier}
    \end{subfigure}
    \caption{(a) Pressure recorded at the nearest infrasound station, 4 km from the crater. (b) Fourier amplitude spectrum of the time series in (a).}
    \label{fig:infrasound}
\end{figure}

We define a non-dimensional pressure as
\begin{equation}
    p_{nd} = \frac{p}{\rho c v_0} .
\end{equation}
Figure \ref{fig:snapshotsVillarrica} shows snapshots of $p_{nd}$.
Figure \ref{fig:spacetime} shows space-time plots of $p_{nd}$ at position $x=0$, i.e., within and above the conduit. One can observe how the initial pressure pulse, propagating upwards, is partly reflected at the top of the conduit.
Figure \ref{fig:timeSeries} shows the recorded pressure at the location corresponding to the nearst infrasound station, for different conduit lengths. Figure \ref{fig:fourier} shows the corresponding frequency spectra. The longer the conduit, the lower the crater resonance frequency and the higher the quality factor.

This simulation illustrates the benefits of the curvilinear multiblock finite difference discretization. The method takes topography into account without using staircasing, which might introduce artifacts. Further, it imposes the solid wall boundary condition accurately and stably.

\section{Conclusions}
We have improved the SBP-SAT method for discretizing the Laplacian on curvilinear multiblock grids in $d$ dimensions. Compared to the previous state-of-the-art method by Virta and Mattsson \cite{VirtaMattsson14}, the new method is significantly less stiff, in particular when the grid is highly skewed. On a simple curvilinear multiblock grid for the unit disk, the largest stable time-step when solving the wave equation with explicit time-stepping is approximately 6 times larger with the new method. Although the method by Virta and Mattsson is typically slightly more accurate on a given grid, our numerical experiments show that the new method is significantly more efficient in terms of CPU time required to achieve a given error tolerance. Further, our results indicate that for more complex grids, which will inevitably contain highly skewed grid cells, the time-step ratio will be even larger in favor of the new method.

The difference compared to the method by Virta and Mattsson lies solely in the SATs for Dirichlet boundary conditions and interface couplings. In essence, the new method is less stiff because the SATs are smaller in magnitude. The reason we were able use smaller SATs and still prove stability is the novel positivity property in Lemma \ref{lemma:R}, which was not used by Virta and Mattsson.

To illustrate the usefulness of the new method, we have applied it to two problems inspired by marine seismic exploration and infrasound monitoring of volcanoes. The seismic exploration problem uses a Cartesian two-block grid with a discontinuity in material parameters at the ocean floor, where the ocean couples to the seabed. The numerical method handles both the fluid-solid coupling and singular source terms (both $\delta$ functions and derivatives of $\delta$ functions), placed right at the fluid-solid interface, without signs of numerical artifacts. The volcano monitoring problem involves non-trivial topography that requires a curvilinear multiblock grid. Again, the discrete solution shows no sign of numerical artifacts.

MATLAB code that reproduces all figures in this paper is available at \\
\mbox{\url{https://bitbucket.org/martinalmquist/laplacian_curvilinear}}.

\section*{Acknowledgements}
We thank Leighton Watson for providing the Villarrica crater topography.
M.\ Almquist gratefully acknowledges support from the Knut and Alice Wallenberg Foundation (Dnr.\ KAW 2016.0498).

\renewcommand{\ddi}{\partial_{i}}
\section*{Appendix: proof of Lemma \ref{lemma:R}}
By construction, $R_{xx}(a)$ is linear in its argument so that, for all functions $a, b$ and scalars $C$,
\begin{equation}
\begin{aligned}
    R_{xx}(a+b) &= R_{xx}(a) + R_{xx}(b), \\
    R_{xx}(Ca) &= C R_{xx}(a).
\end{aligned}
\end{equation}
Further, $R_{xx}$ is symmetric positive semidefinite in the sense that, for all $a(x) \geq 0$,
\begin{equation}
    R_{xx}(a) = R_{xx}^T(a) \geq 0.
\end{equation}
For detailed information about the structure of $R_{xx}$ we refer to \cite{Mattsson11}. It follows from linearity and positive semidefiniteness that if $a(x) \geq b(x)$, then
\begin{equation}
    R_{xx}(a) = R_{xx}(b + a - b) = R_{xx}(b) + R_{xx}(a-b) \geq R_{xx}(b).
\end{equation}
Given grid endpoints $x_{\ell}$ and $x_{r}$, we define cutoff points at the $m_b$th grid point from the boundaries:
\begin{equation}
    x_{\ell, m_b} = x_{\ell}+(m_b-1)h, \quad x_{r, m_b} = x_{r} - (m_b-1)h .
\end{equation}
Let $b_{\ell,min}$ and $b_{r,min}$ denote the minimum of $b(x)$ over the $m_b$ leftmost and rightmost grid points,
\begin{equation}
\begin{aligned}
    b_{\ell,min} &= \min \left( b(x_{\ell}), b(x_{\ell}+h), \ldots, b(x_{\ell,m_b}) \right), \\
    b_{r,min} &= \min \left( b(x_{r,m_b}), b(x_{r,m_b}+h), \ldots, b(x_{r}) \right).
\end{aligned}
\end{equation}
Let $\chi_{\ell,r}$ denote the indicator functions
\begin{equation}
    \chi_{\ell}(x) =
\left\{ \begin{array}{cc}
1, & x \leq x_{\ell,m_b} \\
0, & x > x_{\ell,m_b} \\
\end{array}
\right.,
\quad \chi_{r}(x) =
\left\{ \begin{array}{cc}
0, & x < x_{r,m_b} \\
1, & x \geq x_{r,m_b} \\
\end{array}
\right. .
\end{equation}
Assuming that there are at least $2m_b$ grid points, then $x_{\ell,m_b} < x_{r,m_b}$ and hence, for all nonnegative functions $b$,
\begin{equation}
\begin{aligned}
    b(x) &\geq b(x) \left( \chi_{\ell}(x) + \chi_{r}(x) \right) = b(x) \chi_{\ell}(x) + b(x) \chi_{r}(x) \\
    &\geq b_{\ell,min} \chi_{\ell}(x) + b_{r,min} \chi_{r}(x).
\end{aligned}
\end{equation}
We thus have
\begin{equation}
    R_{xx}(b) \geq b_{\ell,min} R_{xx}\left(\chi_{\ell} \right) + b_{r,min} R_{xx} \left( \chi_{r} \right).
\end{equation}
To determine the constant $\theta_R$ in the property
\begin{equation}
    \bfu^T R_{xx}(b) \bfu \geq h\theta_R b_{\ell,min} (e_{\ell}^T \Delta D_x \bv{u})^2 + h\theta_R b_{r,min} (e_r^T \Delta D_x \bv{u})^2,
\end{equation}
it is thus sufficient to consider $\bfu^T R_{xx}(\chi_{\ell}) \bfu$ or, equivalently, $\bfu^T R_{xx}(\chi_{r}) \bfu$.

\begin{figure}[h]
        \centering
        \includegraphics[width= 0.6 \linewidth]{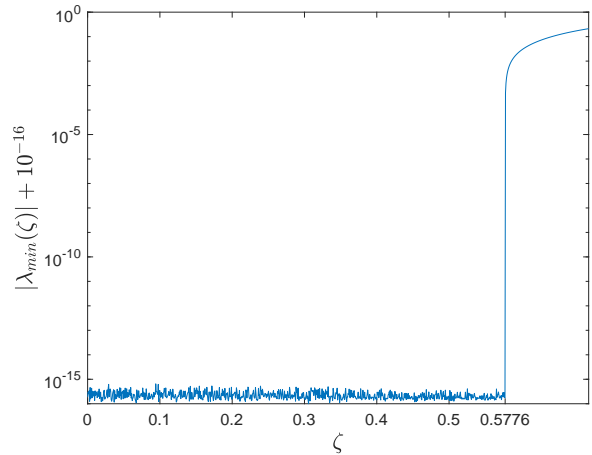}
        \caption{Absolute value of the most negative eigenvalues of $N(\zeta)$, for fourth order accuracy with $m_b=4$. The most negative eigenvalue becomes significantly different from zero at $\zeta = 0.5776$. }
    \label{fig:borrowing}
\end{figure}

We have determined $\theta_R$ by computing eigenvalues of
\begin{equation}
N(\zeta) = R_{xx}(\chi_{\ell}) - \zeta h (e_{\ell}^T \Delta D_x)^T (e_{\ell}^T \Delta D_x)
\end{equation}
for different values of the scalar $\zeta$. Since $R_{xx}(\chi_{\ell})$ is only semidefinite, it has one or more zero eigenvalues. We define $\theta_R$ as the smallest value such that $N(\theta_R)$ has negative eigenvalues. In practice, we compute the smallest value such that the most negative eigenvalue of $N(\theta_R)$ is not within roundoff error from zero. Let $\lambda_{min}(\zeta)$ denote the smallest eigenvalue of $N(\zeta)$. Figure \ref{fig:borrowing} shows $\abs{\lambda_{min}(\zeta)} + 10^{-16}$ on a logarithmic scale, for the fourth order accurate case. The grid spacing is $h=1$ so that roundoff errors in the eigenvalue computations are expected to be of order $10^{-16}$. It is evident that $\lambda_{min}(\zeta)$ is nonzero for $\zeta > 0.5776 = \theta_R$.

We remark that $\theta_R$ is independent of $m$, as long as the total number of grid points $m$ is at least $2 m_b$. However, the value of $\theta_R$ increases with $m_b$ and approaches an asymptotic value, as shown in Table \ref{table:theta_R}. The drawback to having large $m_b$ is that it requires computing the minimum of the coefficient $b$ over many grid points, which can potentially make $b_{\ell,min}$ and $b_{r,min}$ smaller. To minimize the penalty strength required for the stability proof, we would like to maximize the products $b_{\ell,min}\theta_R$ and $b_{r,min}\theta_R$. We thus choose the smallest value of $m_b$ that yields $\theta_R$ near the asymptotic value. Based on Table \ref{table:theta_R} we choose $m_b = 2,4,7$, for orders two, four and six. The same $m_b$ values were chosen in \cite{VirtaMattsson14} when computing similar constants for the matrix $M_{xx}(b) = D_x^T H b D_x + R_{xx}(b)$.

\begin{table}[H]
\centering
\begin{tabular}{cccc}
& second order & fourth order & sixth order \\
\hline
$m_b=1$ & 0 & 0 & 0 \\
$m_b=2$ & \textbf{1.0000} & 0 & 0 \\
$m_b=3$ & 1.0000 & 0.1485 & 0 \\
$m_b=4$ & 1.0000 & \textbf{0.5776} & 0 \\
$m_b=5$ & 1.0000 & 0.5779 & 0 \\
$m_b=6$ & 1.0000 & 0.5779 & 0.2318 \\
$m_b=7$ & 1.0000 & 0.5779 & \textbf{0.3697} \\
$m_b=8$ & 1.0000 & 0.5779 & 0.3697 \\
$m_b=9$ & 1.0000 & 0.5779 & 0.3697 \\
\end{tabular}
\caption{Computed values of $\theta_R$ for different values of $m_b$. The chosen $(m_b,\theta_R)$ pairs are indicated with boldface font. }
\label{table:theta_R}
\end{table}


\end{document}